%% file: main.tex
\documentclass[onefignum,onetabnum]{siamonline220329}

\input{shared}

\input{0_commands.tex}

\begin{document}

\maketitle

\begin{abstract}
  Heteroclinic cycles and networks are structures in dynamical systems composed of invariant sets and connecting heteroclinic orbits, and can be robust in systems with invariant subspaces. The usual method for analysing the stability of heteroclinic cycles and networks is to construct return maps to cross-sections near the network. From these return maps, transition matrices can be defined, whose eigenvalues and eigenvectors can be analysed to determine stability. In this paper, we introduce an extension to this methodology, the \textit{projected map}, which we define by identifying trajectories with, in a certain sense, qualitatively the same dynamics. The projected map is a discrete, piecewise-smooth map of one dimension one fewer than the rank of the transition matrix. We use these maps to describe the dynamics of trajectories near three heteroclinic networks in \(\R^{4}\) with four equilibria. We find in all three cases that the onset of trajectories that switch between cycles of the network can be caused by a fold bifurcation or a \textit{border-collision bifurcation}, where fixed points of the map no longer exist in the corresponding function's domain of definition. We are able to show that trajectories near these three networks are asymptotic to only one subcycle, and cannot switch between subcycles multiple times, in contrast to examples with more than four equilibria, and resolving a 30-year-old claim by Brannath. We are also able to generalise some results to all quasi-simple networks, showing that a border-collision bifurcation of the projected map, corresponding to a condition on the eigenvectors of certain transition matrices, causes a cycle to lose stability.
\end{abstract}



\input{1_introduction.tex}
\input{2_background.tex}

\input{3_projected_map.tex}

\input{4_analysis.tex}

\input{5_generalised.tex}

\input{6_discussion.tex}

\input{7_appendix.tex}

\section*{Acknowledgments}

The authors thank David Simpson for discussions that led to the idea of the projected map. We also thank Sofia Castro and Alexandre Rodrigues for helpful conversations. This research was funded by the Marsden Fund Council from New Zealand Government funding, managed by the Royal Society Te Ap\=arangi (grant no. 21-UOA-048).

\bibliographystyle{siamplain}
\bibliography{main}

\end{document}

\maketitle

\section{Maps near heteroclinic networks}\label{supp:maps}

We list all basic transition matrices and relevant full transition matrices used in the analysis of the three networks we study. Recall that the basic transition matrices have the notation \(m_{ijk}=M(\varphi_{ijk})\), and the full transition matrices have the notation \(M_{j}^{(k)}=M(\Phi_{j}^{(k)})\). The basic maps \(\varphi_{ijk}\colon\poinin[i]{j}\to\poinin[j]{k}\) are the maps between successive incoming cross-sections in a cycle, and the full return maps \(\Phi_{j}^{(k)}\colon\poinin[i]{j}\to\poinin[i]{j}\) are the maps from the incoming cross-section of an equilibrium \(\xi_{j}\) to itself, around the \(\cyc{[k]}\) cycle. The full return maps are the composition of basic maps, and so the full transition matrices are the product of basic transition matrices.

\subsection{The Kirk--Silber network}\label{supp:maps:ks}

The basic transition matrices of the Kirk--Silber network are
\begin{align*}
  m_{312}&=\begin{pmatrix}
    \frac{c_{13}}{e_{12}} & 0 \\[0.5em]
    \frac{c_{14}}{e_{12}} & 1
  \end{pmatrix},\ m_{123}=\begin{pmatrix}
    \frac{c_{21}}{e_{23}} & 0 \\[0.5em]
    -\frac{e_{24}}{e_{23}} & 1
  \end{pmatrix},\ m_{231}=\begin{pmatrix}
    \frac{c_{32}}{e_{31}} & 0 \\[0.5em]
    \frac{t_{34}}{e_{31}} & 1
  \end{pmatrix},\\[0.5em]
  m_{412}&=\begin{pmatrix}
    \frac{c_{13}}{e_{12}} & 1 \\[0.5em]
    \frac{c_{14}}{e_{12}} & 0
  \end{pmatrix},\ m_{124}=\begin{pmatrix}
    0 & \frac{c_{21}}{e_{24}} \\[0.5em]
    1 & -\frac{e_{23}}{e_{24}}
  \end{pmatrix},\ m_{241}=\begin{pmatrix}
    \frac{c_{42}}{e_{41}} & 0 \\[0.5em]
    \frac{t_{43}}{e_{41}} & 1
  \end{pmatrix}.
\end{align*}
The relevant full transition matrices of the Kirk--Silber network are
\begin{align*}
  M_{3}\equiv M_{2}^{(3)}&=m_{312}m_{231}m_{123}=\begin{pmatrix}
    \delta_{3} & 0\\
    \rho_{3} & 1
  \end{pmatrix},\\
  M_{4}\equiv M_{2}^{(4)}&=m_{412}m_{241}m_{124}=\begin{pmatrix}
    1 & \rho_{4}\\
    0 & \delta_{4}
  \end{pmatrix},
\end{align*}
where
\begin{alignat*}{2}
  \delta_{3}=&\frac{c_{13}c_{21}c_{32}}{e_{12}e_{23}e_{31}}>0,&&\delta_{4}=\frac{c_{14}c_{21}c_{42}}{e_{12}e_{24}e_{41}}>0,\\
  \rho_{3}=&-\frac{e_{24}}{e_{23}}+\frac{c_{21}t_{34}}{e_{23}e_{31}}+\frac{c_{14}c_{21}c_{32}}{e_{12}e_{23}e_{31}},\qquad&&\rho_{4}=-\frac{e_{23}}{e_{24}}+\frac{c_{21}t_{43}}{e_{24}e_{41}}+\frac{c_{13}c_{21}c_{42}}{e_{12}e_{24}e_{41}}.
\end{alignat*}
From the full transition matrices \(M_{3}^{(3)}=m_{123}m_{312}m_{231}\) and \(M_{4}^{(4)}=m_{124}m_{412}m_{241}\), we derive
\begin{align*}
  \nu_{3}=&\frac{t_{34}}{e_{31}}+\frac{c_{14}c_{32}}{e_{12}e_{31}}-\frac{c_{13}c_{32}e_{24}}{e_{12}e_{23}e_{31}}\\ 
  \nu_{4}=&\frac{t_{43}}{e_{41}}+\frac{c_{13}c_{42}}{e_{12}e_{41}}-\frac{c_{14}c_{42}e_{23}}{e_{12}e_{24}e_{41}}.
\end{align*}
The remaining full transition matrices are the products \(M_{1}^{(3)}=m_{231}m_{123}m_{312}\) and \(M_{1}^{(4)}=m_{241}m_{124}m_{412}\)

\subsection{The \texorpdfstring{\(\Delta\)}{Delta}-clique network}\label{supp:maps:delta}

The basic transition matrices of the \(\Delta\)-clique network are
\begin{align*}
  m_{412}&=\begin{pmatrix}
    \frac{t_{13}}{e_{12}} & 1\\[0.5em]
    \frac{c_{14}}{e_{12}} & 0\\[0.5em]
  \end{pmatrix},\ 
  m_{123}=\begin{pmatrix}
    \frac{c_{21}}{e_{23}} & 0\\[0.5em]
    -\frac{e_{24}}{e_{23}} & 1\\[0.5em]
  \end{pmatrix},\ 
  m_{124}=\begin{pmatrix}
    0 & -\frac{e_{23}}{e_{24}}\\[0.5em]
    1 & \frac{c_{21}}{e_{24}}\\[0.5em]
  \end{pmatrix},\\[0.5em]
  m_{234}&=\begin{pmatrix}
    1 & \frac{t_{31}}{e_{34}}\\[0.5em]
    0 & \frac{c_{32}}{e_{34}}\\[0.5em]
  \end{pmatrix},\ 
  m_{241}=\begin{pmatrix}
    \frac{c_{42}}{e_{41}} & 0\\[0.5em]
    \frac{c_{43}}{e_{41}} & 1\\[0.5em]
  \end{pmatrix},\ 
  m_{341}=\begin{pmatrix}
    \frac{c_{42}}{e_{41}} & 1\\[0.5em]
    \frac{c_{43}}{e_{41}} & 0\\[0.5em]
  \end{pmatrix},
\end{align*}
The relevant full transition matrices of the \(\Delta\)-clique network are
\begin{align*}
  M_{34}\equiv M_{2}^{(34)}&=m_{412}m_{341}m_{234}m_{123}=\begin{pmatrix}
    \alpha_{1} & \alpha_{2}\\[0.5em]
    \alpha_{3} & \alpha_{4}\\[0.5em]
  \end{pmatrix},\\[0.5em]
  M_{4}\equiv M_{2}^{(4)}&=m_{412}m_{241}m_{124}=\begin{pmatrix}
    1 & \rho_{4}\\[0.5em]
    0 & \delta_{4}\\
  \end{pmatrix},
\end{align*}
where
\begin{align*}
  \alpha_{1}&=\frac{c_{21}c_{43}}{e_{23}e_{41}} + \frac{c_{21}c_{42}t_{13}}{e_{23}e_{41}e_{12}} - \frac{c_{32}e_{24}t_{13}}{e_{34}e_{23}e_{12}} - \frac{c_{43}e_{24}t_{31}}{e_{41}e_{23}e_{34}} - \frac{c_{42}e_{24}t_{13}t_{31}}{e_{41}e_{23}e_{12}e_{34}},\\[0.5em]
  \alpha_{2}&=\frac{c_{32}t_{13}}{e_{34}e_{12}} + \frac{c_{43}t_{31}}{e_{41}e_{34}} + \frac{c_{42}t_{13}t_{31}}{e_{41}e_{12}e_{34}},\\[0.5em]
  \alpha_{3}&=\frac{c_{14}c_{21}c_{42}}{e_{12}e_{23}e_{41}} - \frac{c_{14}c_{32}e_{24}}{e_{12}e_{34}e_{23}} - \frac{c_{14}c_{42}e_{24}t_{31}}{e_{12}e_{41}e_{23}e_{34}},\\[0.5em]
  \alpha_{4}&=\frac{c_{14}c_{32}}{e_{12}e_{34}} + \frac{c_{14}c_{42}t_{31}}{e_{12}e_{41}e_{34}},\\[0.5em]
  \rho_{4}&=-\frac{e_{23}}{e_{24}} + \frac{c_{21}c_{43}}{e_{24}e_{41}} + \frac{c_{21}c_{42}t_{13}}{e_{24}e_{41}e_{12}},
\end{align*}
and \(\delta_{4}\) is the same as for the Kirk--Silber network. The remaining full transition matrices of the \(\Delta\)-clique network can be found by evaluating the following products:
\begin{align*}
  M_{1}^{(34)}&=m_{341}m_{234}m_{123}m_{412},\ M_{3}^{(34)}=m_{123}m_{412}m_{341}m_{234},\ M_{4}^{(34)}=m_{234}m_{123}m_{412}m_{341},\\
  M_{1}^{(4)}&=m_{241}m_{124}m_{412},\ M_{4}^{(4)}=m_{124}m_{412}m_{241}.
\end{align*}
From \(M_{4}^{(4)}\), we obtain the expression
\begin{equation*}
  \nu_{4}=\frac{c_{43}}{e_{41}} + \frac{c_{42}t_{13}}{e_{41}e_{12}} - \frac{c_{14}c_{42}e_{23}}{e_{12}e_{41}e_{24}}.
\end{equation*}

\subsection{The tournament network}\label{supp:maps:tournament}

The basic transition matrices of the tournament network are
\begin{align*}
  m_{312}&=\begin{pmatrix}
    \frac{c_{13}}{e_{12}} & 0 \\[0.5em]
    \frac{c_{14}}{e_{12}} & 1
  \end{pmatrix},\ m_{123}=\begin{pmatrix}
    \frac{c_{21}}{e_{23}} & 0 \\[0.5em]
    -\frac{e_{24}}{e_{23}} & 1
  \end{pmatrix},\ m_{231}=\begin{pmatrix}
    \frac{c_{32}}{e_{31}} & 0 \\[0.5em]
    -\frac{e_{34}}{e_{31}} & 1
  \end{pmatrix},\ m_{234}=\begin{pmatrix}
    1 & \frac{c_{32}}{e_{31}} \\[0.5em]
    0 & -\frac{e_{34}}{e_{31}}
  \end{pmatrix},\\[0.5em]
  m_{412}&=\begin{pmatrix}
    \frac{c_{13}}{e_{12}} & 1 \\[0.5em]
    \frac{c_{14}}{e_{12}} & 0
  \end{pmatrix},\ m_{124}=\begin{pmatrix}
    0 & \frac{c_{21}}{e_{24}} \\[0.5em]
    1 & -\frac{c_{23}}{e_{24}}
  \end{pmatrix},\ m_{241}=\begin{pmatrix}
    \frac{c_{42}}{e_{41}} & 0 \\[0.5em]
    \frac{c_{43}}{e_{41}} & 1
  \end{pmatrix},\ m_{341}=\begin{pmatrix}
    \frac{c_{42}}{e_{41}} & 1 \\[0.5em]
    \frac{c_{43}}{e_{41}} & 0
  \end{pmatrix}.
\end{align*}
The relevant full transition matrices of the tournament network are
\begin{align*}
  M_{3}\equiv M_{2}^{(3)}&=m_{312}m_{231}m_{123}=\begin{pmatrix}
    \delta_{3} & 0\\
    \rho_{3} & 1
  \end{pmatrix},\\
  M_{34}\equiv M_{2}^{(34)}&=m_{412}m_{341}m_{234}m_{123}=\begin{pmatrix}
    \alpha_{1} & \alpha_{2}\\[0.5em]
    \alpha_{3} & \alpha_{4}\\[0.5em]
  \end{pmatrix},\\[0.5em]
  M_{4}\equiv M_{2}^{(4)}&=m_{412}m_{241}m_{124}=\begin{pmatrix}
    1 & \rho_{4}\\[0.5em]
    0 & \delta_{4}\\
  \end{pmatrix},
\end{align*}
where
\begin{align*}
  \alpha_{1}&=\frac{c_{21}c_{43}}{e_{23}e_{41}} + \frac{c_{13}c_{21}c_{42}}{e_{12}e_{23}e_{41}} + \frac{c_{43}e_{24}e_{31}}{e_{41}e_{23}e_{34}} + \frac{c_{13}c_{42}e_{24}e_{31}}{e_{12}e_{41}e_{23}e_{34}} - \frac{c_{13}c_{32}e_{24}}{e_{12}e_{34}e_{23}},\\[0.5em]
  \alpha_{2}&=\frac{c_{32}c_{13}}{e_{12}e_{34}} - \frac{c_{43}e_{31}}{e_{41}e_{34}} - \frac{c_{13}c_{42}e_{31}}{e_{12}e_{41}e_{34}},\\[0.5em]
  \alpha_{3}&=\frac{c_{14}c_{21}c_{42}}{e_{12}e_{23}e_{41}} + \frac{c_{14}c_{42}e_{24}e_{31}}{e_{12}e_{41}e_{23}e_{34}} - \frac{c_{14}c_{32}e_{24}}{e_{12}e_{34}e_{23}},\\[0.5em]
  \alpha_{4}&=\frac{c_{14}c_{32}}{e_{12}e_{34}} - \frac{c_{14}c_{42}e_{31}}{e_{12}e_{41}e_{34}},\\[0.5em]
  \rho_{3}&=-\frac{e_{23}}{e_{23}} - \frac{c_{21}c_{34}}{e_{24}e_{31}} + \frac{c_{14}c_{21}c_{32}}{e_{12}e_{23}e_{31}},\\[0.5em]
  \rho_{4}&=-\frac{e_{23}}{e_{24}} + \frac{c_{21}c_{43}}{e_{24}e_{41}} + \frac{c_{13}c_{21}c_{42}}{e_{12}e_{24}e_{41}},
\end{align*}
and \(\delta_{3}\) and \(\delta_{4}\) are the same as for the Kirk--Silber network. The remaining full transition matrices of the tournament network can be found by evaluating the following products:
\begin{align*}
  M_{1}^{(3)}&=m_{123}m_{312}m_{231},\ M_{3}^{(3)}=m_{123}m_{312}m_{231},\\
  M_{1}^{(34)}&=m_{341}m_{234}m_{123}m_{412},\ M_{3}^{(34)}=m_{123}m_{412}m_{341}m_{234},\ M_{4}^{(34)}=m_{234}m_{123}m_{412}m_{341},\\
  M_{1}^{(4)}&=m_{241}m_{124}m_{412},\ M_{4}^{(4)}=m_{124}m_{412}m_{241}.
\end{align*}
From \(M_{1}^{(3)}\) and \(M_{4}^{(4)}\), we obtain, respectively, the expressions,
\begin{align*}
  \mu_{3}=\frac{c_{14}}{e_{12}} - \frac{c_{13}e_{24}}{e_{12}e_{23}} - \frac{c_{13}c_{21}e_{34}}{e_{12}e_{23}e_{34}},\\
  \nu_{4}=\frac{c_{43}}{e_{41}} + \frac{c_{13}c_{42}}{e_{12} e_{41}} - \frac{c_{14}c_{42}e_{23}}{e_{12}e_{41}e_{24}}.
\end{align*}

%% file: shared.tex
\usepackage{lipsum}
\usepackage{amsfonts}
\usepackage{graphicx}
\usepackage{epstopdf}
\usepackage{algorithmic}
\ifpdf
  \DeclareGraphicsExtensions{.eps,.pdf,.png,.jpg}
\else
  \DeclareGraphicsExtensions{.eps}
\fi

\usepackage{enumitem}
\setlist[enumerate]{leftmargin=.5in}
\setlist[itemize]{leftmargin=.5in}

\newsiamremark{remark}{Remark}
\newsiamremark{hypothesis}{Hypothesis}
\crefname{hypothesis}{Hypothesis}{Hypotheses}
\newsiamthm{claim}{Claim}

\headers{Dynamics near heteroclinic networks in \(\mathbb{R}^{4}\)}{D. C. Groothuizen Dijkema, V. Kirk, and C. M. Postlethwaite}

\title{Analysis of dynamics near heteroclinic networks in \(\R^{4}\) with a projected map}

\author{
  David C Groothuizen Dijkema\thanks{Department of Mathematics, University of Auckland, Auckland, 1142, New Zealand
  (\email{david.groothuizen.dijkema@auckland.ac.nz}, \email{v.kirk@auckland.ac.nz}, \email{c.postlethwaite@auckland.ac.nz}).}
  \and Vivien Kirk\footnotemark[1]
  \and Claire M Postlethwaite\footnotemark[1]
}

\usepackage{amsopn}

%% file: 0_commands.tex
\usepackage{xifthen}

\usepackage{lipsum}

\usepackage{enumitem}
\setlist[enumerate,1]{label={\arabic*.}}
\setlist[enumerate,2]{label={\alph*)},ref={\theenumi\alph*)}}

\usepackage{amsopn}
\usepackage{amssymb}
\usepackage{amsmath}
\usepackage{mathtools}

\usepackage{multirow}

\usepackage{float}
\usepackage{graphicx}
\usepackage{caption}
\usepackage[labelformat=simple]{subcaption}

\usepackage{tikz}
\usetikzlibrary{positioning}
\usetikzlibrary{arrows}
\usetikzlibrary{shapes.geometric}
\usetikzlibrary{arrows.meta}
\usetikzlibrary{decorations.markings}

\usepackage[nameinlink]{cleveref}

\usepackage{xcolor}

\definecolor{myblue}{HTML}{0040FF}
\definecolor{myamber}{HTML}{FFBF00}

\definecolor{xi_orange}{HTML}{FF7F00}
\definecolor{xi_green}{HTML}{66CC00}
\definecolor{xi_blue}{HTML}{0080FF}
\definecolor{xi_pink}{HTML}{FF0080}
\definecolor{xi_dark_green}{HTML}{00B359}
\definecolor{xi_purple}{HTML}{7F00FF}

\DeclareMathOperator{\N}{\mathbb{N}}
\DeclareMathOperator{\Z}{\mathbb{Z}}

\DeclareMathOperator{\R}{\mathbb{R}}

\newcommand{\Rn}{\R^{n}}

\newcommand{\Rp}{\R^{p}}

\def\Rminusscale{0.7}
\newcommand{\negativeR}[1]{\R_{\scalebox{\Rminusscale}[1.0]{--}}^{#1}}

\newcommand{\cyc}[1]{\mathcal{C}_{#1}}
\newcommand{\net}[1]{\mathcal{N}_{#1}}

\newcommand{\poinin}[2][]{\ifthenelse{\isempty{#1}}{\mathbf{H}_{#2}^{\mathrm{in}}}{\mathbf{H}_{#2}^{\mathrm{in},#1}}}
\newcommand{\poinout}[2][]{\ifthenelse{\isempty{#1}}{\mathbf{H}_{#2}^{\mathrm{out}}}{\mathbf{H}_{#2}^{\mathrm{out},#1}}}

\newcommand{\lmax}{\lambda_{\max}}
\newcommand{\wmax}{w_{\max}}

\newcommand{\basin}[1]{\mathcal{B}(#1)}
\newcommand{\dlocbasin}[2]{\mathcal{B}_{#1}(#2)}

\newcommand{\mathdef}[1]{\textbf{\textit{#1}}}

\definecolor{mytextgreen}{HTML}{28A428}
\definecolor{mytextblue}{HTML}{3A53F8}
\definecolor{mytextred}{HTML}{F83A53}

\newlength{\myleftlen}

%% file: 1_introduction.tex
\section{Introduction}\label{sec:intro}

Heteroclinic cycles and networks are flow-invariant structures in dynamical systems. They are composed of invariant sets and heteroclinic orbits between these invariant sets. In dynamical systems with appropriate invariant subspaces, they can occur robustly. This invariance could arise from, for example, symmetries of the dynamical system, or properties of an underlying physical model, such as the invariance of extinction in a population model.

In this paper, we consider only heteroclinic cycles between saddle equilibria. A trajectory attracted to such a cycle spends a period of time in a small neighbourhood of one equilibrium, before switching to the next equilibrium in the cycle. The period of time spent near each equilibrium increases as the trajectory approaches the cycle. Trajectories thus exhibit a particular type of intermittent behaviour. The prototypical example of a heteroclinic cycle is the Guckenheimer-Holmes cycle \cite{guckenheimer_holmes_1988}, containing three equilibria. Since this cycle was initially studied \cite{may_leonard_1975}, and it was shown heteroclinic cycles can be structurally stable \cite{field_1980,dos_Reis_1984,guckenheimer_holmes_1988}, many results regarding the existence, stability, and dynamical properties of various heteroclinic cycles have been derived. See, for example, any of \cite{melbourne_1991,field_swift_1991,krupa_melbourne_1995a,chossat_krupa_melbourne_scheel_1997,krupa_melbourne_2004,postlethwaite_dawes_2005,postlethwaite_dawes_2006,postlethwaite_2010,podvigina_ashwin_2011,podvigina_2012,lohse_2015,podvigina_chossat_2015}, and the references therein.

A connected union of heteroclinic cycles is known as a \textit{heteroclinic network}. Early studies of heteroclinic networks include those of Kirk and Silber \cite{kirk_silber_1994} and Brannath \cite{brannath_1994}. A variety of different examples of heteroclinic networks have since been studied, and some general results about asymptotic stability have been derived. See, for example \cite{aguiar_castro_labouriau_2005,postlethwaite_dawes_2005,kirk_lane_postlethwaite_rucklidge_silber_2010,kirk_postlethwaite_rucklidge_2012,castro_lohse_2014,afraimovich_moses_young_2016,podvigina_castro_labouriau_2020,castro_garrido_da_silva_2022,postlethwaite_rucklidge_2022,castro_ferreira_garrido_da_silva_labouriau_2022,podvigina_2023}, and references therein. The usual approach to analysing heteroclinic cycles and networks is to construct return maps to cross-sections defined near each equilibrium. These return maps approximate the dynamics of trajectories near the cycle or network. The stability of a cycle and the geometry of the cycle's basin of attraction can then be determined, as well as whether trajectories switch from cycling around one subcycle of a network to a different subcycle. These calculations are often performed by moving the return map to logarithmic coordinates, producing a \textit{transition matrix}. Studying heteroclinic cycles and networks with return maps or their transition matrices produces a lower-dimensional discrete map from the continuous flow of the dynamical system, simplifying the analysis. Stability indices---which measure the relative size of the basin of attraction---can also be calculated from the transition matrices \cite{podvigina_ashwin_2011}.

In this paper, we introduce an extension to the usual methods of analysis of heteroclinic networks, which we call the \textit{projected map}. This map is a discrete, piecewise-smooth dynamical system with dimension one fewer than the order of the transition matrix. The projected map allows us to not only determine the stability of heteroclinic cycles in a network, but also more easily classify the behaviour of trajectories near the network. A similar methodology was developed by Peixe and Rodrigues \cite{peixe_rodrigues_2022}, who studied a heteroclinic network on the unit cube in \(\R^{3}\), composed of six equilibria and six cycles, modelled by polymatrix replicator equations.

We study three heteroclinic networks in \(\R^{4}\) composed of four equilibria: the Kirk--Silber, \(\Delta\)-clique, and tournament networks. (See \Cref{fig:R4_networks}.) These three networks are, qualitatively, the only networks of four equilibria in \(\R^{4}\) composed of quasi-simple cycles of length three or four \cite{castro_lohse_2016a}. Analysis of the projected map of each heteroclinic network allows us to classify the dynamics of trajectories in all sufficiently small neighbourhoods of the network, for all parameter values of the system.

In \cite{brannath_1994}, Brannath studied these networks as dynamical systems on the tetrahedron and considered whether trajectories could be asymptotic to more than one cycle of any of these networks. They showed, in agreement with the results of \cite{kirk_silber_1994}, that this cannot happen near the Kirk--Silber network. Brannath was also able to show, assuming certain relations between parameter values, that this was also the case for the \(\Delta\)-clique network and the tournament network. They conjectured generally that trajectories are asymptotic to at most one cycle for all parameter values, writing \textit{``...we could ask for the existence of orbits whose \(\omega\)-limits contain more than one heteroclinic cycle... The author thinks that there should be a quite general (but non-trivial) argument which excludes such a `switching' for networks on \(S_{3}\)...''} \cite{brannath_1994}. In this paper, we are able to resolve this claim. Our analysis enables us to conclude that a trajectory asymptotic to one of these three networks is asymptotic to at most one subcycle, and cannot regularly or irregularly cycle between cycles, assuming only that the network has no positive globally transverse eigenvalues. In comparison, analysis of networks in \(\R^{5}\) and \(\R^{6}\) shows examples of networks where trajectories switch between cycles of the network an infinite number of times \cite{postlethwaite_rucklidge_2022,postlethwaite_dawes_2005,podvigina_2023,castro_lohse_2016b}.

\input{figure_1.tex}
Our analysis of the Kirk--Silber reproduces results previously derived in \cite{kirk_silber_1994} by analysing return maps. However, we consider the additional case where trajectories switch from an unstable cycle to a stable cycle. Analysis of trajectories near the \(\Delta\)-clique network includes some subtleties regarding the stability of fixed points of the projected map. The analysis of the tournament network is the most subtle. Although this network only contains three subcycles, these cycles form three subnetworks, one equivalent to the Kirk--Silber network, and two equivalent to the \(\Delta\)-clique network. These complications do not lead to particularly rich dynamics, but numerous different cases for various different relationships between parameters of the dynamical system must be considered. We therefore give a more qualitative overview of the analysis of the tournament network. Stability of cycles in the tournament network under some assumptions on the parameters of the dynamical system has previously been studied by Castro, Ferreira, and Labouriau \cite{castro_ferreira_labouriau_2024}. For all three networks, we find that the onset of switching between cycles often arises as a result of a \textit{border-collision bifurcation} of the projected map.

The remainder of this paper proceeds as follows. \Cref{sec:background} defines the three different heteroclinic networks we will study, and provides relevant definitions. We overview how return maps are constructed in \cref{sec:ret_proj_map} by considering the Kirk--Silber network, and then define the projected map of this network. In \cref{sec:anal}, we use the projected map to study the dynamics of trajectories near the three networks we consider. We generalise some of our results to arbitrary heteroclinic networks in \cref{sec:proj_map_BCB}. \Cref{sec:disc} concludes. \Cref{app:well_defness} gives some technical details of the definition of the projected map.

%% file: figure_1.tex
\begin{figure}
  \centering
  \newlength{\thickness}
  \setlength{\thickness}{0.5mm}

  \newlength{\xdisp}
  \setlength{\xdisp}{1.3cm}
  \newlength{\ydisp}
  \setlength{\ydisp}{1cm}
  \newlength{\ydispp}
  \setlength{\ydispp}{2.2cm}

  \begin{subfigure}[t]{0.33\textwidth}
    \centering
    \begin{tikzpicture}
      \node [circle,minimum size=2.5mm,inner sep=0pt,draw=xi_pink,fill=xi_pink,    label={[font=\normalsize, label distance=-0.5mm]above:\(\xi_{1}\)}]                                        (one)   {};
      \node [circle,minimum size=2.5mm,inner sep=0pt,draw=xi_green,fill=xi_green,  label={[font=\normalsize, label distance=-0.5mm]below:\(\xi_{2}\)},below = \ydispp of one]                 (two)   {};
      \node [circle,minimum size=2.5mm,inner sep=0pt,draw=xi_blue,fill=xi_blue,    label={[font=\normalsize, label distance=-1.0mm]left: \(\xi_{3}\)},below left  = \ydisp and \xdisp of one] (three) {};
      \node [circle,minimum size=2.5mm,inner sep=0pt,draw=xi_orange,fill=xi_orange,label={[font=\normalsize, label distance=-0.5mm]right:\(\xi_{4}\)},below right = \ydisp and \xdisp of one] (four)  {};

      \draw [-{Straight Barb},line width=\thickness,shorten >=5pt,shorten <=5pt] (one)   to [] node[below] {} (two);
      \draw [-{Straight Barb},line width=\thickness,shorten >=5pt,shorten <=5pt] (two)   to [] node[below] {} (three);
      \draw [-{Straight Barb},line width=\thickness,shorten >=5pt,shorten <=5pt] (three) to [] node[below] {} (one);
      \draw [-{Straight Barb},line width=\thickness,shorten >=5pt,shorten <=5pt] (two)  to [] node[below] {} (four);
      \draw [-{Straight Barb},line width=\thickness,shorten >=5pt,shorten <=5pt] (four) to [] node[below] {} (one);
    \end{tikzpicture}
    \caption{The Kirk--Silber network.}
    \label{fig:R4_networks:kirk_silber_net}
  \end{subfigure}%
  \begin{subfigure}[t]{0.33\textwidth}
    \centering
    \begin{tikzpicture}
      \node [circle,minimum size=2.5mm,inner sep=0pt,draw=xi_orange,fill=xi_orange,label={[font=\normalsize, label distance=-0.5mm]above:\(\xi_{4}\)}]                                         (four)  {};
      \node [circle,minimum size=2.5mm,inner sep=0pt,draw=xi_green,fill=xi_green,  label={[font=\normalsize, label distance=-0.5mm]below:\(\xi_{2}\)},below = \ydispp of four]                 (two)   {};
      \node [circle,minimum size=2.5mm,inner sep=0pt,draw=xi_pink,fill=xi_pink,    label={[font=\normalsize, label distance=-1.0mm]left: \(\xi_{1}\)},below left  = \ydisp and \xdisp of four] (one)   {};
      \node [circle,minimum size=2.5mm,inner sep=0pt,draw=xi_blue,fill=xi_blue,    label={[font=\normalsize, label distance=-0.5mm]right:\(\xi_{3}\)},below right = \ydisp and \xdisp of four] (three) {};

      \draw [-{Straight Barb},line width=\thickness,shorten >=5pt,shorten <=5pt] (one)   to [] node[below] {} (two);
      \draw [-{Straight Barb},line width=\thickness,shorten >=5pt,shorten <=5pt] (two)   to [] node[below] {} (three);
      \draw [-{Straight Barb},line width=\thickness,shorten >=5pt,shorten <=5pt] (two)   to [] node[below] {} (four);
      \draw [-{Straight Barb},line width=\thickness,shorten >=5pt,shorten <=5pt] (three) to [] node[below] {} (four);
      \draw [-{Straight Barb},line width=\thickness,shorten >=5pt,shorten <=5pt] (four)  to [] node[below] {} (one);
    \end{tikzpicture}
    \caption{The \(\Delta\)-clique network.}
    \label{fig:R4_networks:delta_clique_net}
  \end{subfigure}%
  \begin{subfigure}[t]{0.33\textwidth}
    \centering
    \begin{tikzpicture}
      \node [circle,minimum size=2.5mm,inner sep=0pt,draw=xi_pink,fill=xi_pink,    label={[font=\normalsize, label distance=-0.5mm]above:\(\xi_{1}\)}]                                        (one)   {};
      \node [circle,minimum size=2.5mm,inner sep=0pt,draw=xi_green,fill=xi_green,  label={[font=\normalsize, label distance=-0.5mm]below:\(\xi_{2}\)},below = \ydispp of one]                 (two)   {};
      \node [circle,minimum size=2.5mm,inner sep=0pt,draw=xi_blue,fill=xi_blue,    label={[font=\normalsize, label distance=-1.0mm]left: \(\xi_{3}\)},below left  = \ydisp and \xdisp of one] (three) {};
      \node [circle,minimum size=2.5mm,inner sep=0pt,draw=xi_orange,fill=xi_orange,label={[font=\normalsize, label distance=-0.5mm]right:\(\xi_{4}\)},below right = \ydisp and \xdisp of one] (four)  {};

      \draw [-{Straight Barb},line width=\thickness,shorten >=5pt,shorten <=5pt] (one)   to [] node[below] {} (two);
      \draw [-{Straight Barb},line width=\thickness,shorten >=5pt,shorten <=5pt] (two)   to [] node[below] {} (three);
      \draw [-{Straight Barb},line width=\thickness,shorten >=5pt,shorten <=5pt] (three) to [] node[below] {} (one);
      \draw [-{Straight Barb},line width=\thickness,shorten >=5pt,shorten <=5pt] (two)   to [] node[below] {} (four);
      \draw [-{Straight Barb},line width=\thickness,shorten >=5pt,shorten <=5pt] (four)  to [] node[below] {} (one);
      \draw [-{Straight Barb},line width=\thickness,shorten >=5pt,shorten <=5pt] (three) to [] node[below] {} (four);
    \end{tikzpicture}
    \caption{The tournament network.}
    \label{fig:R4_networks:tournament_net}
  \end{subfigure}
  \caption{Diagrammatic representations of the three heteroclinic networks between four equilibria in \(\R^{4}\). Coloured vertices represent hyperbolic saddle equilibria, and directed edges represent robust heteroclinic orbits.}
  \label{fig:R4_networks}
\end{figure}
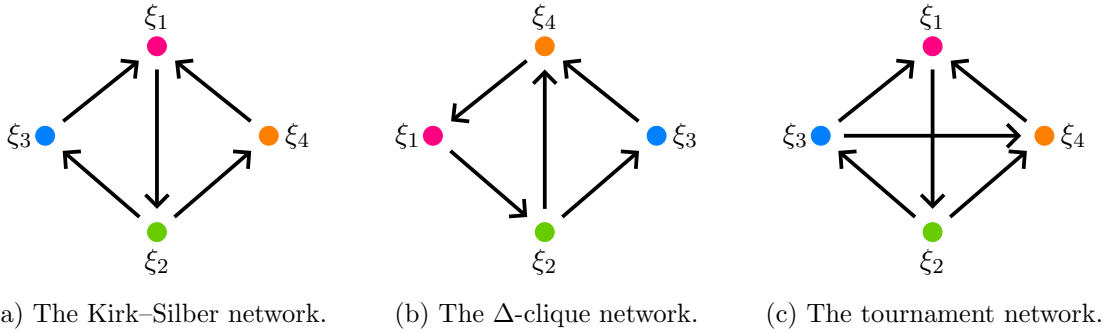

%% file: 2_background.tex
\section{Preliminaries and definitions}\label{sec:background}

In this section, we define the three heteroclinic networks we consider in this paper, and provide relevant definitions to classify the dynamics of trajectories near heteroclinic networks and the stability of heteroclinic cycles.

\subsection{Heteroclinic cycles and networks}\label{sec:background:ssec:cycs_nets}

We consider continuous-time dynamical systems of the form
\begin{equation}\label{eqn:dyn_sys}
  \dot{x}=f(x),
\end{equation}
where \(x\in\R^{4}\) and \(f\colon\R^{4}\to\R^{4}\) is an (at least) \(C^{1}\) smooth vector field.

For a point \(x_{0}\in\Rn\), let \(\phi(x_{0},t)\) be the flow generated by \cref{eqn:dyn_sys} through \(x_{0}\). Let \(\xi_{j}\) and \(\xi_{k}\) be two equilibria of \cref{eqn:dyn_sys}. A \mathdef{heteroclinic orbit} \(\phi_{jk}(t)\) is a solution to \cref{eqn:dyn_sys} such that \(\phi_{jk}(t)\to\xi_{j}\) as \(t\to-\infty\) and \(\phi_{jk}(t)\to\xi_{k}\) as \(t\to\infty\). We will also write \(\xi_{j}\to\xi_{k}\) for the heteroclinic orbit \(\phi_{jk}(t)\).

\begin{definition}
  A \mathdef{heteroclinic cycle} is a set \(\cyc{}\subseteq\Rn\) that is the union of a set of \(m\) distinct hyperbolic saddles \(\left\{\xi_{1},\ldots,\xi_{m}\right\}\) and a set of heteroclinic orbits \(\left\{\phi_{12}(t),\ldots,\phi_{m1}(t)\right\}\), where \(\phi_{jk}(t)\) is a heteroclinic orbit from \(\xi_{j}\) to \(\xi_{k}\).
\end{definition}
When discussing the equilibria of heteroclinic cycles, we will implicitly assume all operations involving indices are taken modulo \(m\); that is, \(\xi_{m+1}\equiv\xi_{1}\).

\begin{definition}
  Let \(\cyc{1},\ldots,\cyc{s}\), for \(s\geq 2\), be a collection of heteroclinic cycles. Then the union \(\net{}=\bigcup_{j=1}^{s}\cyc{j}\) is a \mathdef{heteroclinic network} if, for every pair of equilibria, there exists a sequence of heteroclinic orbits connecting them; that is, if \(\xi_{a},\xi_{b}\in\net{}\), then there exists a set of equilibria \(\left\{\xi_{\ell_{1}},\dots,\xi_{\ell_{p}}\right\}\subseteq\net{}\) and heteroclinic orbits \(\left\{\phi_{\ell_{1}\ell_{2}}(t),\dots,\phi_{\ell_{p}\ell_{1}}(t)\right\}\subseteq\net{}\), such that \(\xi_{\ell_{1}}\equiv\xi_{a}\), \(\xi_{\ell_{p}}\equiv\xi_{b}\), and \(\phi_{\ell_{i}\ell_{j}}(t)\) is a heteroclinic orbit from \(\xi_{\ell_{i}}\) to \(\xi_{\ell_{j}}\).
\end{definition}

In this paper, we are interested in three networks in \(\R^{4}\) that are each the union of cycles containing three or four equilibria. Specifically, these cycles are the \(B_{3}^{-}\) cycle and the \(C_{4}^{-}\) cycle in the classification of cycles in \(\R^{4}\) by Krupa and Melbourne \cite{krupa_melbourne_2004}. The three heteroclinic networks we consider are:
\begin{enumerate}
  \item The \textit{Kirk--Silber network} \cite{kirk_silber_1994}, the union of two cycles: \(\cyc{[3]}=\xi_{1}\to\xi_{2}\to\xi_{3}\to\xi_{1}\) and \(\cyc{[4]}=\xi_{1}\to\xi_{2}\to\xi_{4}\to\xi_{1}\);
  \item The \textit{\(\Delta\)-clique network}, the union of two cycles: \(\cyc{[4]}=\xi_{1}\to\xi_{2}\to\xi_{4}\to\xi_{1}\) and \(\cyc{[34]}=\xi_{1}\to\xi_{2}\to\xi_{3}\to\xi_{4}\to\xi_{1}\); and
  \item The \textit{tournament network}, the union of three cycles: \(\cyc{[3]}=\xi_{1}\to\xi_{2}\to\xi_{3}\to\xi_{1}\), \(\cyc{[4]}=\xi_{1}\to\xi_{2}\to\xi_{4}\to\xi_{1}\), and \(\cyc{[34]}=\xi_{1}\to\xi_{2}\to\xi_{3}\to\xi_{4}\to\xi_{1}\).
\end{enumerate}
\Cref{fig:R4_networks} gives diagrammatic representations of these three networks. The Kirk--Silber network was first studied in \cite{kirk_silber_1994}. The directed triangle of heteroclinic orbits  \(\xi_{2}\to\xi_{4}\) and \(\xi_{2}\to\xi_{3}\to\xi_{4}\) in the second network was termed a \(\Delta\)-clique in \cite{ashwin_castro_lohse_2020}, and so we name it the \textit{\(\Delta\)-clique network}. The graph of the third network is the only one of the three networks that forms a tournament in the graph-theoretic sense, and so we name it the \textit{tournament network}.

These three networks exist in dynamical systems where, for instance, \cref{eqn:dyn_sys} takes the form
\begin{equation}\label{eqn:general_form_sq}
  \dot{x}_{j}=x_{j}\bigg(1-\left\lVert x\right\rVert_{2}^{2}+\sum_{k\neq j}\alpha_{kj}x_{k}^{2}\bigg).
\end{equation}
It is straight-forward to verify that these equations have a saddle equilibrium \(\xi_{j}\) on each of the four coordinate axes \(L_{j}\) at unit distance from the origin. Moreover, they are \(\Z_{2}^{4}\)-equivariant, where the symmetries are reflections in each of the four coordinate-hyperplanes. These symmetries ensure that the coordinate axes \(L_{j}\) and coordinate hyperplanes \(P_{jk}\) are flow-invariant. Specific choices in the case of each network of the sign of the parameters \(\alpha_{kj}\) ensures that certain equilibria are a saddle or a sink in certain subspaces, and thus that the desired heteroclinic orbits exist between them. In \cref{eqn:ks_dyn_sys}, \cref{eqn:d_clique_dyn_sys}, and \cref{eqn:tournament_dyn_sys} we give the specific form of \cref{eqn:general_form_sq} that contains each network.

Saddle-sink connections in invariant subspaces are robust to sufficiently small perturbations that preserve these invariances. As a result, the three networks are robust, and occur as codimension-zero phenomena; that is, they exist for an open set of parameter values. Other examples of ODEs also exist that contain robust examples of the three networks we consider.

\begin{table}[t]
  \centering
  \caption{Classification of eigenvalues of \(Df(\xi_{k})\) \cite{krupa_melbourne_2004}.}\label{tbl:eivals}
  \vspace{5pt}
  \begin{tabular}{l l}
    \hline
    Eigenvalue class & Eigenvector subspace\\
    \hline
    Radial      (\(r\)) & \(V_{k}(r)=L_{k}\equiv P_{jk}\cap P_{kl}\)\\
    Contracting (\(c\)) & \(V_{k}(c)=P_{jk}\cap L_{k}^{\perp}\)\\
    Expanding   (\(e\)) & \(V_{k}(e)=P_{kl}\cap L_{k}^{\perp}\)\\
    Transverse  (\(t\)) & \(V_{k}(t)=(P_{jk}+P_{kl})^{\perp}\)\\
    \hline
  \end{tabular}
\end{table}

Let \(\xi_{k}\) be an equilibrium of a heteroclinic cycle. In \cite{krupa_melbourne_2004}, Krupa and Melbourne classify the eigenvalues of \(Df(\xi_{k})\) as radial, contracting, expanding, or transverse, depending on the subspace in which their corresponding eigenspace lies; this classification is given in \Cref{tbl:eivals}. For the purposes of this classification, we assume that the \(\xi_{j}\) is the previous equilibrium in the cycle before \(\xi_{k}\), and that \(\xi_{l}\) is the next equilibrium after \(\xi_{k}\).

The Kirk--Silber network and \(\Delta\)-clique network contain equilibria that are a part of only one cycle. We assume that the transverse eigenvalues of these equilibria---which are \textit{globally transverse} in the terminology of \cite{podvigina_castro_labouriau_2020,castro_garrido_da_silva_2022}---are negative. The sign of every other eigenvalue is fixed by the choices of the parameters \(\alpha_{kj}\) needed to ensure the existence of the given heteroclinic network.

\subsection{Itineraries of trajectories near heteroclinic networks}

Let \(\net{}\) be one of the three heteroclinic networks we study in this paper. For a trajectory sufficiently close to \(\net{}\), we define its \mathdef{itinerary} as the indices of the sequence of equilibria the trajectory visits; that is, \(\mathcal{W}=\left(a_{n}\right)\), where \(a_{n}\in\left\{1,2,3,4\right\}\). If the trajectory leaves a small neighbourhood of the network, the itinerary is finite, whereas a trajectory that is asymptotic to the network has an itinerary indexed by \(\N\). An itinerary \(\mathcal{W}\) is \mathdef{eventually periodic with minimal period \(p\) and with root sequence \(w\in \left\{1,2,3,4\right\}^{p}\)} if there exists a positive integer \(N\) such that \(\mathcal{W}_{N+j}=w_{j}\) for \(1\leq j\leq p\) and, for all \(n>N\), \(\mathcal{W}_{n+p}=\mathcal{W}_{n}\), and there does not exist a \(p^{\prime}<p\) such that \(\mathcal{W}_{n+p^{\prime}}=\mathcal{W}_{n}\) \cite{postlethwaite_rucklidge_2022}.

If \(\cyc{}=\xi_{j_{1}}\to\xi_{j_{2}}\to\cdots\to\xi_{j_{m}}\to\xi_{j_{1}}\) is a cycle of \(\mathcal{N}\), let \(\mathcal{X}=(j_{1}j_{2}\ldots j_{m})\). If a trajectory is asymptotic to \(\cyc{}\), then its itinerary is eventually periodic with minimal period \(m\) and root sequence some circular shift of \(\mathcal{X}\). We will say that \(\mathcal{X}\), and every circular shift thereof, are the root sequences of \(\cyc{}\).

\begin{definition}
  If \(\cyc{1}\) and \(\cyc{2}\) are cycles of \(\mathcal{N}\), then a trajectory \(\phi(x_{0},t)\) \mathdef{switches} between \(\cyc{1}\) and \(\cyc{2}\) if there exists a root sequence \(\mathcal{X}_{1}\) of \(\cyc{1}\) and \(\mathcal{X}_{2}\) of \(\cyc{2}\), such that the itinerary of \(\phi(x_{0},t)\) has \((\mathcal{X}_{1}\mathcal{X}_{2})\) as a subsequence.
\end{definition}
Therefore, for a trajectory to have switched between cycles, it must, starting in a small neighbourhood of an equilibrium they have in common, make at least one full excursion around one cycle before returning to that equilibrium and then make at least one full excursion around the other cycle before returning to that equilibrium. We emphasise that our definition of \textit{switching} is a property of a particular \textit{trajectory}, and not a property of the heteroclinic network. This definition is unrelated to the various notions of \textit{switching at a node}, \textit{switching along a connection}, or \textit{finite switching} and \textit{infinite switching near a network} used in papers such as \cite{aguiar_castro_labouriau_2005,castro_lohse_2016b,castro_garrido_da_silva_2022}. Here, \textit{switching} is used as it was in \cite{kirk_silber_1994}. Sustained switching between cycles of a network has been studied in papers such as \cite{postlethwaite_dawes_2005,postlethwaite_rucklidge_2022,podvigina_2023}.

\subsection{Stability of invariant sets}\label{sec:background:ssec:stab}

Stability properties of cycles and networks can be subtle, because positive eigenvalues affect the geometry of the basin of attraction of the cycle or network. Again, let \(\net{}\) be one of the three heteroclinic networks.

The \mathdef{basin of attraction} of \(\net{}\), \(\basin{\net{}}\), is the set of all \(x\in\Rn\) such that \(\omega(x)\subseteq\net{}\):
\begin{equation*}
  \basin{\net{}}=\left\{x\in\Rn\mid\omega(x)\subseteq \net{}\right\},
\end{equation*}
where \(\omega(x)\) is the \(\omega\)-limit set of \(x\). The basin of attraction of \(\net{}\) may include points that initially move away from \(\net{}\). Therefore, we also define, for all \(\delta>0\), the \mathdef{\(\delta\)-local basin of attraction} of \(\net{}\), \(\dlocbasin{\delta}{\net{}}\), as the set of all \(x\in\Rn\) such that \(\omega(x)\subseteq \net{}\) and that remain within distance \(\delta\) of \(\net{}\) in forward time:
\begin{equation*}
  \dlocbasin{\delta}{\net{}}=\left\{x\in\Rn\mid\omega(x)\subseteq \net{}\textrm{ and }d(\phi(x,t),\net{})<\delta\textrm{ for all }t\geq 0\right\}.
\end{equation*}
In \cite[Definition~2]{podvigina_2012}, Podvigina introduced the definition of \mathdef{fragmentary asymptotic stability} (f.a.s.): \(\net{}\) is f.a.s. if \(\mu(\dlocbasin{\delta}{\net{}})>0\) for all \(\delta>0\), where \(\mu\) is the Lebesgue measure.

If \(\mu(\basin{\net{}})=0\), but \(\basin{\net{}}\setminus \net{}\neq\emptyset\), \(\net{}\) is \mathdef{almost completely unstable} (a.c.u.), and if there exists a neighbourhood \(U\) of \(\net{}\) such that, for all \(x\in U\setminus\mathcal{N}\), \(\phi(x,t)\notin U\) for some \(t>0\), then \(\net{}\) is \mathdef{completely unstable} (c.u.).

%% file: 3_projected_map.tex
\section{Return maps and the projected map of the Kirk--Silber network}\label{sec:ret_proj_map}

In this section, we construct the projected map of the Kirk--Silber network. We will first outline the usual process of analysing the stability of heteroclinic cycles and networks \cite{kirk_silber_1994,postlethwaite_2010} in \cref{sec:ret_proj_map:ssec:ret_maps}. We give an explanation of the purpose of the projected map in \cref{sec:ret_proj_map:ssec:explanation}, and then derive the projected map itself in \cref{sec:ret_proj_map:ssec:proj_map}. The process outlined in this section can be used to derive the projected map of the other two heteroclinic networks we study.

\subsection{Return maps near the Kirk--Silber network}\label{sec:ret_proj_map:ssec:ret_maps}

For simplicity of the calculations in this section, we assume, without loss of generality, that \cref{eqn:dyn_sys} takes the form in \cref{eqn:general_form_sq}.

Consider, for example, \(\xi_{3}=(0,0,1,0)\). We define local coordinates by translating \(\xi_{3}\) to the origin and setting \(u_{3}=x_{3}-1\). Rescaling the coordinates as necessary to ensure we are in a neighbourhood of approximate linear flow, we define two cross-sections in a small neighbourhood of \(\xi_{3}\):
\begin{equation*}
  \poinin[2]{3}=\left\{\left(x_{1},x_{2},u_{3},x_{4}\right)\bigm\vert|u_{3}|<1,x_{2}=1,0<x_{1},x_{4}<1\right\}
\end{equation*}
and
\begin{equation*}
  \poinout[1]{3}=\left\{\left(x_{1},x_{2},u_{3},x_{4}\right)\bigm\vert|u_{3}|<1,x_{1}=1,0<x_{2},x_{4}<1\right\}.
\end{equation*}
The cross-section \(\poinin[2]{3}\) is transverse to the orbit \(\xi_{2}\to\xi_{3}\), and \(\poinout[1]{3}\) is transverse to \(\xi_{3}\to\xi_{1}\).

In the system \cref{eqn:dyn_sys}, there is an attracting invariant sphere \cite{field_1996}, and so the radial coordinate \(u_{3}\) will not affect the stability of the cycle \cite{kirk_silber_1994,krupa_1997,podvigina_2012}. Thus, we can consider two-dimensional subsets of the cross-sections, and the relevant coordinates are \(\left(x_{1},x_{4}\right)\in\poinin[2]{3}\) and \(\left(x_{2},x_{4}\right)\in\poinout[1]{3}\). The local map \(\psi_{231}\colon\poinin[2]{3}\to\poinout[1]{3}\) approximates the dynamics of trajectories between these two cross-sections using the linear flow near \(\xi_{3}\) and the residence time of a trajectory between the two sections, and is given by
\begin{equation*}
  \psi_{231}(x_{1},x_{4})=\left(x_{1}^{\frac{c_{32}}{e_{31}}},x_{4}\,x_{1}^{\frac{t_{34}}{e_{31}}}\right).
\end{equation*}

To approximate the dynamics near the orbit \(\xi_{3}\to\xi_{1}\), we define a cross-section in a small neighbourhood of \(\xi_{1}\), transverse to the heteroclinic orbit \(\xi_{3}\to\xi_{1}\):
\begin{equation*}
  \poinin[3]{1}=\left\{\left(u_{1},x_{2},x_{3},x_{4}\right)\bigm\vert|u_{1}|<1,x_{3}=1,0<x_{2},x_{4}<1\right\}.
\end{equation*}
Again, we can neglect the radial coordinate, now \(u_{1}\), and the relevant coordinates are \(\left(x_{2},x_{4}\right)\in\poinin[3]{1}\), the same as those on \(\poinout[1]{3}\). By the flow-invariance of the coordinate hyperplanes, the global map \(\Psi_{31}\colon\poinout[1]{3}\to\poinin[3]{1}\) can be written, to lowest order, as
\begin{equation*}
  \Psi_{31}\left(x_{2},x_{4}\right)=\left(\alpha_{31}x_{2},\beta_{31}x_{4}\right),
\end{equation*}
where \(\alpha_{31}\) and \(\beta_{31}\) are order \(1\) positive constants.

We then define the \textit{basic map} \cite{garrido_da_silva_castro_2019} \(\varphi_{231}=\Psi_{31}\circ\psi_{231}\colon\poinin[2]{3}\to\poinin[3]{1}\). This process is repeated for all equilibria and appropriate cycles, giving a total of six basic maps. These maps are then composed to give \textit{full return maps} \(\Phi_{j}\colon\poinin[i]{j}\to\poinin[i]{j}\). Full details are given in the Supplementary Material (section \ref{supp:maps}). We write \(\Phi_{j}^{(3)}\) or \(\Phi_{j}^{(4)}\) for a return map to an incoming cross-section defined near \(\xi_{j}\) that describes trajectories that make an excursion around the \(\cyc{[3]}\) or \(\cyc{[4]}\) cycle, respectively.

The basic maps and full return maps can also be expressed in logarithmic coordinates with \textit{transition matrices}. For any map \(F\colon\Rp\to\Rp\) that has the form
\begin{equation}\label{eqn:gen_map}
  F\left(x_{1},\ldots,x_{p}\right)=\left(C_{1}x_{1}^{\alpha_{11}}\ldots x_{p}^{\alpha_{1p}},\ldots,C_{p}x_{1}^{\alpha_{p1}}\ldots x_{p}^{\alpha_{pp}}\right),
\end{equation}
we define its transition matrix \(M(F)\) as
\begin{equation*}
  M(F)=\begin{pmatrix}
    \alpha_{11} & \cdots & \alpha_{1p} \\
    \vdots & \ddots & \vdots \\
    \alpha_{p1} & \cdots & \alpha_{pp} \\
  \end{pmatrix}.
\end{equation*}
It is straightforward to verify that, if \(F_{1}\) and \(F_{2}\) are functions of the form in \cref{eqn:gen_map}, then \(M(F_{1}\circ F_{2})=M(F_{1})M(F_{2})\). Full return maps (see \cref{eqn:xi_2_ret_map} for an example) have the form in \cref{eqn:gen_map}, and we write \(M_{j}^{(k)}\equiv M(\Phi_{j}^{(k)})\), and refer to these matrices as \textit{full transition matrices}. For basic maps \(\varphi_{ijk}\), we write \(m_{ijk}\equiv M(\varphi_{ijk})\), and refer to these matrices as \textit{basic transition matrices}. Since each \(\Phi_{j}\) is the composition of basic maps, the full transition matrices are products of basic transition matrices. In \cite{podvigina_2012}, Podvigina gives sufficient conditions for fragmentary asymptotic stability of a heteroclinic cycle using properties of the eigenvalues and eigenvectors of relevant transition matrices.
\begin{theorem}[{\cite[Theorem~5]{podvigina_2012}}]\label{thm:pod_bif}
  For \(1\leq j\leq m\), let \(M_{j}\) be the full transition matrices associated with a heteroclinic cycle \(\cyc{}\). Let \(\lmax\) denote the eigenvalue of \(M_{j}\) with largest absolute value, and \(\wmax\) an associated eigenvector. The heteroclinic cycle \(\cyc{}\) is fragmentarily asymptotically stable if, for all \(j\), the following three conditions hold:

  \begin{enumerate}[label=\normalfont({\Roman*})]
    \item \(\lmax\) is real;\label{cond:podI}
    \item \(\lmax>1\); and\label{cond:podII}
    \item all entries of \(\wmax\) are nonzero and have the same sign.\label{cond:podIII}
  \end{enumerate}
\end{theorem}
In \cite{podvigina_2012}, Podvigina proves that if these conditions are satisfied for all \(M_{k}\) such that \(m_{ijk}\) has a negative entry, then they are satisfied for all full transition matrices \(M_{j}\).

For the Kirk--Silber network, only \(m_{123}\) and \(m_{124}\) have negative entries (see Supplementary Material, section \ref{supp:maps:ks}), and so we only have to check the stability conditions of \Cref{thm:pod_bif} at \(\xi_{3}\) and \(\xi_{4}\). Applying the conditions of \Cref{thm:pod_bif} to the full transition matrices \(M_{3}^{(3)}\) and \(M_{4}^{(4)}\) allows us to conclude that the \(\cyc{[3]}\) cycle is f.a.s. if and only if \(\delta_{3}>1\) and \(\nu_{3}>0\), and that the \(\cyc{[4]}\) cycle is f.a.s. if and only if \(\delta_{4}>1\) and \(\nu_{4}>0\), where
\begin{alignat*}{2}
  \delta_{3}=&\frac{c_{13}c_{21}c_{32}}{e_{12}e_{23}e_{31}}>0,&&\delta_{4}=\frac{c_{14}c_{21}c_{42}}{e_{12}e_{24}e_{41}}>0,\\
  \nu_{3}=&\frac{t_{34}}{e_{31}}+\frac{c_{14}c_{32}}{e_{12}e_{31}}-\frac{c_{13}c_{32}e_{24}}{e_{12}e_{23}e_{31}},\qquad&&\nu_{4}=\frac{t_{43}}{e_{41}}+\frac{c_{13}c_{42}}{e_{12}e_{41}}-\frac{c_{14}c_{42}e_{23}}{e_{12}e_{24}e_{41}}.
\end{alignat*}
These conditions for stability were derived by Kirk and Silber in \cite{kirk_silber_1994}.

We now proceed to explain our new \textit{projected map}. The equilibrium \(\xi_{2}\) has more than outgoing heteroclinic orbit, and so we say it is a \mathdef{splitting equilibrium}. For the Kirk--Silber network, \(\xi_{2}\) is the only splitting equilibrium. As such, where a trajectory strikes \(\poinin[1]{2}\) determines which cycle the trajectory next cycles around. Therefore, to study the itinerary of trajectories near the Kirk--Silber network, we focus our analysis on \(\poinin[1]{2}\).

We define the incoming cross-section at \(\xi_{2}\) as
\begin{equation*}
  \poinin[1]{2}=\left\{(x_{1},u_{2},x_{3},x_{4})\bigm\vert|u_{2}|<1,x_{1}=1,0<x_{3},x_{4}<1\right\}.
\end{equation*}
The relevant coordinates on this section are \(x_{3}\) and \(x_{4}\). The curve \(x_{3}^{e_{24}}=x_{4}^{e_{23}}\) divides \(\poinin[1]{2}\) into trajectories that leave a small neighbourhood of \(\xi_{2}\) in the direction of \(\xi_{3}\) or \(\xi_{4}\).
Generically, we can consider the subsets of \(\poinin[1]{2}\) defined by
\begin{equation*}
  \Gamma_{3}=\left\{\left(x_{3},x_{4}\right)\in\poinin[1]{2}\mid\left(1-\epsilon\right)x_{3}^{\frac{e_{24}}{e_{23}}}>x_{4}\right\}\textrm{ and }
  \Gamma_{4}=\left\{\left(x_{3},x_{4}\right)\in\poinin[1]{2}\mid\left(1-\epsilon\right)x_{4}^{\frac{e_{23}}{e_{24}}}>x_{3}\right\},
\end{equation*}
where \(0<\epsilon<1\). The constant \(\left(1-\epsilon\right)\) excludes from \(\Gamma_{3}\) and \(\Gamma_{4}\) all points that lie too close to the separatrix dividing the basins of attraction of the \(\cyc{[3]}\) and \(\cyc{[4]}\) cycles, and so trajectories that do not strike \(\poinin[2]{3}\) or \(\poinin[2]{4}\). We refer to the subset of \(\poinin[1]{2}\) between \(\Gamma_{3}\) and \(\Gamma_{4}\) as the \textit{excluded cusp}, defined as \(\Gamma_{c}\coloneq\poinin[1]{2}\setminus\left(\Gamma_{3}\cup\Gamma_{4}\right)\). We refer to the curve defined by \(x_{3}^{e_{24}}=x_{4}^{e_{23}}\) as the \textit{switching curve}, and denote it \(\Sigma_{s}\). Furthermore, the boundary between \(\Gamma_{c}\) and \(\Gamma_{3}\) is denoted \(\Sigma_{c}^{+}\), and the boundary between \(\Gamma_{c}\) and \(\Gamma_{4}\) is denoted \(\Sigma_{c}^{-}\).

\input{figure_2.tex}
The sets \(\Gamma_{3}\) and \(\Gamma_{4}\) are the domains of definition of the local maps \(\psi_{123}\colon\Gamma_{3}\to\poinout[3]{2}\) and \(\psi_{124}\colon\Gamma_{4}\to\poinout[4]{2}\), respectively, and so also of the full return maps \(\Phi_{2}^{(3)}\colon\Gamma_{3}\to\poinin[1]{2}\) and \(\Phi_{2}^{(4)}\colon\Gamma_{4}\to\poinin[1]{2}\). The full return map \(\Phi_{2}\colon\Gamma_{3}\cup\Gamma_{4}\to\poinin[1]{2}\) is, to lowest order,
\begin{equation}\label{eqn:xi_2_ret_map}
  \Phi_{2}(x_{3},x_{4})=\begin{cases}
    \Phi_{2}^{(3)}(x_{3},x_{4})=(A_{32}x_{3}^{\delta_{3}},B_{32}x_{3}^{\rho_{3}}x_{4}), & \textrm{if } \left(x_{3},x_{4}\right)\in\Gamma_{3},\\
    \Phi_{2}^{(4)}(x_{3},x_{4})=(A_{42}x_{3}x_{4}^{\rho_{4}},B_{42}x_{4}^{\delta_{4}}), & \textrm{if } \left(x_{3},x_{4}\right)\in\Gamma_{4},
  \end{cases}
\end{equation}
where \(A_{32}\), \(B_{32}\), \(A_{42}\), and \(B_{42}\) are order \(1\) constants, and
\begin{equation*}
    \rho_{3}=-\frac{e_{24}}{e_{23}}+\frac{c_{21}t_{34}}{e_{23}e_{31}}+\frac{c_{14}c_{21}c_{32}}{e_{12}e_{23}e_{31}}\qquad\textrm{ and }\qquad\rho_{4}=-\frac{e_{23}}{e_{24}}+\frac{c_{21}t_{43}}{e_{24}e_{41}}+\frac{c_{13}c_{21}c_{42}}{e_{12}e_{24}e_{41}}.
\end{equation*}
As shown by Kirk and Silber \cite{kirk_silber_1994}, if \(\delta_{3}>1\), then \(\nu_{3}>0\) implies \(\rho_{3}>0\), and if \(\delta_{4}>1\), then \(\nu_{4}>0\) implies \(\rho_{4}>0\). From here, our primary focus will be on this return map, and so we write simply \(\Phi\equiv\Phi_{2}\), \(\Phi_{3}\equiv\Phi_{2}^{(3)}\), and \(\Phi_{4}\equiv\Phi_{2}^{(4)}\). Examples of the action of the map \(\Phi\) are shown in \Cref{fig:ks_maps_normal_coords} for both \(\nu_{3}>0\) (panel \subref{fig:ks_maps_normal_coords:subfig:pre_BCB})---where there is no switching between cycles---and \(\nu_{3}<0\) (panel \subref{fig:ks_maps_normal_coords:subfig:post_BCB})---where there is switching from the \(\cyc{[3]}\) to \(\cyc{[4]}\) cycle.

The return maps \(\Phi_{3}\) and \(\Phi_{4}\) have transition matrices
\begin{equation*}
  M_{3}=\begin{pmatrix}
    \delta_{3} & 0 \\
    \rho_{3} & 1
  \end{pmatrix}
  \qquad\textrm{ and }\qquad
  M_{4}=\begin{pmatrix}
    1 & \rho_{4} \\
    0 & \delta_{4}
  \end{pmatrix}.
\end{equation*}
Let \(X_{3}\coloneq\log x_{3}<0\) and \(X_{4}\coloneq\log x_{4}<0\). These matrices act on subsets of \(\negativeR{2}\)---the negative orthant of \(\R^{2}\)---that correspond to the sets \(\Gamma_{3}\) and \(\Gamma_{4}\) in logarithmic coordinates. We label these subsets as \(\mathcal{D}_{3}\) and \(\mathcal{D}_{4}\), respectively. We define \(\mathcal{D}_{c}\) as \(\Gamma_{c}\) in logarithmic coordinates. Setting \(w_{s}=(1,\frac{e_{24}}{e_{23}})\), which we call the \textit{switching vector}, we define the \textit{switching subspace} \(W_{s}\) as the span of \(w_{s}\), which is \(\Sigma_{s}\) in logarithmic coordinates. We also write \(W_{s}^{+}\) and \(W_{s}^{-}\) for the one-dimensional affine subspaces that correspond to the curves \(\Sigma_{c}^{+}\) and \(\Sigma_{c}^{-}\). The lines \(W_{s}^{+}\) and \(W_{s}^{-}\) are parallel to \(W_{s}\) and do not pass through the origin.

\input{figure_3.tex}
We can then define a map \(M\colon\mathcal{D}_{3}\cup\mathcal{D}_{4}\to\negativeR{2}\) by
\begin{equation*}
  M(X_{3},X_{4})=\begin{cases}
    M_{3}(X_{3},X_{4})=\begin{pmatrix}\delta_{3}&0\\\rho_{3}&1\end{pmatrix}\begin{pmatrix}X_{3}\\X_{4}\end{pmatrix}, & \textrm{if } \left(X_{3},X_{4}\right)\in\mathcal{D}_{3},\\[1.8em]
    M_{4}(X_{3},X_{4})=\begin{pmatrix}1&\rho_{4}\\0&\delta_{4}\end{pmatrix}\begin{pmatrix}X_{3}\\X_{4}\end{pmatrix}, & \textrm{if } \left(X_{3},X_{4}\right)\in\mathcal{D}_{4}.
  \end{cases}
\end{equation*}

We write \(w_{3}^{*}\) for the eigenvector \(\left(\delta_{3}-1,\rho_{3}\right)\) of \(M_{3}\) with eigenvalue \(\delta_{3}\), and \(w_{4}^{*}\) for the eigenvector \(\left(\rho_{4},\delta_{4}-1\right)\) of \(M_{4}\) with eigenvalue \(\delta_{4}\). We write \(W_{3}^{*}\) for the eigenspace of \(M_{3}\) that is the span of \(w_{3}^{*}\), and similarly for \(W_{4}^{*}\). We also write \(w_{3}^{-}=(0,1)\) and \(W_{3}^{-}\) for its span, and \(w_{4}^{-}=(1,0)\) and \(W_{4}^{-}\) for its span. These eigenvectors both have eigenvalue \(1\).

\subsection{Motivating the projected map}\label{sec:ret_proj_map:ssec:explanation}

Write \(V_{a}\) for the one-dimensional linear subspace of \(\R^{2}\) defined by \(X_{4}=aX_{3}\), where \(a>0\). Then, all vectors \(X\in V_{a}\) in the domain of \(M\) are mapped into the same one-dimensional subspace, \(V_{h(a)}\), where
\begin{equation*}
  h(a)=\begin{cases}
    h_{3}(a)=\frac{a+\rho_{3}}{\delta_{3}},& \textrm{ if } a>\frac{e_{24}}{e_{23}},\\[1em]
    h_{4}(a)=\frac{a\delta_{4}}{1+a\rho_{4}},& \textrm{ if } a<\frac{e_{24}}{e_{23}}.
  \end{cases}
\end{equation*}
The map \(h_{3}(a)\) has a fixed point \(a_{3}^{*}=\frac{\rho_{3}}{\delta_{3}-1}\), and \(V_{a_{3}^{*}}\) is the eigenspace \(W_{3}^{*}\). The map \(h_{4}(a)\) has a fixed point \(a_{4}^{*}=\frac{\delta_{4}-1}{\rho_{4}}\), and \(V_{a_{4}^{*}}\) is the eigenspace \(W_{4}^{*}\).

These properties of \(M_{3}\) and \(M_{4}\) as linear maps have analogues in the return maps \(\Phi_{3}\) and \(\Phi_{4}\), which was previously considered in \cite{postlethwaite_2010}. For \(a>0\), consider the one-dimensional subspace defined by \(X_{4}=aX_{3}\), restricted to the negative orthant. Taking exponentials gives a curve \(x_{4}=x_{3}^{a}\) in \(\poinin[1]{2}\). We note that any point \(\left(q_{3},q_{4}\right)\in\poinin[2]{1}\) lies on the curve \(x_{4}=x_{3}^{a}\), where \(a=\log q_{4}/\log q_{3}\). All points in the curve \(x_{4}=x_{3}^{a}\) that are in the domain of \(\Phi\) are mapped into the curve \(x_{4}=x_{3}^{h(a)}\). Write \(\Sigma_{3}^{*}\) and \(\Sigma_{4}^{*}\) for the curves defined by \(x_{4}=x_{3}^{a_{3}^{*}}\) and \(x_{4}=x_{3}^{a_{4}^{*}}\), respectively. These curves are invariant under the action of \(\Phi_{3}\) and \(\Phi_{4}\), respectively. The eigenspaces \(W_{3}^{*}\) and \(W_{4}^{*}\) correspond to these invariant curves in logarithmic coordinates.

If we are interested in not just the stability properties of the Kirk--Silber network or its component cycles, but also in classifying the behaviour of nearby trajectories, we can determine the itinerary of trajectories by which domain of the return map \(\Phi\) they pass through. However, we do not need to study all trajectories near the network, because trajectories in the same curve of the form \(x_{4}=x_{3}^{a}\) have the same itinerary. We instead study a map that identifies every point on each curve, or, equivalently, every vector in each linear subspace restricted to \(\negativeR{2}\). We do not study the map \(h\), but instead derive a map from the transition matrices \(M_{3}\) and \(M_{4}\), which we call the \textit{projected map}.

\subsection{The projected map of the Kirk--Silber network}\label{sec:ret_proj_map:ssec:proj_map}

We define the set
\begin{equation*}
  S=\left\{\left(X_{3},X_{4}\right)\in\negativeR{2}\mid X_{3}+X_{4}=-1\right\}.
\end{equation*}
By identifying any point in \(S\) with its first coordinate, we can, without ambiguity, identify the set \(S\) with the open interval \(\left(-1,0\right)\subseteq\negativeR{}\).

We define the projection \(\Pi\colon\negativeR{2}\to S\) of any vector \(X\in\negativeR{2}\) onto \(S\) by taking the first coordinate of the intersection of the line segment \(S\) with the subspace containing \(X\):
\begin{equation}\label{eqn:projected_point}
  \Pi\colon X=\left(X_{3},X_{4}\right)\mapsto \frac{-1}{\left(1,1\right)\cdot X}e_{1}\cdot X=\frac{-X_{3}}{X_{3}+X_{4}},
\end{equation}
where \(e_{1}\) is the vector \(\left(1,0\right)\). \Cref{fig:projection} depicts this projection. For all \(X\in\negativeR{2}\) and \(\alpha>0\),
\begin{equation*}
  \Pi(\alpha X)=\Pi(X).
\end{equation*}
%
% Therefore, all vectors in the same linear subspace are represented by the same point in \(S\).

\input{figure_4.tex}
We define the projected map as the induced action of the map \(M\) on \(S\). Let
\begin{equation*}
  \vartheta_{s}\coloneqq\Pi(w_{s})=\frac{-1}{1+\frac{e_{24}}{e_{23}}}.
\end{equation*}
For a point \(\vartheta\in S\), if \(\vartheta_{s}<\vartheta\), then
\begin{equation}\label{eqn:m3_mapped_point}
  M\left(\vartheta,-1-\vartheta\right)=M_{3}\left(\vartheta,-1-\vartheta\right)=\begin{pmatrix}
    \delta_{3} & 0 \\
    \rho_{3} & 1
  \end{pmatrix}\begin{pmatrix}
    \vartheta \\
    -1-\vartheta
  \end{pmatrix}=\begin{pmatrix}
    \delta_{3}\vartheta \\
    (\rho_{3}-1)\vartheta - 1
  \end{pmatrix}
\end{equation}
and, if \(\vartheta<\vartheta_{s}\), then 
\begin{equation}\label{eqn:m4_mapped_point}
  M\left(\alpha\left(\vartheta,-1-\vartheta\right)\right)=M_{4}\left(\vartheta,-1-\vartheta\right)=\begin{pmatrix}
    1 & \rho_{4} \\
    0 & \delta_{4}
  \end{pmatrix}\begin{pmatrix}
    \vartheta \\
    -1-\vartheta
  \end{pmatrix}=\begin{pmatrix}
    (1-\rho_{4})\vartheta-\rho_{4} \\
    \delta_{4}(-1-\vartheta)
  \end{pmatrix}
\end{equation}
We therefore define the intervals \(\Theta_{3}=\left(\vartheta_{s},0\right)\) and \(\Theta_{4}=\left(-1,\vartheta_{s}\right)\), and so \(S=\Theta_{4}\cup\left\{\vartheta_{s}\right\}\cup\Theta_{3}\).

With the action of \(M_{3}\) in \cref{eqn:m3_mapped_point} and \(M_{4}\) in \cref{eqn:m4_mapped_point}, we apply the projection in \cref{eqn:projected_point}, deriving maps
\begin{equation*}
  f_{3}=\Pi(M_{3}(\vartheta,-1-\vartheta))\colon\Theta_{3}\to S
  \qquad\textrm{ and }\qquad
  f_{4}=\Pi(M_{4}(\vartheta,-1-\vartheta))\colon\Theta_{4}\to S.
\end{equation*}
We combine these maps into the projected map \(f\colon\Theta_{3}\cup\Theta_{4}\to S\):
\begin{equation}\label{eqn:ks_proj_map}
  f(\vartheta)=\begin{cases}
    f_{3}(\vartheta)=\dfrac{-\delta_{3}\vartheta}{\left(\delta_{3}+\rho_{3}-1\right)\vartheta-1}, & \textrm{if } \vartheta\in\Theta_{3},\\[2em]
    f_{4}(\vartheta)=\dfrac{\left(1-\rho_{4}\right)\vartheta-\rho_{4}}{\left(\delta_{4}+\rho_{4}-1\right)\vartheta+\delta_{4}+\rho_{4}}, & \textrm{if } \vartheta\in\Theta_{4},
  \end{cases}
\end{equation}
The projected map \(f\) is a piecewise-smooth dynamical system, and the point \(\vartheta_{s}\) is its \textit{switching manifold}. \Cref{fig:ks_proj_map} gives schematic representations of this projected map for a variety of relations between parameters.

In \cref{app:well_defness}, we show that this map is well-defined and describes the dynamics of trajectories in all sufficiently small neighbourhoods of the Kirk--Silber network.

%% file: figure_2.tex
\begin{figure}
  \centering
  \begin{subfigure}[t]{0.5\linewidth}
    \centering
    \includegraphics[width=0.78\linewidth]{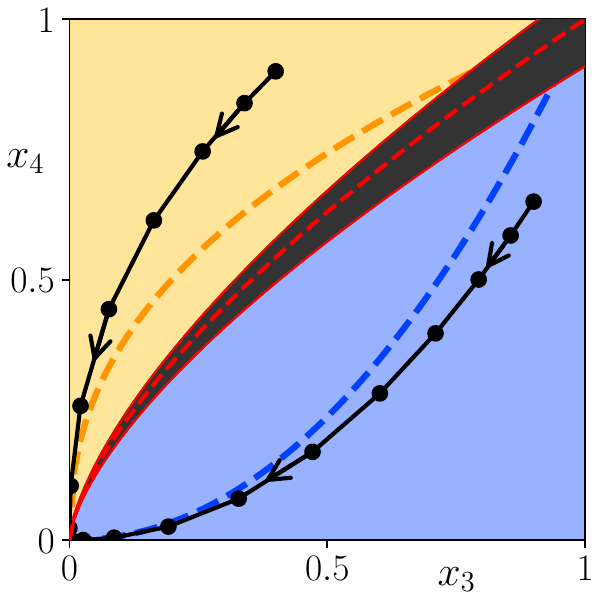}
    \caption{\(\nu_{3}>0\).}
    \label{fig:ks_maps_normal_coords:subfig:pre_BCB}
  \end{subfigure}%
  \begin{subfigure}[t]{0.5\linewidth}
    \centering
    \includegraphics[width=0.78\linewidth]{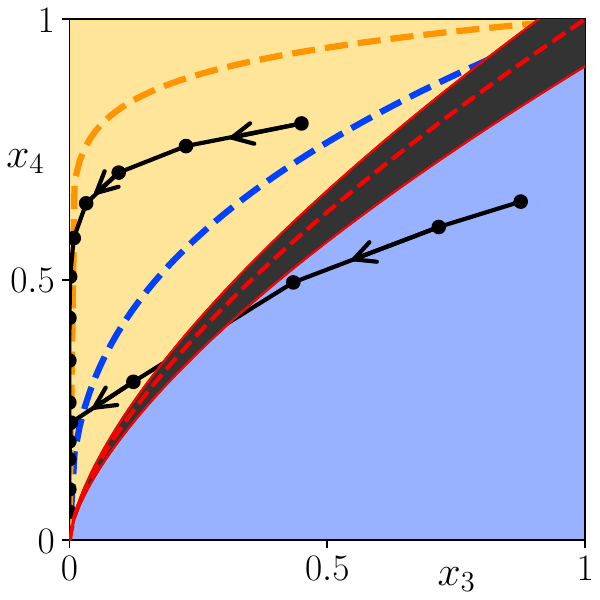}
    \caption{\(\nu_{3}<0\).}
    \label{fig:ks_maps_normal_coords:subfig:post_BCB}
  \end{subfigure}
  \caption{Examples of the action of the return map \(\Phi\colon\poinin[1]{2}\to\poinin[1]{2}\). Points in the same orbit have been joined for clarity only. The domain \(\Gamma_{3}\) is shaded blue, \(\Gamma_{4}\) orange, and the excluded cusp \(\Gamma_{c}\) black. The dashed blue and orange lines are the curves \(\Sigma_{3}^{*}\) and \(\Sigma_{4}^{*}\) (defined in \cref{sec:ret_proj_map:ssec:explanation}). The switching curve \(\Sigma_{s}\) is a dashed red line, and \(\Sigma_{c}^{-}\) and \(\Sigma_{c}^{+}\) are solid red lines.}
  \label{fig:ks_maps_normal_coords}
\end{figure}

%% file: figure_3.tex
\begin{figure}
  \centering
  \begin{subfigure}[t]{0.5\linewidth}
    \centering
    \includegraphics[width=0.79\linewidth]{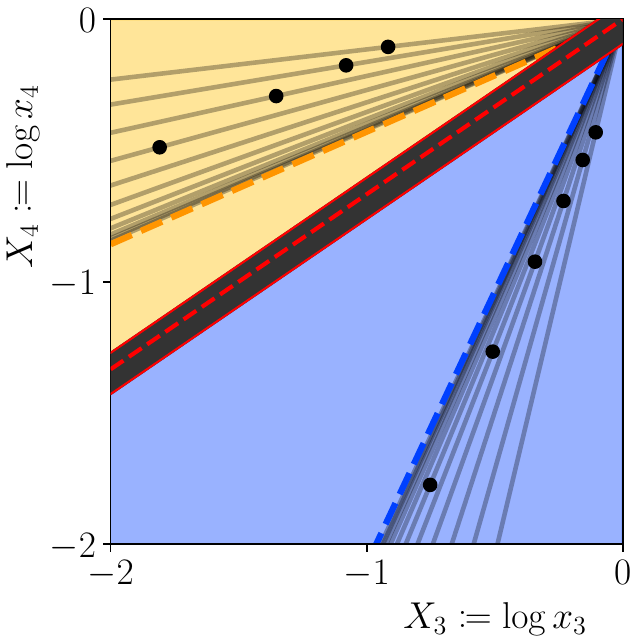}
    \caption{\(\nu_{3}>0\).}
    \label{fig:ks_maps_log_coords:subfig:pre_BCB}
  \end{subfigure}%
  \begin{subfigure}[t]{0.5\linewidth}
    \centering
    \includegraphics[width=0.79\linewidth]{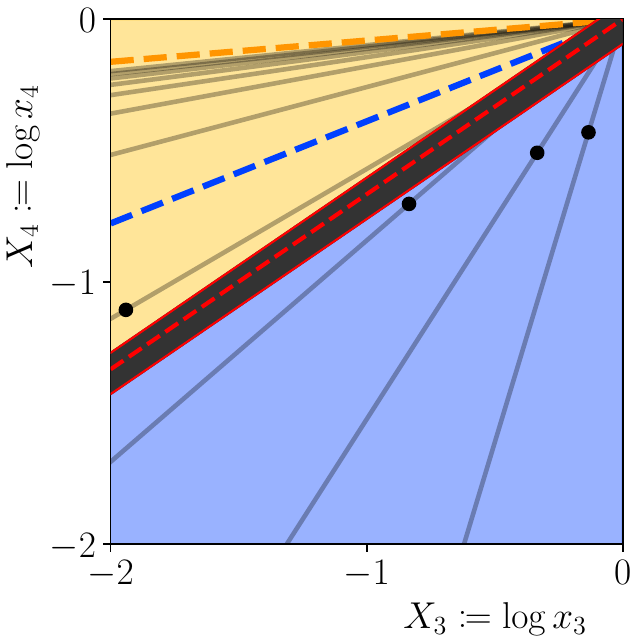}
    \caption{\(\nu_{3}<0\).}
    \label{fig:ks_maps_log_coords:subfig:post_BCB}
  \end{subfigure}
  \caption{Examples of the action of the map \(M\colon\negativeR{2}\to\negativeR{2}\), where vectors extend from the origin to the black dots, and the linear subspace containing each vector is a solid grey line. The domain of \(\mathcal{D}_{3}\) is shaded blue, \(\mathcal{D}_{4}\) orange, and \(\mathcal{D}_{c}\) black. The dashed blue and orange lines are the eigenspaces \(W_{3}^{*}\) and \(W_{4}^{*}\). The vectors can be seen to converge to these eigenspaces under iteration of \(M\). The switching subspace \(W_{s}\) is a dashed red line, and the two affine subspaces \(W_{c}^{-}\) and \(W_{c}^{+}\) are solid red lines. For clarity, \subref{fig:ks_maps_log_coords:subfig:post_BCB} only shows the logarithm of the coordinates of the trajectory in \Cref{fig:ks_maps_normal_coords:subfig:post_BCB} that switches from the \(\cyc{[3]}\) to \(\cyc{[4]}\) cycle.}
  \label{fig:ks_maps_log_coords}
\end{figure}

%% file: figure_4.tex
%
\begin{figure}
  \centering
  \includegraphics[width=0.36\linewidth]{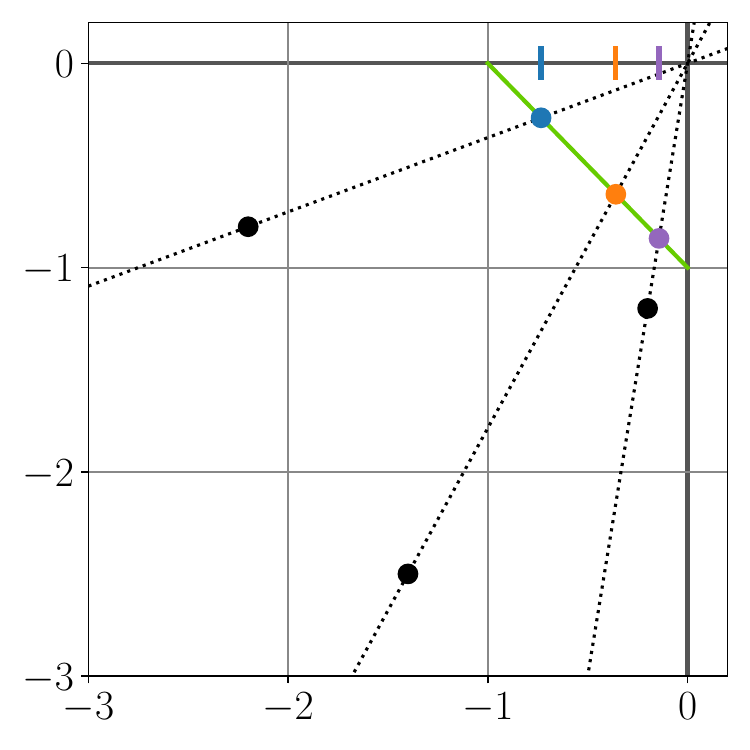}
  \caption{Three examples of determining a vector's projection onto the line segment \(S\), which is shown as a green line. The three black points represent vectors, and the linear subspaces containing them are the dotted lines. The intersection of each linear subspace and \(S\) as a line segment is given by a coloured point, and the horizontal components are the three correspondingly coloured ticks on the \(X_{3}\)-axis}
  \label{fig:projection}
\end{figure}

%% file: 4_analysis.tex
\section{Analysis of projected maps}\label{sec:anal}

In this section, we use the projected maps of the Kirk--Silber, \(\Delta\)-clique, and tournament networks to classify the dynamics of trajectories near each network.

\input{4_1_kirk_silber.tex}
\input{4_2_delta_clique.tex}

\input{4_3_tournament.tex}

%% file: 4_1_kirk_silber.tex
\subsection{The Kirk--Silber network}\label{sec:anal:ssec:ks}

We first define the following two subsets of \(\poinin[1]{2}\):
\begin{align*}
  \Gamma_{3}^{-}&=\left\{\left(x_{3},x_{4}\right)\in\Gamma_{3}\mid\left(1-\epsilon\right)x_{3}^{\frac{e_{24}}{e_{23}}}>x_{4}>x_{3}^{\frac{\rho_{3}}{\delta_{3}-1}}\right\}\subseteq\Gamma_{3},\textrm{ if }\delta_{3}<0\textrm{ and }\nu_{3}<0;\textrm{ and }\\
  \Gamma_{4}^{-}&=\left\{\left(x_{3},x_{4}\right)\in\Gamma_{4}\mid\left(1-\epsilon\right)x_{4}^{\frac{e_{23}}{e_{24}}}>x_{3}>x_{4}^{\frac{\rho_{4}}{\delta_{4}-1}}\right\}\subseteq\Gamma_{4},\textrm{ if }\delta_{4}<0\textrm{ and }\nu_{4}<0.
\end{align*}

In the remainder of this section, we will use the projected map to classify the dynamics of trajectories near the Kirk--Silber network. In particular, we will prove the following theorem. We state these results for only the \(\cyc{[3]}\) cycle, and the equivalent results for the \(\cyc{[4]}\) cycle can be stated by swapping the numbers \(3\) and \(4\) in each of these statements. Our proof focuses on \(f_{3}\), though we highlight where properties of \(f_{4}\) are different.
\begin{theorem}\label{thm:ks}
  \begin{enumerate}[label=\normalfont({\alph*})]
    \item
    If \(\delta_{3}>1\) and \(\nu_{3}>0\), then the \(\cyc{[3]}\) cycle is f.a.s. All trajectories starting in \(\Gamma_{3}\) are asymptotic to the \(\cyc{[3]}\) cycle, with itinerary \((231)^{\infty}\). \label{thm:ks:cyc3:fas}

    \item
    If \(\delta_{3}>1\) and \(\nu_{3}<0\), then the \(\cyc{[3]}\) cycle is a.c.u. All trajectories starting in \(\Gamma_{3}\) are eventually mapped into \(\Gamma_{4}\). If \(\delta_{4}>1\), these trajectories switch from the \(\cyc{[3]}\) cycle to the \(\cyc{[4]}\) cycle, to which they are then asymptotic. For all positive integers \(n\), there exists an open set of trajectories starting in \(\Gamma_{3}\) with the itinerary \((231)^{n}(241)^{\infty}\). \label{thm:ks:cyc3:switch}

    \item
    If \(\delta_{3}<1\), then the \(\cyc{[3]}\) cycle is c.u. If \(\nu_{3}<0\), then all trajectories starting in \(\Gamma_{3}^{-}\) sufficiently close to the network are eventually mapped into \(\Gamma_{4}\). If \(\delta_{4}>1\), these trajectories switch from the \(\cyc{[3]}\) cycle to the \(\cyc{[4]}\) cycle, to which they are then asymptotic. For all positive integers \(n\), there exists an open set of trajectories in \(\Gamma_{3}\) with the itinerary \((231)^{n}(241)^{\infty}\). All other trajectories starting near the cycle leave a neighbourhood of the cycle. \label{thm:ks:cyc3:res}
  \end{enumerate}
\end{theorem}

The results in \ref{thm:ks:cyc3:fas} and \ref{thm:ks:cyc3:switch} were proven in \cite{kirk_silber_1994}, though not using the terminology we have used here. We rederive these results to demonstrate how trajectories near a heteroclinic network can be analysed with its projected map. We also derive new results in \ref{thm:ks:cyc3:res} concerning switching from a c.u. cycle.

\input{figure_5.tex}
\subsubsection{Fixed points and admissibility}\label{sec:anal:ssec:ks:sssec:fp}

Both \(f_{3}\) and \(f_{4}\) have only one fixed point that can lie in \(S\),
\begin{equation*}
  \vartheta_{3}^{*}=\frac{1-\delta_{3}}{\delta_{3}+\rho_{3}-1}
  \qquad\textrm{ and }\qquad
  \vartheta_{4}^{*}=\frac{-\rho_{4}}{\delta_{4}+\rho_{4}-1}.
\end{equation*}the
These fixed points correspond to the eigenspaces \(W_{3}^{*}\) and \(W_{4}^{*}\).

Since \(f_{3}\) and \(f_{4}\) are only defined on a subset of \(S\), formal solutions to the equations \(f_{3}\left(\vartheta\right)=\vartheta\) and \(f_{4}\left(\vartheta\right)=\vartheta\) may exist outside \(\Theta_{3}\) or \(\Theta_{4}\), respectively. Consider, for example, \(\vartheta_{3}^{*}\). If \(\vartheta_{3}^{*}\in\Theta_{3}\), we say it is \textit{admissible}, and if \(\vartheta_{3}^{*}\in\Theta_{4}\), we say it is \textit{virtual}. We generically say that a \textit{border-collision bifurcation} occurs when \(\vartheta_{3}^{*}=\vartheta_{s}\) and so \(\vartheta_{3}^{*}\) passes from \(\Theta_{3}\) to \(\Theta_{4}\). If \(\vartheta_{3}^{*}\notin S\), we also say \(\vartheta_{3}^{*}\) is \textit{virtual}. We now determine when these fixed points are admissible or virtual, and when border-collision bifurcations occur.

We assume that \(\delta_{3}+\rho_{3}>1\) and \(\delta_{4}+\rho_{4}>1\). Solving \(\vartheta_{3}^{*}>\vartheta_{s}\) gives
\begin{equation*}
  \rho_{3}-\frac{e_{24}}{e_{23}}\left(\delta_{3}-1\right)\equiv\frac{c_{21}}{e_{23}}\nu_{3}>0.
\end{equation*}
Since \(c_{21}\) and \(e_{23}\) are both positive, \(\vartheta_{3}^{*}\) is admissible if and only if \(\nu_{3}>0\) and it is virtual if and only if \(\nu_{3}<0\). A border-collision bifurcation occurs at \(\nu_{3}=0\). Equivalent results hold for \(\vartheta_{4}^{*}\) and \(\nu_{4}\).

In the case that \(\delta_{3}+\rho_{3}<1\), \(\vartheta_{3}^{*}\) is virtual if and only if \(\nu_{3}>0\), and admissible if and only if \(\nu_{3}<0\). The same holds for \(\vartheta_{4}^{*}\) when \(\delta_{4}+\rho_{4}<1\).

The maps \(f_{3}\) and \(f_{4}\) each have one additional fixed point: \(\vartheta_{3}^{-}=0\) and \(\vartheta_{4}^{-}=-1\), respectively. These correspond to the eigenspaces \(W_{3}^{-}\) and \(W_{4}^{-}\), respectively. These points are not in the domains of the maps \(f_{3}\) and \(f_{4}\), but it will be useful to consider them in our analysis.

\subsubsection{Stability}\label{sec:anal:ssec:kirk_silber:sssec:stab}

To determine the stability of \(\vartheta_{3}^{*}\) and \(\vartheta_{4}^{*}\), we calculate that
\begin{equation*}
  Df_{3}(\vartheta)=\frac{\delta_{3}}{\left(\left(\delta_{3}+\rho_{3}-1\right)\vartheta-1\right)^{2}}
  \qquad\textrm{ and }\qquad
  Df_{4}(\vartheta)=\frac{\delta_{4}}{\left(\left(\delta_{4}+\rho_{4}-1\right)\vartheta+\delta_{4}+\rho_{4}\right)^2}.
\end{equation*}
Since \(\delta_{3}\) and \(\delta_{4}\) are both positive, \(Df_{3}(\vartheta)>0\) and \(Df_{4}(\vartheta)>0\) for all \(\vartheta\in S\), and so both functions are strictly increasing.

Evaluating \(Df_{3}(\vartheta)\) and \(Df_{4}(\vartheta)\) at \(\vartheta_{3}^{*}\) and \(\vartheta_{4}^{*}\) gives
\begin{equation*}
  Df_{3}(\vartheta_{3}^{*})=\frac{1}{\delta_{3}}
  \qquad\textrm{ and }\qquad
  Df_{4}(\vartheta_{4}^{*})=\frac{1}{\delta_{4}}.
\end{equation*}
Therefore, \(\vartheta_{3}^{*}\) is asymptotically stable when \(\delta_{3}>1\), and \(\vartheta_{4}^{*}\) is asymptotically stable when \(\delta_{4}>1\). At the point of stability loss, \(\vartheta_{3}^{*}=\vartheta_{3}^{-}\) and \(\vartheta_{4}^{*}=\vartheta_{4}^{-}\). Additionally, \(Df_{3}(\vartheta_{3}^{-})=\delta_{3}\) and \(Df_{4}(\vartheta_{4}^{-})=\delta_{4}\). Thus, if \(f_{3}\) and \(f_{4}\) were studied as functions of the real line, these bifurcations would be transcritical bifurcations.

\subsubsection{Equivalence of the border-collision bifurcation and Podvigina's third stability condition}\label{sec:anal:ssec:kirk_silber:sssec:BCB_PodIII_equiv}

When studying a heteroclinic cycle, it may be that there is only one equilibrium at which an open set of initial conditions exists that does not return to that equilibrium under the flow of the system. As such, instability may only be evident at one equilibrium of the cycle. Thus, return maps to cross-sections defined near every equilibrium need to be considered. Since the projected map \(f\) is derived from the transition matrices \(M_{3}\) and \(M_{4}\) only, we do not know \textit{a priori} that stability of the cycles cannot necessarily be determined by studying this map alone. However, it is in fact possible to establish an equivalence between the border-collision bifurcation and stability conditions at \(\xi_{3}\) and \(\xi_{4}\); this equivalence allows us to deduce if a trajectory is asymptotic to the \(\cyc{[3]}\) or \(\cyc{[4]}\) cycle from examination of the projected map only.

If \(w\) is an eigenvector of \(M_{3}\), then \(m_{123}w\) is an eigenvector of \(M_{3}^{(3)}\), and, if \(w\) is an eigenvector of \(M_{3}^{(3)}\), then \(m_{123}^{-1}w\) is an eigenvector of \(M_{3}\). The key observation is that the image of the switching vector \(w_{s}=\left(1,\frac{e_{24}}{e_{23}}\right)\) under \(m_{123}\) is
\begin{equation*}
  m_{123}w_{s}=\begin{pmatrix*}[r]
    \frac{c_{21}}{e_{23}} & 0 \\[0.5em]
    -\frac{e_{24}}{e_{23}} & 1
  \end{pmatrix*}\begin{pmatrix}
    1 \\[0.5em]
    \frac{e_{24}}{e_{23}}
  \end{pmatrix}=\begin{pmatrix}
    \frac{c_{21}}{e_{23}} \\[0.5em]
    0
  \end{pmatrix}.
\end{equation*}
Therefore, if \(w_{\max}\) of \(M_{3}\) is in \(W_{s}\), \(w_{\max}\) of \(M_{3}^{(3)}\) will not satisfy condition \ref{cond:podIII} of \Cref{thm:pod_bif}. Similarly, at the bifurcation \(\nu_{3}=0\), \(w_{\max}\) of \(M_{3}^{(3)}\) is \((\delta_{3}-1,0)\), and we find that
\begin{equation*}
  m_{123}^{-1}w_{\max}=\begin{pmatrix*}[r]
    \frac{e_{23}}{c_{21}} & 0 \\[0.5em]
    \frac{e_{24}}{c_{21}} & 1
  \end{pmatrix*}\begin{pmatrix}
    \delta_{3}-1 \\[0.5em] 0
  \end{pmatrix}=\frac{e_{23}}{c_{21}}\left(\delta_{3}-1\right)\begin{pmatrix}
    1 \\[0.5em] \frac{e_{24}}{e_{23}}
  \end{pmatrix}\in W_{s}.
\end{equation*}
With a similar argument for the matrix \(m_{124}\), we have proved the following proposition.
\begin{proposition}\label{prop:ks:bcb_podIII_equiv}
  Suppose \(\delta_{3}>1\). The fixed point \(\vartheta_{3}^{*}\) is admissible if and only if \(w_{\max}\) of \(M_{3}^{(3)}\) satisfies Podvigina's third condition for stability. Suppose \(\delta_{4}>1\). The fixed point \(\vartheta_{4}^{*}\) is admissible if and only if \(w_{\max}\) of \(M_{4}^{(4)}\) satisfies Podvigina's third condition for stability.
\end{proposition}
We can then conclude that, if \(\left(f^{n}(\vartheta)\right)\) is an orbit of \(f\) that is asymptotic to an admissible \(\vartheta_{3}^{*}\), then all corresponding trajectories near the Kirk--Silber network are asymptotic to the origin in \textit{all} cross-sections of the \(\cyc{[3]}\) cycle. A similar argument applies to the \(\cyc{[4]}\) cycle.

\subsubsection{Discontinuity of the projected map of the Kirk--Silber network at \texorpdfstring{\(\vartheta_{s}\)}{theta s}}

We now consider the limiting behaviour of \(f_{3}\) and \(f_{4}\) as \(\vartheta\to\vartheta_{s}\). The projected map is not defined at \(\vartheta_{s}\) because of the method used to construct the return map \(\Phi\). The functions \(f_{3}\) and \(f_{4}\) can nevertheless be formally evaluated at \(\vartheta_{s}\).

We consider the projected map \(f\) to be continuous if
\begin{equation*}
  \lim_{\vartheta\searrow\vartheta_{s}}f_{3}(\vartheta)=\lim_{\vartheta\nearrow\vartheta_{s}}f_{4}(\vartheta).
\end{equation*}
In the case of the Kirk--Silber network, we find that
\begin{equation*}
  \lim_{\vartheta\searrow\vartheta_{s}}f_{3}(\vartheta)=\frac{-c_{13}c_{32}}{c_{32}\left(c_{13}+c_{14}\right)+e_{12}t_{34}},
  \qquad\textrm{ and }\qquad
  \lim_{\vartheta\nearrow\vartheta_{s}}f_{4}(\vartheta)=\frac{-\left(c_{13}c_{32} + e_{12}t_{43}\right)}{c_{42}\left(c_{13}+c_{14}\right)+e_{12}t_{43}},
\end{equation*}
which, generically, are not equal. Therefore, the projected map of the Kirk--Silber network is discontinuous. This discontinuity can be seen in all three examples in \Cref{fig:ks_proj_map}.

\subsubsection{Dynamics when \texorpdfstring{\(\nu_{3}>0\)}{v3>0} and \texorpdfstring{\(\nu_{4}>0\)}{v4>0}}\label{sec:anal:ssec:kirk_silber:sssec:admiss}

We first consider the dynamics of \(f\) when both \(\vartheta_{3}^{*}\) and \(\vartheta_{4}^{*}\) are admissible.

We first assume \(\delta_{3}>1\). For all \(\vartheta\in\Theta_{3}\), \(f(\vartheta)\in\Theta_{3}\). If \(\vartheta>\vartheta_{3}^{*}\), then \(f(\vartheta)<\vartheta\), and if \(\vartheta<\vartheta_{3}^{*}\), then \(f(\vartheta)>\vartheta\). Hence, \(\vartheta_{3}^{*}\) is globally attracting in \(\Theta_{3}\) if \(\delta_{3}>1\). A similar analysis shows that \(\vartheta_{4}^{*}\) is globally attracting in \(\Theta_{4}\) if \(\delta_{4}>1\). Hence, we have proven \ref{thm:ks:cyc3:fas} of \Cref{thm:ks}. Examples of orbits under iteration of the projected map for \(\delta_{3},\delta_{4}>1\) and \(\nu_{3},\nu_{4}>0\) are shown in \Cref{fig:ks_proj_map:subfig:pre_bcb}.

When \(\delta_{3}<1\), \(\vartheta_{3}^{*}>\vartheta_{3}^{-}=0\), and \(f_{3}\) has no fixed points in \(\Theta_{3}\). However, for all \(\vartheta\in\Theta_{3}\), \(f(\vartheta)\in\Theta_{3}\) and \(f(\vartheta)>\vartheta\). Therefore, \(\lim_{n\to\infty}f^{n}(\vartheta)=\vartheta_{3}^{-}=0\), and so there is no \(\vartheta\in\Theta_{3}\) asymptotic to a point in \(S\). The \(\cyc{[3]}\) cycle is thus c.u. Hence, we have proven the first part of \ref{thm:ks:cyc3:res} of \Cref{thm:ks}. Examples of orbits under iteration of the projected map when \(\delta_{3}<1\) are shown in \Cref{fig:ks_proj_map:subfig:post_res}.

\subsubsection{Dynamics when \texorpdfstring{\(\nu_{3}<0\)}{v3<0} or \texorpdfstring{\(\nu_{4}<0\)}{v4<0}}\label{sec:anal:ssec:kirk_silber:sssec:virt}

As noted in section \ref{sec:ret_proj_map:ssec:ret_maps}, Kirk and Silber show that these values cannot simultaneously be negative. We will analyse the dynamics of \(f\) when \(\nu_{3}<0\), and the analysis in the case of \(\nu_{4}<0\) is similar.

We first assume that \(\delta_{3}>1\). Since \(\nu_{3}<0\), \(\vartheta_{3}^{*}\) is virtual, and \(\vartheta_{3}^{*}<\vartheta_{s}\). For all \(\vartheta\in\Theta_{3}\), \(f(\vartheta)<\vartheta\). However, there now exist \(\vartheta\in\Theta_{3}\) such that \(f(\vartheta)\notin\Theta_{3}\). For all non-negative integers \(n\), define \(\mathcal{E}_{n}\) as the \(n\)th pre-image of \(\vartheta_{s}\) under \(f_{3}\); in particular, \(f(\mathcal{E}_{n+1})=\mathcal{E}_{n}\) and \(f^{n}(\mathcal{E}_{n})=\vartheta_{s}\). The exact value of \(\mathcal{E}_{n}\) is
\begin{equation}\label{eqn:theta_s_pre_image}
  \mathcal{E}_{n}=\frac{\vartheta_{s}}{\delta_{3}^{n}+\vartheta_{s}\left(\delta_{3}^{n}-1+\rho_{3}\sum\limits_{k=0}^{n-1}\delta_{3}^{k}\right)}.
\end{equation}
Note that \(\mathcal{E}_{0}=\vartheta_{s}\), and, for all \(n\), \(\mathcal{E}_{n}<\mathcal{E}_{n+1}<0\), and \(\lim_{n\to\infty}\mathcal{E}_{n}=\vartheta_{3}^{-}=0\).

For all \(\vartheta\in\left(\mathcal{E}_{0},\mathcal{E}_{1}\right)\), \(f(\vartheta)\in\Theta_{4}\). Moreover, for all \(\vartheta\in\left(\mathcal{E}_{n},\mathcal{E}_{n+1}\right)\), \(f(\vartheta)\in\left(\mathcal{E}_{n-1},\mathcal{E}_{n}\right)\) and so \(f^{n+1}(\vartheta)\in\Theta_{4}\). Therefore, all points in \(\Theta_{3}\) that are not a pre-image of \(\vartheta_{s}\) are eventually mapped into \(\Theta_{4}\) after a finite number of iterations.

If \(\delta_{4}>1\), \(\vartheta_{4}^{*}\) is globally attracting in \(\Theta_{4}\), and so all \(\vartheta\in\Theta_{3}\), except the points \(\mathcal{E}_{n}\), are asymptotic to \(\vartheta_{4}^{*}\). However, by \cite[Theorem~2.4]{krupa_melbourne_2004}, the \(\cyc{[3]}\) cycle remains asymptotically stable in the subspace \(S_{123}\) defined by \(x_{4}=0\), and so \(\cyc{[3]}\) is a.c.u. Hence, we have completed the proof of \ref{thm:ks:cyc3:switch} of \Cref{thm:ks}. Examples of orbits under iteration of the projected map that switch from \(\Theta_{3}\) to \(\Theta_{4}\) are shown in \Cref{fig:ks_proj_map:subfig:post_bcb}.

The same analysis can be done in the case of \(\nu_{4}<0\), showing that all points in \(\Theta_{4}\), except those that are a pre-image of \(\vartheta_{s}\), are eventually mapped into \(\Theta_{3}\).

If \(\delta_{3}<1\), then since \(\nu_{3}<0\), \(\vartheta_{3}^{*}\) is admissible but repelling in \(\Theta_{3}\). As such, all \(\vartheta>\vartheta_{3}^{*}\) are asymptotic to \(\vartheta_{3}^{-}\), as was the case when \(\nu_{3}>0\). However, \(\lim_{n\to\infty}\mathcal{E}_{n}=\vartheta_{3}^{*}\), for the \(\mathcal{E}_{n}\) defined in \cref{eqn:theta_s_pre_image}. Thus, all \(\vartheta\in(\vartheta_{s},\vartheta_{3}^{*})\) are mapped into \(\Theta_{4}\) after a finite number of iterations. The subset \(\vartheta\in(\vartheta_{s},\vartheta_{3}^{*})\subseteq\Theta_{3}\) is the projection of the subset \(\Gamma_{3}^{-}\subseteq\poinin[1]{2}\). We have therefore completed the proof of \ref{thm:ks:cyc3:res} of \Cref{thm:ks}. An example of the projected map with \(\delta_{3}<1\) and \(\nu_{3}<0\) is shown in \Cref{fig:ks_proj_map:subfig:post_res_bcb}.

Hence, we have completed the proof of \Cref{thm:ks}.

\subsubsection{Numerical examples}

We conclude our study of the Kirk--Silber network by presenting some numerical examples of the phenomena we have described in this section. We consider the system of ODEs
\begin{equation}\label{eqn:ks_dyn_sys}
  \begin{aligned}
    \dot{x}_{1}&=x_{1}\left(1-\left\lVert x\right\rVert_{2}^{2}-c_{21}x_{2}^{2}+e_{31}x_{3}^{2}+e_{41}x_{4}^{2}\right),\\
    \dot{x}_{2}&=x_{2}\left(1-\left\lVert x\right\rVert_{2}^{2}+e_{12}x_{1}^{2}-c_{32}x_{3}^{2}-c_{42}x_{4}^{2}\right),\\
    \dot{x}_{3}&=x_{3}\left(1-\left\lVert x\right\rVert_{2}^{2}-c_{13}x_{1}^{2}+e_{23}x_{2}^{2}-t_{43}x_{4}^{2}\right),\\
    \dot{x}_{4}&=x_{4}\left(1-\left\lVert x\right\rVert_{2}^{2}-c_{14}x_{1}^{2}+e_{24}x_{2}^{2}-t_{34}x_{3}^{2}\right),
  \end{aligned}
\end{equation}
where \(c_{jk}\), \(e_{jk}\), and \(t_{jk}\) are positive real numbers. System \cref{eqn:ks_dyn_sys} is \(\Z_{2}^{4}\)-equivariant, and contains a Kirk--Silber network. The four equilibria of the network lie at unit distance from the origin on the four coordinate axes.

To minimise roundoff error when computing \(f\) for  trajectories close to the network, we follow \cite{postlethwaite_dawes_2005} and use the change of coordinates \(X_{j}\coloneq\log x_{j}\), and integrate the resulting ODEs \(\dot{X}_{j}\) with a variety of parameter values. \Cref{fig:ks_examples} shows three different trajectories near the Kirk--Silber network. The timeseries in \Cref{fig:ks_examples:C3_cycle} is an example of a trajectory asymptotic to the \(\cyc{[3]}\) cycle. The example in \Cref{fig:ks_examples:switch} has the same initial conditions as \Cref{fig:ks_examples:C3_cycle}, but the parameters are now such that \(\nu_{3}<0\), and the trajectory switches from the a.c.u. \(\cyc{[3]}\) cycle to the \(\cyc{[4]}\) cycle. This trajectory has itinerary \((123)^{2}(124)^{\infty}\). \Cref{fig:ks_examples:res_switch} shows an example of switching from a c.u. \(\cyc{[3]}\) cycle to a stable \(\cyc{[4]}\) cycle, with itinerary \((231)^{3}(241)^{\infty}\).

\input{figure_6.tex}

%% file: figure_5.tex
\begin{figure}
  \centering
  \begin{subfigure}[t]{0.5\linewidth}
    \centering
    \includegraphics[width=\linewidth]{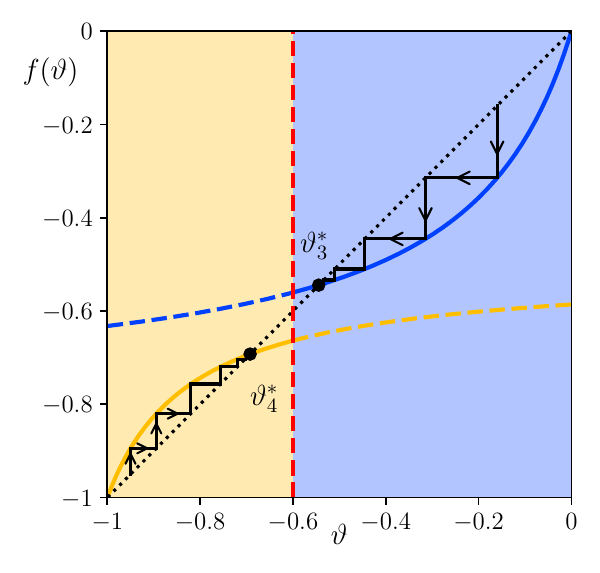}
    \vspace{-8mm}
    \caption{\(\delta_{3}>1\) and \(\nu_{3}>0\)}
    \label{fig:ks_proj_map:subfig:pre_bcb}
  \end{subfigure}%
  \begin{subfigure}[t]{0.5\linewidth}
    \centering
    \includegraphics[width=\linewidth]{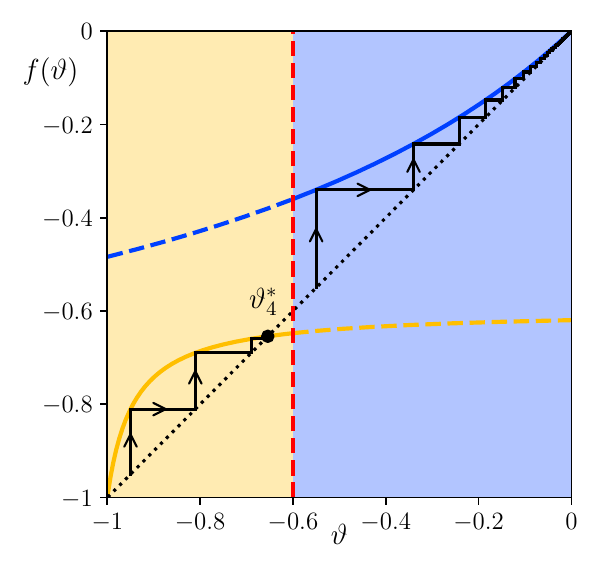}
    \vspace{-8mm}
    \caption{\(\delta_{3}<1\) and \(\nu_{3}>0\)}
    \label{fig:ks_proj_map:subfig:post_res}
  \end{subfigure}
  \begin{subfigure}[t]{0.5\linewidth}
    \centering
    \includegraphics[width=\linewidth]{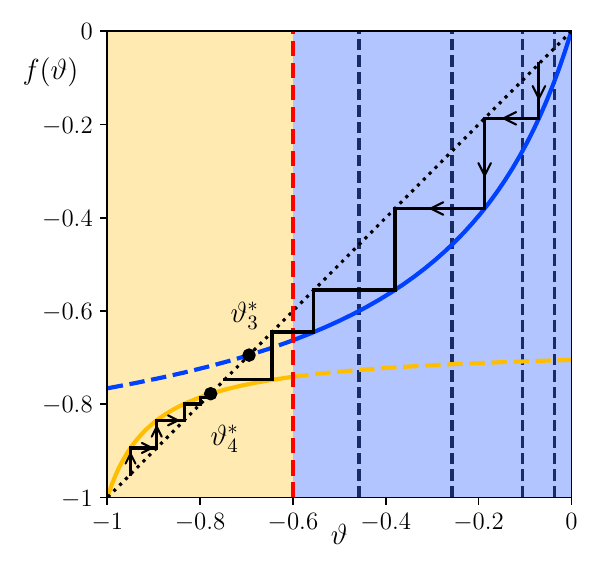}
    \vspace{-8mm}
    \caption{\(\delta_{3}>1\) and \(\nu_{3}<0\).}
    \label{fig:ks_proj_map:subfig:post_bcb}
  \end{subfigure}%
  \begin{subfigure}[t]{0.5\linewidth}
    \centering
    \includegraphics[width=\linewidth]{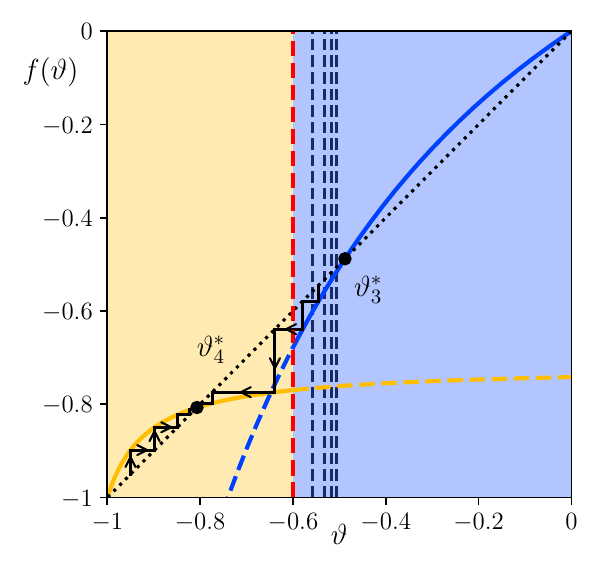}
    \vspace{-8mm}
    \caption{\(\delta_{3}<1\) and \(\nu_{3}<0\).}
    \label{fig:ks_proj_map:subfig:post_res_bcb}
  \end{subfigure}
  \caption{Schematic representations of the projected map of the Kirk--Silber network as cobweb plots, for different relations between parameters, as indicated. In all examples, the domain \(\Theta_{3}\) is shaded blue and the domain \(\Theta_{4}\) orange. The dashed red line is \(\vartheta_{s}\). The dotted black line is \(f\left(\vartheta\right)=\vartheta\). Each function \(f_{3}\) and \(f_{4}\) is shown as solid in its domain of definition and dashed otherwise. In all four figures, \(\delta_{4}>1\) and \(\nu_{4}>0\). In \subref{fig:ks_proj_map:subfig:post_bcb} and \subref{fig:ks_proj_map:subfig:post_res_bcb}, the values \(\mathcal{E}_{n}\) (see \cref{eqn:theta_s_pre_image}) are indicated by dashed black lines for \(n=1,2,3,4\).}
  \label{fig:ks_proj_map}
\end{figure}

%% file: figure_6.tex
\begin{figure}
  \begin{subfigure}[c]{\linewidth}
    \includegraphics[width=0.46\linewidth]{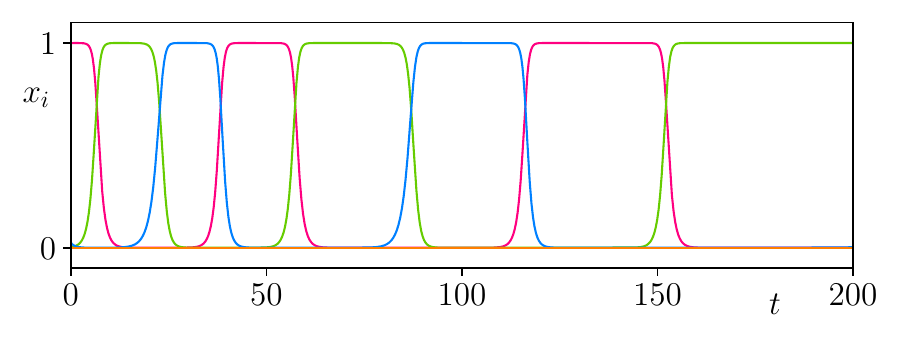}%
    \includegraphics[width=0.53\linewidth]{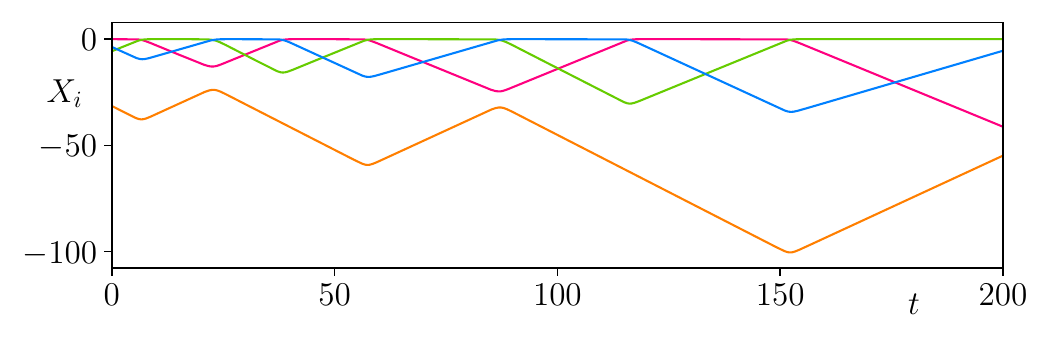}
    \vspace{-3mm}
    \caption{\(\delta_{3}>1\) and \(\nu_{3}>0\)}
    \label{fig:ks_examples:C3_cycle}
  \end{subfigure}
  \begin{subfigure}[c]{\linewidth}
    \includegraphics[width=0.46\linewidth]{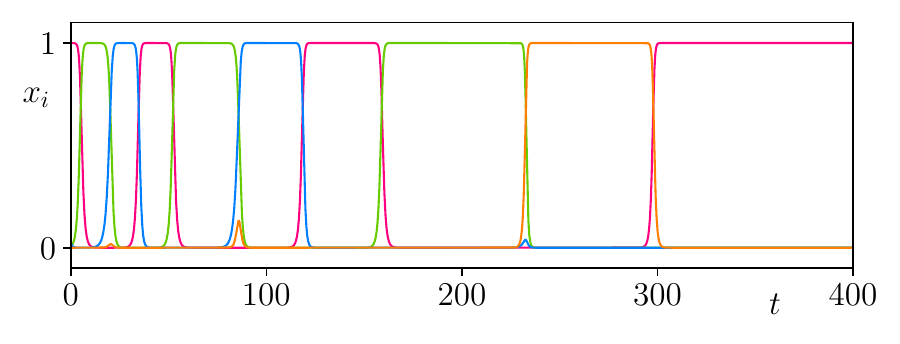}%
    \includegraphics[width=0.53\linewidth]{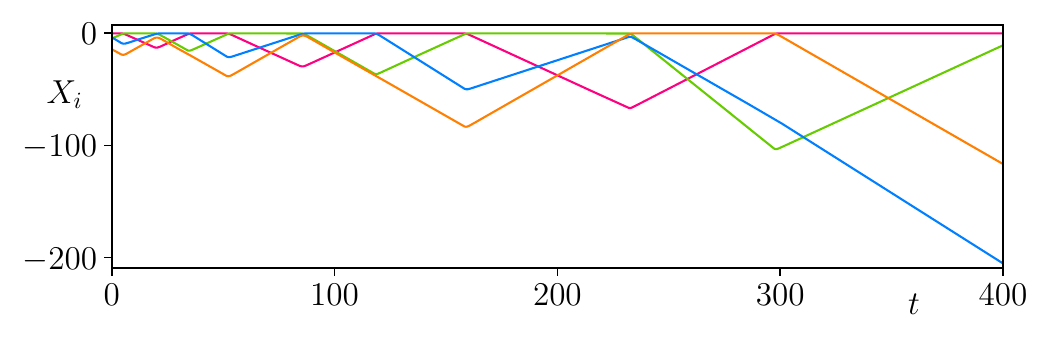}
    \vspace{-3mm}
    \caption{\(\delta_{3}>1\) and \(\nu_{3}<0\)}
    \label{fig:ks_examples:switch}
  \end{subfigure}
  \begin{subfigure}[c]{\linewidth}
    \includegraphics[width=0.46\linewidth]{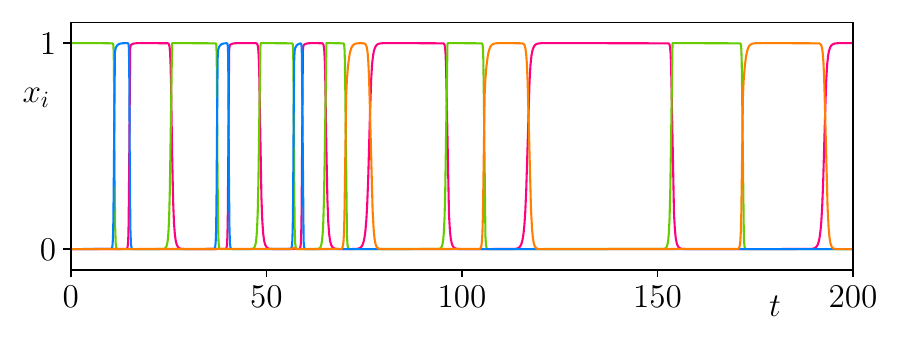}%
    \includegraphics[width=0.53\linewidth]{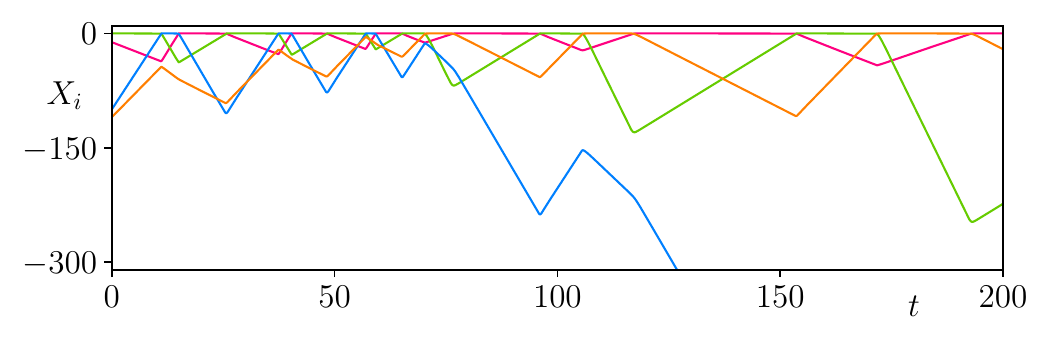}
    \vspace{-3mm}
    \caption{\(\delta_{3}<1\) and \(\nu_{3}<0\)}
    \label{fig:ks_examples:res_switch}
  \end{subfigure}
  \caption{Example timeseries of trajectories near the Kirk--Silber network. The colours red, blue, green, and orange correspond to \(x_{1}\), \(x_{2}\), \(x_{3}\), and \(x_{4}\), respectively, and similarly for the \(X_{i}\). In all three examples, \(\delta_{4}>1\) and \(\nu_{4}>0\), and the values of \(\delta_{3}\) and \(\nu_{3}\) are as indicated.}
  \label{fig:ks_examples}
\end{figure}

%% file: 4_2_delta_clique.tex
\subsection{The \texorpdfstring{\(\Delta\)}{Delta}-clique network}\label{sec:anal:ssec:delta_clique}

We now study the \(\Delta\)-clique network using its projected map. Following the process outlined in \cref{sec:ret_proj_map:ssec:proj_map}, we derive the projected map of the \(\Delta\)-clique network, \(f\colon\Theta_{34}\cup\Theta_{4}\to S\):
\begin{equation}\label{eqn:delta_clique_proj_map}
  f(\vartheta)=\begin{cases}
    f_{34}(\vartheta)=\dfrac{\left(\alpha_{2}-\alpha_{1}\right)\vartheta + \alpha_{2}}{\left(\alpha_{1}-\alpha_{2}+\alpha_{3}-\alpha_{4}\right)\vartheta-\alpha_{2}-\alpha_{4}}, & \textrm{if } \vartheta\in\Theta_{34},\\[1.46em]
    f_{4}(\vartheta)=\dfrac{\left(1-\rho_{4}\right)\vartheta-\rho_{4}}{\left(\delta_{4}+\rho_{4}-1\right)\vartheta+\delta_{4}+\rho_{4}}, & \textrm{if } \vartheta\in\Theta_{4},
  \end{cases}
\end{equation}
where
\begin{align*}
  \alpha_{1}&=\frac{c_{21}c_{43}}{e_{23}e_{41}} + \frac{c_{21}c_{42}t_{13}}{e_{23}e_{41}e_{12}} - \frac{c_{32}e_{24}t_{13}}{e_{34}e_{23}e_{12}} - \frac{c_{43}e_{24}t_{31}}{e_{41}e_{23}e_{34}} - \frac{c_{42}e_{24}t_{13}t_{31}}{e_{41}e_{23}e_{12}e_{34}},\\[0.5em]
  \alpha_{2}&=\frac{c_{32}t_{13}}{e_{34}e_{12}} + \frac{c_{43}t_{31}}{e_{41}e_{34}} + \frac{c_{42}t_{13}t_{31}}{e_{41}e_{12}e_{34}},
  \qquad\alpha_{3}=\frac{c_{14}c_{21}c_{42}}{e_{12}e_{23}e_{41}} - \frac{c_{14}c_{32}e_{24}}{e_{12}e_{34}e_{23}} - \frac{c_{14}c_{42}e_{24}t_{31}}{e_{12}e_{41}e_{23}e_{34}},\\[0.5em]
  \alpha_{4}&=\frac{c_{14}c_{32}}{e_{12}e_{34}} + \frac{c_{14}c_{42}t_{31}}{e_{12}e_{41}e_{34}},\qquad\mathrm{ and }
  \qquad\rho_{4}=-\frac{e_{23}}{e_{24}} + \frac{c_{21}c_{43}}{e_{24}e_{41}} + \frac{c_{21}c_{42}t_{13}}{e_{24}e_{41}e_{12}}.
\end{align*}
The definition of \(\delta_{4}\) is the same as for the Kirk--Silber network. The derivation of these expressions can be found in the Supplementary Material (section \ref{supp:maps:delta}).

\input{figure_7.tex}
The maps \(f_{34}\colon\Theta_{34}\to S\) and \(f_{4}\colon\Theta_{4}\to S\) are the projected maps of the \(\cyc{[34]}\) and \(\cyc{[4]}\) cycles. The map \(f_{4}\) is analogous to that for the Kirk--Silber network, with the exception that \(\rho_{4}\) has a slightly different definition as a result of the different topologies of the Kirk--Silber and \(\Delta\)-clique networks. 

The domain \(\Theta_{34}\) of \(f_{34}\) corresponds to the set \(\Gamma_{34}\subseteq\poinin[1]{2}\) of trajectories that leave a neighbourhood of \(\xi_{2}\) in the direction of \(\xi_{3}\), and therefore cycle around the \(\cyc{[34]}\) cycle. The domain \(\Theta_{4}\) corresponds to those trajectories that cycle around the \(\cyc{[4]}\) cycle. These domains are \(\Theta_{34}=\left(\vartheta_{s},0\right)\) and \(\Theta_{4}=\left(-1,\vartheta_{s}\right)\).

For the purposes of describing the dynamics that can be observed near the \(\Delta\)-clique network, we introduce the following expression for the \(\cyc{[4]}\) cycle:
\begin{equation*}
  \nu_{4}=\frac{c_{43}}{e_{41}}+\frac{c_{42}t_{13}}{e_{41}e_{12}}-\frac{c_{14}c_{42}e_{23}}{e_{12}e_{24}e_{41}},
\end{equation*}
and the following expressions for the \(\cyc{[34]}\) cycle:
\begin{align*}
  \tau_{34}&=\alpha_{1}+\alpha_{4}=\mathrm{tr}M_{34},\qquad\delta_{34}=\alpha_{1}\alpha_{4}-\alpha_{2}\alpha_{3}=\det M_{34}=\frac{c_{14}c_{21}c_{32}c_{43}}{e_{12}e_{23}e_{34}e_{41}}>0,\\
  \omega_{34}&=\tau_{34}^{2}-4\delta_{34},\qquad\lambda_{34}^{\pm}=\frac{1}{2}\left(\tau_{34}\pm\sqrt{\omega_{34}}\right),\qquad\zeta_{34}=\alpha_{1}-\alpha_{2}+\alpha_{3}-\alpha_{4},\quad\mathrm{ and }\\
  \beta_{34}&=\alpha_{3}+\frac{1}{2}\left(\alpha_{1}-\alpha_{4}\right)+\frac{e_{24}}{e_{23}}\left(\alpha_{2}+\frac{1}{2}\left(\alpha_{4}-\alpha_{1}\right)\right).
\end{align*}
The expression for \(\nu_{4}\) here is slightly different to that in \cref{sec:ret_proj_map}, again due to the different topologies of the two networks, but bifurcations when \(\nu_{4}=0\) correspond to the same dynamical phenomena.

In the case that \(\omega_{34},\tau_{34},\nu_{4}>0\) and \(\zeta_{34},\beta_{34}<0\), we also introduce the following two subsets of \(\Gamma_{34}\):
\begin{align*}
  \Gamma_{34}^{-}&=\left\{\left(x_{1},x_{4}\right)\in\Gamma_{4}\mid(1-\epsilon)x_{3}^{\frac{e_{24}}{e_{23}}}>x_{4}>x_{3}^{\frac{\lambda_{34}^{+}-\alpha_{1}-\alpha_{3}}{-\lambda_{34}^{+}+\alpha_{2}+\alpha_{4}}}\right\},\\
  \Gamma_{34}^{+}&=\left\{\left(x_{1},x_{4}\right)\in\Gamma_{4}\mid x_{4}\leq x_{3}^{\frac{\lambda_{34}^{+}-\alpha_{1}-\alpha_{3}}{-\lambda_{34}^{+}+\alpha_{2}+\alpha_{4}}}\right\}=\Gamma_{34}\setminus\Gamma_{34}^{-}.
\end{align*}
In the case that \(\delta_{4}<1\) and \(\nu_{4}<0\), we define \(\Gamma_{4}^{-}\) as in \cref{sec:anal:ssec:ks}. In the remainder of this section, we will use the projected map to prove the following theorem about dynamics near the \(\Delta\)-clique network.

\begin{theorem}\label{thm:delta_clique}
  \begin{enumerate}[label=\normalfont({\alph*})]
    \item
    If \(\omega_{34}>0\) and \(\tau_{34}>\min\left\{2,1+\delta_{34}\right\}\), and either
    \begin{enumerate}[label=\normalfont({\roman*})]
      \item \(\zeta_{34}<0\) and \(\beta_{34}<0\); or \label{thm:delta_clique:cyc34:fas:I}
      \item \(\zeta_{34}<0\), \(\beta_{34}>0\), and \(\nu_{4}<0\); or \label{thm:delta_clique:cyc34:fas:II}
      \item \(\zeta_{34}>0\) and \(\nu_{4}<0\) \label{thm:delta_clique:cyc34:fas:III}
    \end{enumerate}
    then the \(\cyc{[34]}\) cycle is f.a.s. In cases \ref{thm:delta_clique:cyc34:fas:II} and \ref{thm:delta_clique:cyc34:fas:III}, and in case \ref{thm:delta_clique:cyc34:fas:I} if \(\nu_{4}<0\), all trajectories starting in \(\Gamma_{34}\) are asymptotic to the \(\cyc{[34]}\) cycle, with itinerary \((2341)^{\infty}\). In case \ref{thm:delta_clique:cyc34:fas:I} if \(\nu_{4}>0\), all trajectories starting in \(\Gamma_{34}^{+}\) are asymptotic to the \(\cyc{[34]}\) cycle, with itinerary \((2341)^{\infty}\), and all trajectories starting in \(\Gamma_{34}^{-}\) are mapped into \(\Gamma_{4}\). If \(\delta_{4}>1\), these trajectories switch from the \(\cyc{[34]}\) cycle to the \(\cyc{[4]}\) cycle, to which they are then asymptotic, and, for all positive integers \(n\), there exists an open set of trajectories in \(\Gamma_{34}^{-}\) with the itinerary \((2341)^{n}(241)^{\infty}\). \label{thm:delta_clique:cyc34:fas}

    \item
    If \(\omega_{34}>0\) and \(\tau_{34}>\min\left\{2,1+\delta_{34}\right\}\), and either
    \begin{enumerate}[label=\normalfont({\roman*})]
      \item \(\zeta_{34}<0\), \(\beta_{34}>0\), and \(\nu_{4}>0\); or \label{thm:delta_clique:cyc34:BCB:I}
      \item \(\zeta_{34}>0\) and \(\nu_{4}>0\) \label{thm:delta_clique:cyc34:BCB:II}
    \end{enumerate}
    then the \(\cyc{[34]}\) cycle is c.u., and all trajectories starting in \(\Gamma_{34}\)
    are mapped into \(\Gamma_{4}\). If \(\delta_{4}>1\), these trajectories switch from the \(\cyc{[34]}\) cycle to the \(\cyc{[4]}\) cycle, to which they are then asymptotic, and there exists an open set of trajectories in \(\Gamma_{34}\) with the itinerary \((2341)^{n}(241)^{\infty}\) for all positive integers \(n\) less than some \(N^{*}\).\label{thm:delta_clique:cyc34:BCB}

    \item
    If \(\omega_{34}>0\) and \(0<\tau_{34}<\min\left\{2,1+\delta_{34}\right\}\), then the \(\cyc{[34]}\) cycle is c.u. If any of the cases in \ref{thm:delta_clique:cyc34:BCB} hold, then all trajectories starting in \(\Gamma_{34}\) sufficiently close to the network are mapped into \(\Gamma_{4}\). If \(\delta_{4}>1\), these trajectories switch from the \(\cyc{[34]}\) cycle to the \(\cyc{[4]}\) cycle, to which they are then asymptotic, and there exists an open set of trajectories in \(\Gamma_{34}\) with the itinerary \((2341)^{n}(241)^{\infty}\) for all positive integers \(n\) less than some \(N^{*}\). If \ref{thm:delta_clique:cyc34:fas}\ref{thm:delta_clique:cyc34:fas:I} holds with \(\nu_{4}>0\), all trajectories starting in \(\Gamma_{34}^{-}\) sufficiently close to the network are mapped into \(\Gamma_{4}\). If \(\delta_{4}>1\), these trajectories switch from the \(\cyc{[34]}\) cycle to the \(\cyc{[4]}\) cycle, to which they are then asymptotic, and, for all positive integers \(n\), there exists an open set of trajectories in \(\Gamma_{34}^{-}\) with the itinerary \((2341)^{n}(241)^{\infty}\). All other trajectories starting near \(\cyc{[34]}\) leave a neighbourhood of \(\cyc{[34]}\). If \ref{thm:delta_clique:cyc34:fas}\ref{thm:delta_clique:cyc34:fas:I} holds with \(\nu_{4}<0\), or if \ref{thm:delta_clique:cyc34:fas}\ref{thm:delta_clique:cyc34:fas:II} or \ref{thm:delta_clique:cyc34:fas}\ref{thm:delta_clique:cyc34:fas:III} holds, all trajectories starting near \(\cyc{[34]}\) leave a neighbourhood of \(\cyc{[34]}\).\label{thm:delta_clique:cyc34:res}

    \item
    If \(\omega_{34}<0\) or \(\tau_{34}<0\), then the \(\cyc{[34]}\) cycle is c.u. All trajectories starting in \(\Gamma_{34}\) sufficiently close to the network are mapped into \(\Gamma_{4}\). If \(\delta_{4}>1\), these trajectories switch from the \(\cyc{[34]}\) cycle to the \(\cyc{[4]}\) cycle, to which they are then asymptotic, and there exists an open set of trajectories in \(\Gamma_{34}\) with the itinerary \((2341)^{n}(241)^{\infty}\) for all positive integers \(n\) less than some \(N^{*}\). All other trajectories starting near \(\cyc{[34]}\) leave a neighbourhood of \(\cyc{[34]}\).\label{thm:delta_clique:cyc34:trace_fold_switching}

    \item
    If \(\delta_{4}>1\) and \(\nu_{4}>0\), then the \(\cyc{[4]}\) cycle is f.a.s. All trajectories starting in \(\Gamma_{4}\) are asymptotic to the \(\cyc{[4]}\) cycle, with itinerary \((241)^{\infty}\). \label{thm:delta_clique:cyc4:fas}

    \item
    If \(\delta_{4}>1\) and \(\nu_{4}<0\), then the \(\cyc{[4]}\) cycle is a.c.u. and all trajectories starting in \(\Gamma_{4}\) are mapped into \(\Gamma_{34}\). If any of the condition in \ref{thm:delta_clique:cyc34:fas} hold, these trajectories switch from the \(\cyc{[4]}\) cycle to the \(\cyc{[34]}\) cycle, to which they are then asymptotic, and, for all positive integers \(n\), there exists an open set of trajectories in \(\Gamma_{4}\) with the itinerary \((241)^{n}(2341)^{\infty}\). \label{thm:delta_clique:cyc4:BCB}

    \item
    If \(\delta_{4}<1\), then the \(\cyc{[4]}\) cycle is c.u. If \(\nu_{4}<0\), then all trajectories starting in \(\Gamma_{4}^{-}\) sufficiently close to the network are mapped into \(\Gamma_{34}\). If the \(\cyc{[34]}\) cycle is f.a.s., these trajectories switch from the \(\cyc{[4]}\) cycle to the \(\cyc{[34]}\) cycle, to which they are then asymptotic, and, for all positive integers \(n\), there exists an open set of trajectories in \(\Gamma_{4}\) with the itinerary \((241)^{n}(2341)^{\infty}\). All other trajectories starting near \(\cyc{[4]}\) leave a neighbourhood of \(\cyc{[4]}\). \label{thm:delta_clique:cyc4:res}
  \end{enumerate}
\end{theorem}

Note that \ref{thm:delta_clique:cyc34:fas}\ref{thm:delta_clique:cyc34:fas:I} and \ref{thm:delta_clique:cyc4:fas} indicate a qualitative difference between the dynamics near the Kirk--Silber and the \(\Delta\)-clique networks: if there is bistability of the two cycles of the \(\Delta\)-clique network, then there must be switching from the \(\cyc{[34]}\) cycle to the \(\cyc{[4]}\) cycle, whereas there is no switching near the Kirk--Silber network if there is bistability of the two cycles.

The proofs of \ref{thm:delta_clique:cyc4:fas}, \ref{thm:delta_clique:cyc4:BCB}, and \ref{thm:delta_clique:cyc4:res} are the same as for the Kirk--Silber network, and so we omit them. The other parts of the theorem are proved in the following sections.

\Cref{fig:delta_clique_bif_set} gives a bifurcation set that summarises some of the results of \Cref{thm:delta_clique}.

\subsubsection{Fixed points, existence, and admissibility}\label{sec:anal:ssec_delta_clique:sssec:fp_exist_admiss}

The fixed points of the map \(f_{4}\), and their admissibility, is the same as for the Kirk--Silber network in \cref{sec:anal:ssec:ks:sssec:fp}.

The fixed points of the map \(f_{34}\) are
\begin{equation*}
  \vartheta_{34}^{*}=\frac{-\lambda_{34}^{+}+\alpha_{2}+\alpha_{4}}{\zeta_{34}}
  \qquad\textrm{ and }\qquad
  \vartheta_{34}^{-}=\frac{-\lambda_{34}^{-}+\alpha_{2}+\alpha_{4}}{\zeta_{34}}.
\end{equation*}
These fixed points exist if and only if \(\omega_{4}\geq 0\), and, if \(\omega_{4}=0\),
\begin{equation*}
  \vartheta_{34}^{*}=\vartheta_{34}^{-}=\vartheta_{34}^{c}\coloneqq \frac{2\alpha_{2}+\alpha_{4}-\alpha_{1}}{2\zeta_{34}}.
\end{equation*}
In the next subsection, we confirm that the bifurcation at \(\omega_{34}=0\) is a fold bifurcation. Assuming \(\omega_{34}=0\), we say that the fold bifurcation is virtual if \(\vartheta_{34}^{c}\) is virtual, and we say the fold bifurcation is admissible if \(\vartheta_{34}^{c}\) is admissible. If \(\tau_{34}<0\) or \(\zeta_{34}>0\), then, if \(\omega_{34}=0\), the fold is virtual. However, if \(\tau_{34}>0\) and \(\zeta_{34}<0\), then, if \(\omega_{34}=0\), the fold is virtual if and only if \(\beta_{34}>0\), and the fold is admissible if and only if \(\beta_{34}<0\).

If \(\tau_{34}<0\), then both \(\vartheta_{34}^{*}\) and \(\vartheta_{34}^{-}\) are virtual. If \(\tau_{34}>0\), the admissibility of \(\vartheta_{34}^{*}\) and \(\vartheta_{4}^{-}\) depends on \(\zeta_{34}\), \(\beta_{34}\) and \(\nu_{4}\). We present these results in \cref{tbl:delta_clique_admiss}. Since the admissibility of these fixed points depends on \(\nu_{4}\), at the point of the border-collision bifurcation, where \(\nu_{4}=0\), \(\vartheta_{4}^{*}\) will also undergo a border-collision bifurcation. This relationship can be seen from the continuity of the projected map of \(\Delta\)-clique network (see \cref{sec:anal:ssec:delta_clique:sssec:cont}): if \(f_{4}(\vartheta_{s})=\vartheta_{s}\), then \(f_{34}(\vartheta_{s})=\vartheta_{s}\) also and there is a border-collision bifurcation of a fixed point of \(f_{34}\).

With the same reasoning as in \cref{sec:anal:ssec:kirk_silber:sssec:BCB_PodIII_equiv} we derive the following proposition:
% , which allows us to study the stability of the heteroclinic cycles \(\cyc{[34]}\) and \(\cyc{[4]}\) at the splitting equilibrium \(\xi_{2}\).
%
\begin{proposition}\label{prop:d_clique:bcb_podIII_equiv}
  Suppose \(\omega_{34}>0\) and \(\tau_{34}>\min\left\{2,1+\delta_{34}\right\}\). The fixed point \(\vartheta_{34}^{*}\) is admissible if and only if \(w_{\max}\) of \(M_{3}^{(34)}\) satisfies Podvigina's third condition for stability. Suppose \(\delta_{4}>1\). The fixed point \(\vartheta_{4}^{*}\) is admissible if and only if \(w_{\max}\) of \(M_{4}^{(4)}\) satisfies Podvigina's third condition for stability.
\end{proposition}

\input{tbl_admiss.tex}

\subsubsection{Stability}\label{sec:anal:ssec:delta_clique:sssec:stab}

The stability of fixed points of the map \(f_{4}\) is the same as for the Kirk--Silber network in \cref{sec:anal:ssec:kirk_silber:sssec:stab}.

Calculating \(Df_{34}(\vartheta)\) gives
\begin{equation*}
  Df_{34}(\vartheta)=\frac{\delta_{34}}{\left(\zeta_{34}\vartheta-\alpha_{2}-\alpha_{4}\right)^{2}},
\end{equation*}
Since \(\delta_{34}>0\), \(Df_{34}\) is strictly positive. Evaluating at \(\vartheta_{34}^{*}\) and \(\vartheta_{34}^{-}\) gives
\begin{equation*}
  Df_{34}(\vartheta_{34}^{*})=\frac{\delta_{34}}{\tau_{34}\lambda_{34}^{+}-\delta_{34}}
  \qquad\textrm{ and }\qquad
  Df_{34}(\vartheta_{34}^{-})=\frac{\delta_{34}}{\tau_{34}\lambda_{34}^{-}-\delta_{34}}.
\end{equation*}
At the bifurcation point \(\omega_{4}=0\), where \(\vartheta_{34}^{*}=\vartheta_{34}^{-}=\vartheta_{34}^{c}\), \(Df_{34}(\vartheta_{34}^{*})=Df_{34}(\vartheta_{34}^{-})=1\), and therefore this bifurcation is generically a fold bifurcation.

Assume \(\omega_{34}>0\). If \(\tau_{34}>0\), \(0<Df_{34}(\vartheta_{34}^{*})<1\), and \(\vartheta_{34}^{*}\) is asymptotically stable, while \(1<Df_{34}(\vartheta_{34}^{-})\), and \(\vartheta_{34}^{-}\) is unstable. If \(\tau_{34}<0\), \(1<Df_{34}(\vartheta_{34}^{*})\), and \(\vartheta_{34}^{*}\) is unstable, while \(0<Df_{34}(\vartheta_{34}^{-})<1\), and \(\vartheta_{34}^{-}\) is asymptotically stable. We give bifurcation diagrams of the fixed points, and their admissibility, in \Cref{fig:delta_clique_bif_diag:subfig:admiss} when the fold is admissible, and in \Cref{fig:delta_clique_bif_diag:subfig:virt} when the fold is virtual.

Stability of a fixed point is equivalent to the eigenspace it represents having the eigenvalue with greatest absolute value. In the case of \(\cyc{[3]}\) and \(\cyc{[4]}\) cycles, this condition also corresponds to Podvigina's second stability condition \(\lambda_{\max}>1\). However, in the case of the \(\cyc{[34]}\) cycle, if \(\tau_{34}>0\) and so \(\vartheta_{34}^{*}\) is stable, \(\lambda_{\max}=\lambda_{34}^{+}\), whereas Podvigina's second condition is satisfied only if \(\tau_{34}>\min\{2,1+\delta_{34}\}\); that is, \(\lambda_{34}^{+}>1\).

All full transition matrices of the \(\cyc{[4]}\) cycle have an eigenvector with eigenvalue \(1\). Thus, when the loss of stability of this cycle is associated with Podvigina's second stability condition, there will also be a change in the stability of fixed points of the projected map, as the eigenvalue with greatest absolute value will also change. In the case of the \(\cyc{[34]}\) cycle, the eigenvalue \(\lambda_{34}^{+}\) can become smaller than \(1\) but remain \(\lambda_{\max}\), and so there is no change in the stability of the fixed points of the projected map. For this reason, we also need to impose the condition \(\tau_{34}>\min\left\{2,1+\delta_{34}\right\}\), which is not captured in the dynamics of the projected map.

\input{figure_8.tex}
\subsubsection{Continuity of the projected map of the \texorpdfstring{\(\Delta\)}{Delta}-clique network at \texorpdfstring{\(\vartheta_{s}\)}{}}\label{sec:anal:ssec:delta_clique:sssec:cont}

In the case of the \(\Delta\)-clique network, we find that
\begin{equation*}
  \lim_{\vartheta\searrow\vartheta_{s}}f_{34}(\vartheta)=\frac{-\left(c_{42}t_{13}+e_{12}c_{43}\right)}{c_{42}\left(t_{13}+c_{14}\right)+e_{12}c_{43}},
  \qquad\textrm{ and }\qquad
  \lim_{\vartheta\nearrow\vartheta_{s}}f_{4}(\vartheta)=\frac{-\left(c_{42}t_{13}+e_{12}c_{43}\right)}{c_{42}\left(t_{13}+c_{14}\right)+e_{12}c_{43}}.
\end{equation*}
Therefore, the projected map of the \(\Delta\)-clique network is continuous, which can be seen in all three examples in \Cref{fig:delta_clique_proj_map}, where the two functions have the same value at \(\vartheta_{s}\), given by the dashed red line.

\subsubsection{Dynamics with \texorpdfstring{\(\omega_{34}<0\)}{omega34<0} or \texorpdfstring{\(\tau_{34}<0\)}{tau34<0}}\label{sec:anal:ssec:delta_clique:sssec:fold_trace_switching}

Note that \(f(0),f(\vartheta_{s})\in S\), regardless of parameter values, and so \(f(0)<0\) and \(f(\vartheta_{s})>-1\). When \(\omega_{34}<0\) or \(\tau_{34}<0\), \(f_{34}\) either has no fixed points, or has no fixed points in \(\Theta_{34}\), respectively. Therefore, \(f_{34}(\vartheta)<\vartheta\) for all \(\vartheta\in\Theta_{34}\) and so, for all \(\vartheta\), there exists a positive integer \(n\) such that \(f_{34}^{n}(\vartheta)\in\Theta_{4}\). As in \cref{sec:anal:ssec:kirk_silber:sssec:virt}, we define, for all non-negative integers \(n\), \(\mathcal{E}_{n}\) as the \(n\)\textsuperscript{th} pre-image of \(\vartheta_{s}\) in \(\Theta_{34}\); that is, \(f_{4}^{n}\left(\mathcal{E}_{n}\right)=\vartheta_{s}\). In \cref{sec:anal:ssec:kirk_silber:sssec:virt}, there are simple expressions for \(\mathcal{E}_{n}\) due to the simplified form of the matrices \(M_{3}\) and \(M_{4}\) of the Kirk--Silber network. The nature of the matrix \(M_{34}\) of the \(\Delta\)-clique network means there is no simplified expression for \(\mathcal{E}_{n}\), beyond expressions involving matrix multiplication. Since there are no \(f_{34}\)-invariant sets in \(\Theta_{34}\) when \(\omega_{34}<0\) or \(\tau_{34}<0\), there is some \(N^{*}\) such that \(\mathcal{E}_{n}>0\) for all \(n\geq N^{*}\). We have therefore proven \ref{thm:delta_clique:cyc34:trace_fold_switching} of \cref{thm:delta_clique}.

An example of the dynamics of the projected map under these parameter values is shown in \Cref{fig:delta_clique_proj_map:subfig:post_fold}.

\input{figure_9.tex}
\subsubsection{Dynamics with a virtual fold}\label{sec:anal:ssec:delta_clique:sssec:virtual_fold}

We now consider dynamics when \(\omega_{34}>0\), \(\tau_{34}>0\), and when the fold bifurcation is virtual.

Under these parameter values, if \(\nu_{4}>0\), then both \(\vartheta_{34}^{*}\) and \(\vartheta_{34}^{-}\) are virtual, but \(\vartheta_{4}^{*}\) is admissible. For all \(\vartheta\in\Theta_{34}\), \(f_{34}(\vartheta)<\vartheta\). Therefore, for every \(\vartheta\in\Theta_{34}\), there is some \(n\) such that \(f^{n}(\vartheta)\in\Theta_{4}\). With the same reasoning as in \cref{sec:anal:ssec:delta_clique:sssec:fold_trace_switching}, the number of iterations \(n\) is bounded above by some \(N^{*}\), proving \ref{thm:delta_clique:cyc34:BCB}.

At \(\nu_{4}=0\), both \(\vartheta_{34}^{*}\) and \(\vartheta_{4}^{*}\) undergo a border-collision bifurcation and coincide with \(\vartheta_{s}\). For \(\nu_{4}<0\), \(\vartheta_{34}^{*}\) is admissible and \(\vartheta_{4}^{*}\) virtual. Under these parameter values, \(f_{34}(\vartheta)>\vartheta\) for all \(\vartheta\in\left(\vartheta_{s},\vartheta_{34}^{*}\right)\), and \(f_{34}(\vartheta)<\vartheta\) for all \(\vartheta\in\left(\vartheta_{34}^{*},0\right)\). Therefore, \(\vartheta_{34}^{*}\) is globally attracting in \(\Theta_{34}\). If \(\lambda_{34}^{+}>1\)---or, equivalently, \(\tau_{34}>\min\left\{2,1+\delta_{34}\right\}\)---the \(\cyc{[34]}\) cycle is f.a.s. We have therefore proven \ref{thm:delta_clique:cyc34:fas}\ref{thm:delta_clique:cyc34:fas:II} and \ref{thm:delta_clique:cyc34:fas}\ref{thm:delta_clique:cyc34:fas:III} of \cref{thm:delta_clique}.

In the case that \(0<\tau_{34}<\min\{2,1+\delta_{34}\}\), \(0<\lambda_{\max}<1\), and the \(\cyc{[34]}\) is c.u. \cite{podvigina_2012}. Trajectories beginning near the cycle leave a neighbourhood of the cycle if \(\nu_{4}<0\). If \(\nu_{4}>0\), trajectories beginning sufficiently close to the network will switch from the \(\cyc{[34]}\) to \(\cyc{[4]}\) cycle before they leave a neighbourhood of the network, completing part of the proof of \ref{thm:delta_clique:cyc34:res}.

An example of the dynamics of the projected map under these parameter values is shown in \Cref{fig:delta_clique_proj_map:subfig:post_bcb}. Note that the projected map looks qualitatively the same under these parameter values if the \(\cyc{[34]}\) cycle is unstable and \(0<\tau_{34}<\min\{2,1+\delta_{34}\}\). Therefore, in this case, an additional check must be completed to derive dynamics of trajectories near the network from the projected map.

\subsubsection{Dynamics with an admissible fold}\label{sec:anal:ssec:delta_clique:sssec:admissible_fold}

We last consider dynamics when \(\omega_{34}>0\), \(\tau_{34}>0\), and when the fold bifurcation is admissible.

If \(\nu_{4}>0\), then \(\vartheta_{34}^{-}\) and \(\vartheta_{34}^{*}\) are both admissible, and \(\vartheta_{34}^{-}<\vartheta_{34}^{*}\). For all \(\vartheta\in\left(\vartheta_{s},\vartheta_{34}^{-}\right)\cup\left(\vartheta_{34}^{*},0\right)\), \(f(\vartheta)<\vartheta\), and for all \(\vartheta\in\left(\vartheta_{34}^{-},\vartheta_{34}^{*}\right)\), \(f(\vartheta)>\vartheta\). Therefore, \(\vartheta_{34}^{*}\) is attracting in \(\left(\vartheta_{34}^{-},0\right)\), which corresponds to \(\Gamma_{34}^{+}\). Again, we require \(\tau_{34}>\min\left\{2,1+\delta_{34}\right\}\) for the \(\cyc{[34]}\) cycle to be f.a.s. For all \(\vartheta\in\left(\vartheta_{s},\vartheta_{34}^{-}\right)\), there is an \(n\) such that \(f^{n}(\vartheta)\in\Theta_{4}\). We again write \(\mathcal{E}_{n}\) for the \(n\)th preimage of \(\vartheta_{s}\). For these parameter values, \(\mathcal{E}_{n}<\mathcal{E}_{n+1}\), and \(\lim_{n\to\infty}\mathcal{E}_{n}=\vartheta_{34}^{-}\), and again there is no simplified expression. However, since this limit is a fixed point, \(n\) is not bounded above. The interval \(\left(\vartheta_{s},\vartheta_{34}^{-}\right)\) corresponds to the set \(\Gamma_{34}^{-}\).

At \(\nu_{4}=0\), both \(\vartheta_{34}^{-}\) and \(\vartheta_{4}^{*}\) undergo a border-collision bifurcation and coincide with \(\vartheta_{s}\). For \(\nu_{4}<0\), both \(\vartheta_{34}^{-}\) and \(\vartheta_{4}^{*}\) are virtual, but \(\vartheta_{34}^{*}\) remains admissible and is now globally attracting in \(\Theta_{34}\). Again, if \(\tau_{34}>\min\left\{2,1+\delta_{34}\right\}\), the \(\cyc{[34]}\) cycle is f.a.s, completing the proof of \ref{thm:delta_clique:cyc34:fas} of \cref{thm:delta_clique}. The case where \(0<\tau_{34}<\left\{2,1+\delta_{34}\right\}\) is the same as in \cref{sec:anal:ssec:delta_clique:sssec:virtual_fold}, except now, if \(\nu_{4}>0\), trajectories beginning in \(\Gamma_{34}^{-}\) sufficiently close to the network switch between cycles, completing the proof of \ref{thm:delta_clique:cyc34:res} of \cref{thm:delta_clique}, which also completes the proof of \cref{thm:delta_clique}.

An example of the dynamics of the projected map under these parameter values is shown in \Cref{fig:delta_clique_proj_map:subfig:pre_bcb}. Again, note that the projected map looks qualitatively the same under these parameter values if the \(\cyc{[34]}\) cycle is unstable and \(0<\tau_{34}<\min\{2,1+\delta_{34}\}\).

\subsubsection{Numerical examples}

\input{figure_10.tex}
We consider the ODEs
\begin{equation}\label{eqn:d_clique_dyn_sys}
  \begin{aligned}
    \dot{x}_{1}&=x_{1}\left(1-\left\lVert x\right\rVert_{2}^{2}-c_{21}x_{2}^{2}-t_{31}x_{3}^{2}+e_{41}x_{4}^{2}\right),\\
    \dot{x}_{2}&=x_{2}\left(1-\left\lVert x\right\rVert_{2}^{2}+e_{12}x_{1}^{2}-c_{32}x_{3}^{2}-c_{42}x_{4}^{2}\right),\\
    \dot{x}_{3}&=x_{3}\left(1-\left\lVert x\right\rVert_{2}^{2}-t_{13}x_{1}^{2}+e_{23}x_{2}^{2}-c_{43}x_{4}^{2}\right),\\
    \dot{x}_{4}&=x_{4}\left(1-\left\lVert x\right\rVert_{2}^{2}-c_{14}x_{1}^{2}+e_{24}x_{2}^{2}+e_{34}x_{3}^{2}\right),
  \end{aligned}
\end{equation}
where \(c_{jk}\), \(e_{jk}\), and \(t_{jk}\) are positive real numbers. System \cref{eqn:ks_dyn_sys} is \(\Z_{2}^{4}\)-equivariant, and contains a \(\Delta\)-clique network.

We again use the change of coordinates \(X_{i}\coloneq\log x_{j}\), and integrate the resulting ODEs \(\dot{X}_{j}\) with a variety of parameter values. Three different examples are found in \Cref{fig:delta_clique_examples}. Parameter values in \Cref{fig:delta_clique_examples:C34_cycle_fas,fig:delta_clique_examples:C34_cycle_fas_switch} are identical, and the projected map in these examples looks like \Cref{fig:delta_clique_proj_map:subfig:pre_bcb}. Both cycles are f.a.s., with a subset of \(\Gamma_{34}\) switching to the \(\cyc{[4]}\) cycle. These trajectories have itinerary \((2341)^{\infty}\) and \((2341)^{2}(241)^{\infty}\), respectively. The trajectory in \cref{sub@fig:delta_clique_examples:fold_switch} switches between cycles with itinerary \((2341)^{2}(241)^{\infty}\). Under these parameter values, the projected map looks like \Cref{fig:delta_clique_proj_map:subfig:post_fold}.

%% file: figure_7.tex
\begin{figure}
  \centering
  \begin{subfigure}[t]{0.49\linewidth}
    \centering
    \includegraphics[width=0.98523\linewidth]{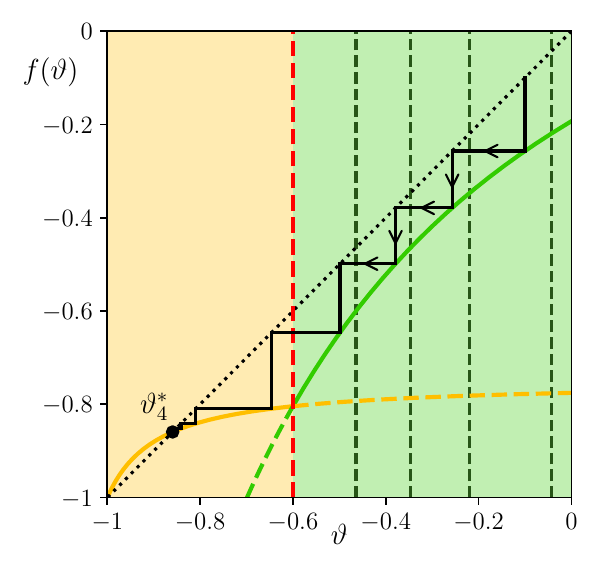}
    \vspace{-4mm}
    \caption{\(\omega_{34}<0\).}
    \label{fig:delta_clique_proj_map:subfig:post_fold}
  \end{subfigure}
  \begin{subfigure}[t]{0.49\linewidth}
    \centering
    \includegraphics[width=0.98523\linewidth]{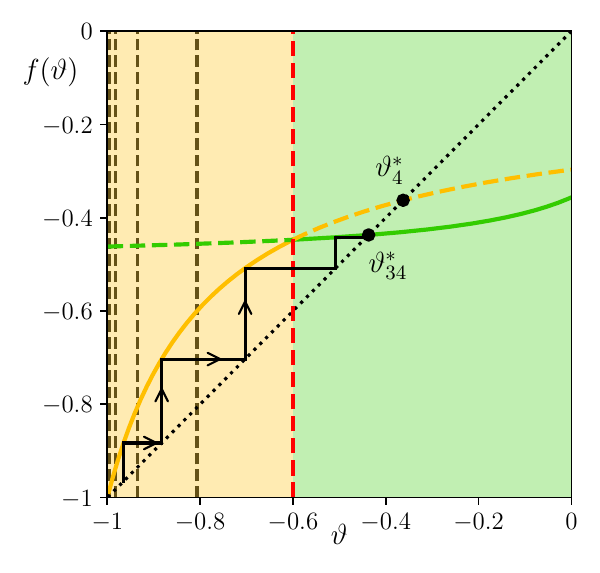}
    \vspace{-4mm}
    \caption{Virtual fold, with \(\nu_{4}<0\).}
    \label{fig:delta_clique_proj_map:subfig:post_bcb}
  \end{subfigure}
  \begin{subfigure}[t]{0.49\linewidth}
    \centering
    \includegraphics[width=0.98523\linewidth]{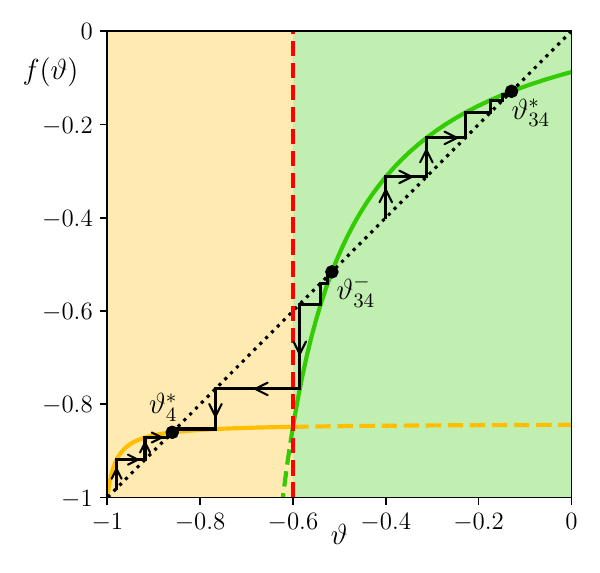}
    \vspace{-4mm}
    \caption{Admissible fold, with \(\nu_{4}>0\).}
    \label{fig:delta_clique_proj_map:subfig:pre_bcb}
  \end{subfigure}
  \begin{subfigure}[t]{0.49\linewidth}
    \centering
    \includegraphics[width=0.98523\linewidth]{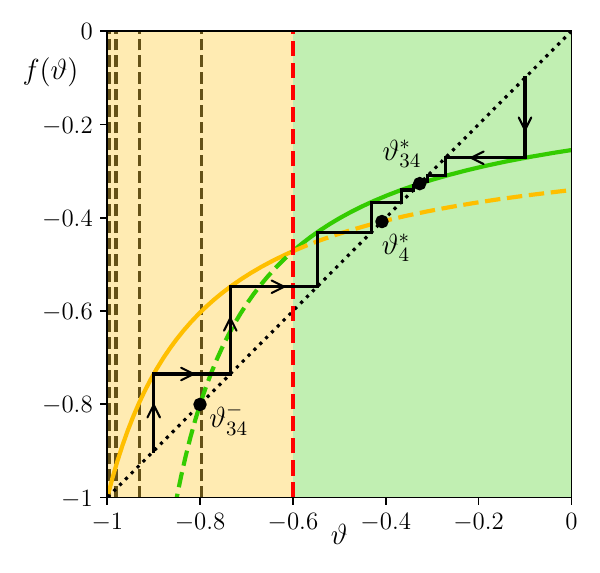}
    \vspace{-4mm}
    \caption{Admissible fold, with \(\nu_{4}<0\).}
    \label{fig:delta_clique_proj_map:subfig:admiss_fold_post_bcb}
  \end{subfigure}
  \caption{Schematic representations of the projected map of the \(\Delta\)-clique network as cobweb plots, for different relations between parameters. In all examples, the domain \(\Theta_{34}\) is shaded green and the domain \(\Theta_{4}\) orange. The dashed red line is \(\vartheta_{s}\). The dotted black line is \(y=\vartheta\). Each function \(f_{34}\) and \(f_{4}\) is shown as solid in its domain of definition and dashed otherwise. In all three figures, \(\delta_{4}>1\) and \(\tau_{34}>0\). The values of other expressions are indicated, as well as whether \(\vartheta_{34}^{*}\) and \(\vartheta_{34}^{-}\) have been created in a virtual or admissible fold. In \subref{fig:delta_clique_proj_map:subfig:post_fold}, \subref{fig:delta_clique_proj_map:subfig:post_bcb}, and \subref{fig:delta_clique_proj_map:subfig:admiss_fold_post_bcb}, the values \(\mathcal{E}_{n}\) are indicated by dashed black lines for \(n=1,2,3,4\). In \subref{fig:delta_clique_proj_map:subfig:pre_bcb}, the interval \(\left(\vartheta_{s},\vartheta_{34}^{-}\right)\) is also partitioned by \(\mathcal{E}_{n}\), but we do not plot these lines for visual clarity. The definition of virtual and admissible fold can be found in \cref{sec:anal:ssec_delta_clique:sssec:fp_exist_admiss}}
  \label{fig:delta_clique_proj_map}
\end{figure}

%% file: tbl_admiss.tex
%
% \begin{center}
  \begin{table}
    \centering
    \caption{Admissibility of \(\vartheta_{34}^{*}\) and \(\vartheta_{34}^{-}\) when \(\tau_{34}>0\).}\label{tbl:delta_clique_admiss}
    \vspace{5pt}
    \begin{tabular}{c c | c c}
      \hline
      \hline
                                          &                                     & \(\nu_{4}>0\)                   & \(\nu_{4}<0\) \\
      \hline
      \multirow{4}{4em}{\(\zeta_{34}<0\)} & \multirow{2}{4em}{\(\beta_{34}>0\)} & \(\vartheta_{34}^{*}\): virtual & \(\vartheta_{34}^{*}\): admissible\\
                                          &                                     & \(\vartheta_{34}^{-}\): virtual & \(\vartheta_{34}^{-}\): virtual\\[0.5em]
                                          & \multirow{2}{4em}{\(\beta_{34}<0\)} & \(\vartheta_{34}^{*}\): admissible & \(\vartheta_{34}^{*}\): admissible\\
                                          &                                     & \(\vartheta_{34}^{-}\): admissible & \(\vartheta_{34}^{-}\): virtual\\
      \hline
      \multirow{2}{4em}{\(\zeta_{34}>0\)} & \multirow{2}{1.3em}{---} & \(\vartheta_{34}^{*}\): virtual & \(\vartheta_{34}^{*}\): admissible\\
                                          &                          & \(\vartheta_{34}^{-}\): virtual & \(\vartheta_{34}^{-}\): virtual\\
      \hline
      \hline
    \end{tabular}
  \end{table}
% \end{center}
%

%% file: figure_8.tex
\begin{figure}
  \centering
  \includegraphics[width=0.687\linewidth]{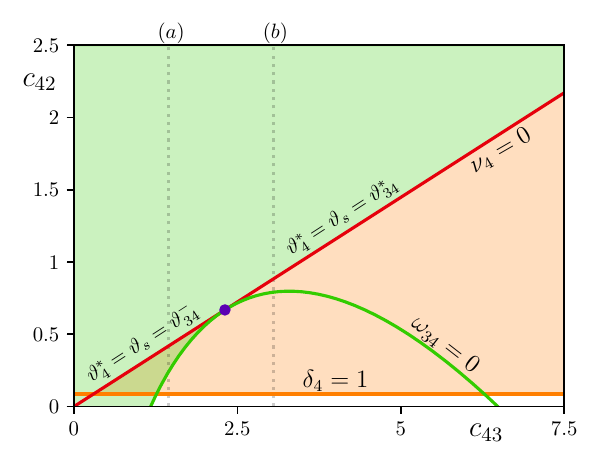}
  \caption{A bifurcation set of the \(\Delta\)-clique network. The set is shaded orange where the \(\cyc{[4]}\) cycle is f.a.s., and shaded green where the \(\cyc{[34]}\) cycle is f.a.s. Solid lines correspond to relations between parameter values as indicated. The green line is the locus of the fold bifurcation of \(\vartheta_{34}^{*}\) and \(\vartheta_{34}^{-}\). The red line is the locus of the border-collision bifurcation of \(\vartheta_{4}^{*}\) and also one of \(\vartheta_{34}^{*}\) and \(\vartheta_{34}^{-}\), and the orange line is the locus of the transcritical bifurcation of \(\vartheta_{4}^{*}\) and \(\vartheta_{4}^{-}\). Faint dotted lines correspond as labelled to parameter values of the bifurcation diagrams in \cref{fig:delta_clique_bif_diag}. For all parameter values in the bifurcation set, \(\tau_{34}>\min\{2,1+\delta_{34}\}\). The purple dot is the location of the codimension-\(2\) point where \(\vartheta_{34}^{c}=\vartheta_{s}\) if \(\omega_{34}=0\).}
  \label{fig:delta_clique_bif_set}
\end{figure}

%% file: figure_9.tex
\begin{figure}
  \centering
  \begin{subfigure}[t]{0.5\linewidth}
    \centering
    \includegraphics[width=0.9\linewidth]{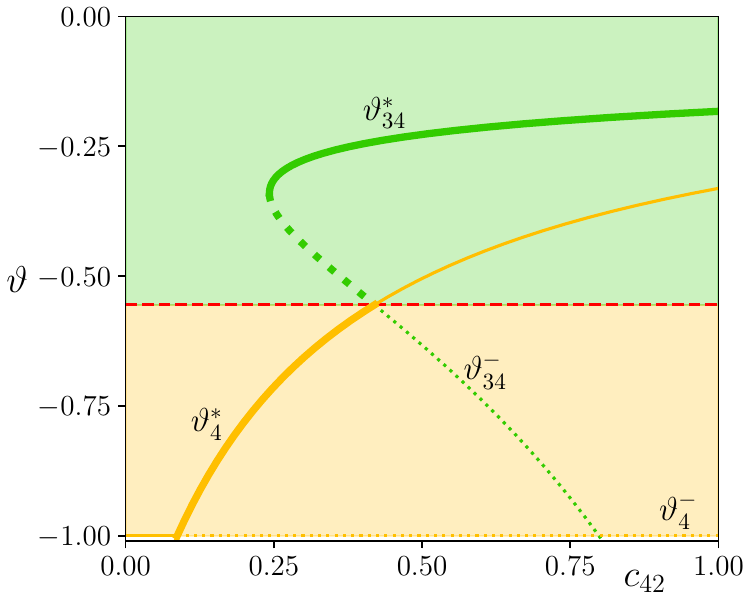}
    % \vspace{-1mm}
    \caption{\(c_{43}=1.45\): admissible fold.}
    \label{fig:delta_clique_bif_diag:subfig:admiss}
  \end{subfigure}%
  \hfill%
  \begin{subfigure}[t]{0.5\linewidth}
    \centering
    \includegraphics[width=0.9\linewidth]{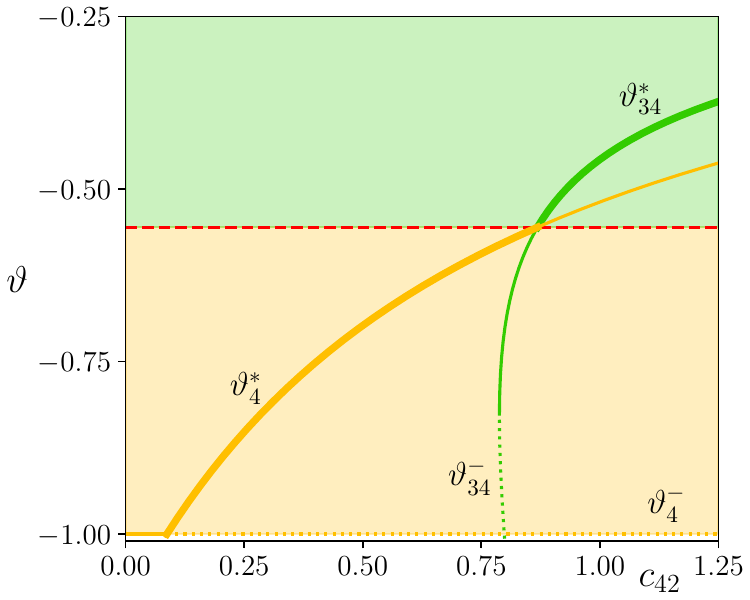}
    % \vspace{-1mm}
    \caption{\(c_{43}=3.05\): virtual fold.}
    \label{fig:delta_clique_bif_diag:subfig:virt}
  \end{subfigure}
  \caption{Bifurcation diagrams of the fixed points of \(f_{4}\) and \(f_{34}\). The domain of \(f_{4}\) is shaded orange, and of \(f_{34}\) is shaded green. Stable fixed points are solid, and unstable dotted. Admissible fixed points are thick lines and virtual fixed points are thin lines. \(\vartheta_{s}\) is the dashed red line. All other parameters are the same as those in \Cref{fig:delta_clique_bif_set}.}
  \label{fig:delta_clique_bif_diag}
\end{figure}

%% file: figure_10.tex
\begin{figure}
  \begin{subfigure}[c]{\linewidth}
    \includegraphics[width=0.46\linewidth]{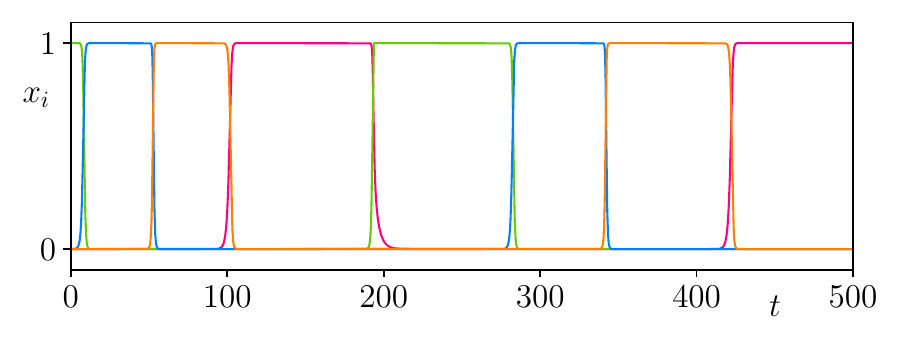}%
    \includegraphics[width=0.53\linewidth]{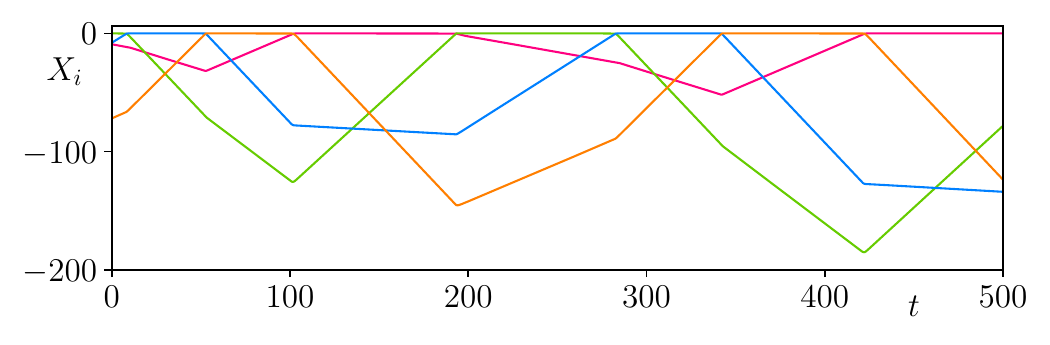}
    \vspace{-3mm}
    \caption{Admissible fold, with \(x(0)\in\Gamma_{34}^{+}\)}
    \label{fig:delta_clique_examples:C34_cycle_fas}
  \end{subfigure}
  \begin{subfigure}[c]{\linewidth}
    \includegraphics[width=0.46\linewidth]{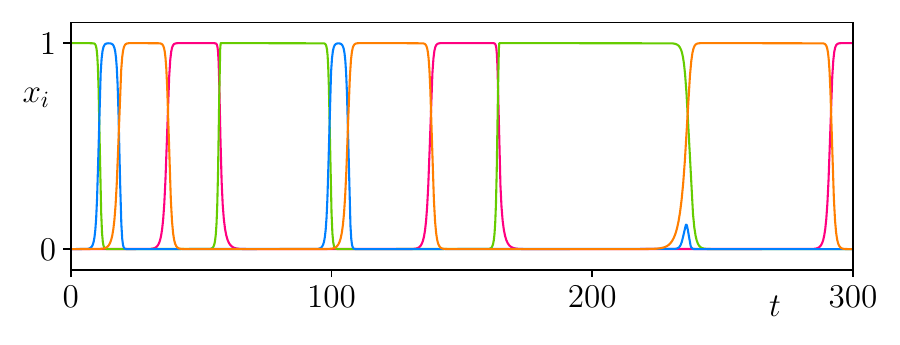}%
    \includegraphics[width=0.53\linewidth]{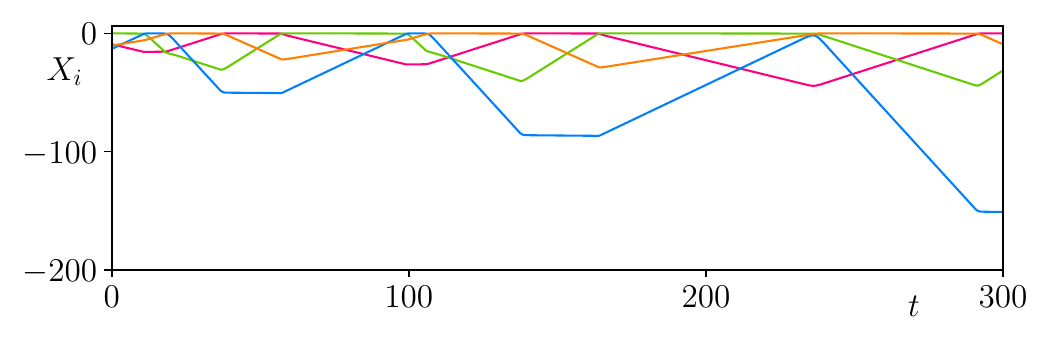}
    \vspace{-3mm}
    \caption{Admissible fold, with \(x(0)\in\Gamma_{34}^{-}\)}
    \label{fig:delta_clique_examples:C34_cycle_fas_switch}
  \end{subfigure}
  \begin{subfigure}[c]{\linewidth}
    \includegraphics[width=0.46\linewidth]{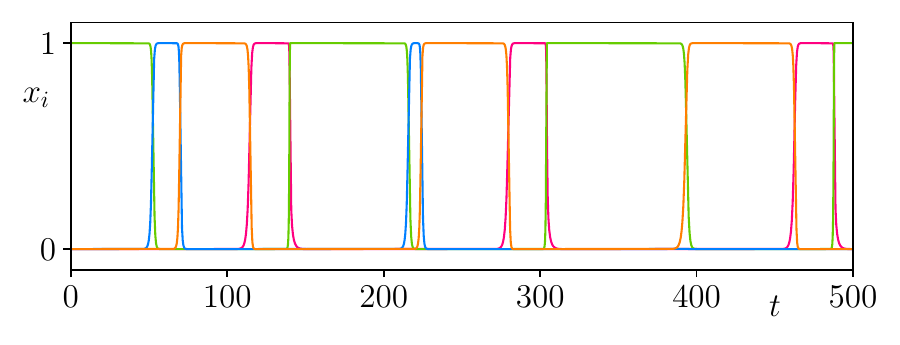}%
    \includegraphics[width=0.53\linewidth]{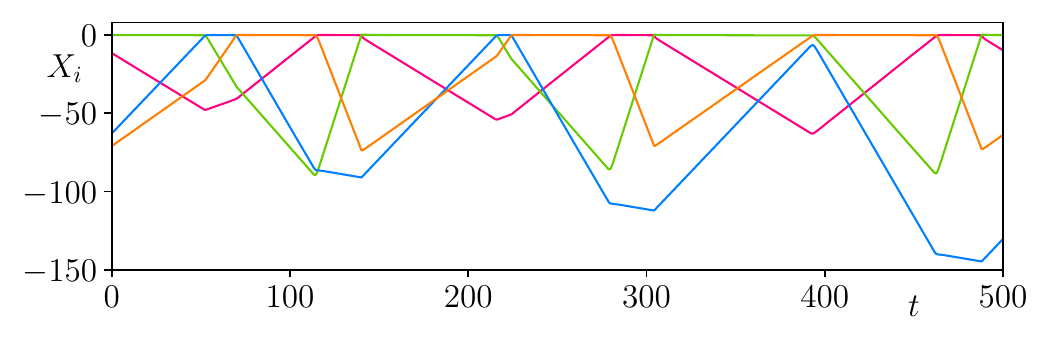}
    \vspace{-3mm}
    \caption{\(\omega_{34}<0\)}
    \label{fig:delta_clique_examples:fold_switch}
  \end{subfigure}
  \caption{Example timeseries of trajectories near the \(\Delta\)-clique network. The colours red, blue, green, and orange correspond to \(x_{1}\), \(x_{2}\), \(x_{3}\), and \(x_{4}\), respectively, and similarly for the \(X_{i}\). In all three examples, \(\delta_{4}>1\) and \(\nu_{4}>0\).}
  \label{fig:delta_clique_examples}
\end{figure}

%% file: 4_3_tournament.tex
\subsection{The tournament network}\label{sec:anal:ssec:tournament}

The tournament network is defined as the union of the cycles \(\cyc{[3]}\), \(\cyc{[4]}\), and \(\cyc{[34]}\). Therefore, the tournament network---unlike the Kirk--Silber and \(\Delta\)-clique networks---contains subnetworks:
\begin{enumerate}
  \item a subnetwork equivalent to the Kirk--Silber network, \(\cyc{[3]}\cup\cyc{[4]}\); and
  \item two subnetworks equivalent to the \(\Delta\)-clique network: \(\cyc{[34]}\cup\cyc{[4]}\), and \(\cyc{[3]}\cup\cyc{[34]}\).
\end{enumerate}
The tournament network has two splitting equilibria: \(\xi_{2}\) and \(\xi_{3}\). All trajectories leaving a small neighbourhood of \(\xi_{2}\) in the direction of \(\xi_{4}\) cycle around \(\cyc{[4]}\). Those trajectories that leave \(\xi_{2}\) in the direction of \(\xi_{3}\) will then leave a small neighbourhood of \(\xi_{3}\) in the direction of \(\xi_{1}\) or \(\xi_{4}\), and therefore cycle around \(\cyc{[3]}\) or \(\cyc{[34]}\), respectively. However, all trajectories will eventually return to \(\xi_{2}\) if the network is attracting.

We consider subsets of points in \(\poinin[2]{3}\) that give rise to trajectories that cycle around \(\cyc{[3]}\) or \(\cyc{[34]}\). We take the pre-image under the basic map \(\varphi_{123}\colon\poinin[1]{2}\to\poinin[2]{3}\) of these two domains, and partition \(\poinin[1]{2}\) into three domains: \(\Gamma_{3}\), \(\Gamma_{4}\), and \(\Gamma_{34}\). These domains correspond to trajectories that, after leaving a small neighbourhood of \(\xi_{2}\), next cycle around the \(\cyc{[3]}\), \(\cyc{[4]}\), or \(\cyc{[34]}\) cycle. This partition allows us to define three corresponding subsets of \(S\):
\begin{equation*}
  \Theta_{3}=\left(\vartheta_{s}^{+},0\right),
  \ 
  \Theta_{34}=\left(\vartheta_{s}^{-},\vartheta_{s}^{+}\right),
  \textrm{ and }\ 
  \Theta_{4}=\left(-1,\vartheta_{s}^{-}\right),
\end{equation*}
where
\begin{equation*}
  \vartheta_{s}^{-}=\frac{-1}{1+\frac{e_{24}}{e_{23}}}
  \qquad\textrm{ and }\qquad
  \vartheta_{s}^{+}=\frac{-1}{1+\left(\frac{c_{21}e_{34}}{e_{23}e_{31}}+\frac{e_{24}}{e_{23}}\right)}.
\end{equation*}
We can then derive the projected map of the tournament network:
\begin{equation}\label{eqn:tournament_proj_map}
  f(\vartheta)=\begin{cases}
    f_{3}(\vartheta)=\dfrac{-\delta_{3}\vartheta}{\left(\delta_{3}+\rho_{3}-1\right)\vartheta-1}, & \textrm{if } \vartheta\in\Theta_{3},\\[1.79em]
    f_{34}(\vartheta)=\dfrac{\left(\alpha_{2}-\alpha_{1}\right)\vartheta + \alpha_{2}}{\zeta_{34}\vartheta-\alpha_{2}-\alpha_{4}}, & \textrm{if } \vartheta\in\Theta_{34},\\[1.79em]
    f_{4}(\vartheta)=\dfrac{\left(1-\rho_{4}\right)\vartheta-\rho_{4}}{\left(\delta_{4}+\rho_{4}-1\right)\vartheta+\delta_{4}+\rho_{4}}, & \textrm{if } \vartheta\in\Theta_{4}.
  \end{cases}
\end{equation}
The expressions for each of the constants in this function are given in the Supplementary Material (section \ref{supp:maps:tournament}).

\input{figure_11.tex}
\input{figure_12.tex}
\input{figure_13.tex}
The complete analysis of this projected map is subtle, and involves many cases. In this section, we will not state a complete description of all dynamics that can be observed near the tournament network, as we did in \Cref{thm:ks,thm:delta_clique}. Instead, we will give a qualitative description of how the dynamics can be analysed, and several examples that show interesting dynamics near the network.

The maps \(f_{3}\), \(f_{34}\), and \(f_{4}\) are analogous to those studied in \cref{sec:anal:ssec:ks,sec:anal:ssec:delta_clique}, though with differing expressions for some constants used due to the change in the global classification of some eigenvalues. As in the case of the Kirk--Silber and \(\Delta\)-clique networks, the \(\cyc{[3]}\) and \(\cyc{[4]}\) cycles are stable if their corresponding map has an admissible and attracting fixed point. In the case of the \(\cyc{[34]}\) cycle, \(\vartheta_{34}^{*}\) must be admissible and stable, but we must additionally impose \(\tau_{34}>\min\{2,1+\delta_{34}\}\). In all three cases, it is straightforward to verify that admissibility of a fixed point corresponds to all transition matrices satisfying Podvigina's third stability condition in \cref{thm:pod_bif}, assuming the first two conditions are also satisfied, as in \cref{prop:ks:bcb_podIII_equiv,prop:d_clique:bcb_podIII_equiv}.

The stability and admissibility of \(\vartheta_{4}^{*}\) is the same as in \cref{sec:anal:ssec:ks,sec:anal:ssec:delta_clique}. The condition \(\nu_{3}=0\) corresponds now to \(\vartheta_{3}^{*}=\vartheta_{s}^{-}\). However, \(\vartheta_{3}^{*}\) is now virtual when \(\vartheta_{3}^{*}<\vartheta_{s}^{+}\), when \(\mu_{3}<0\), where
\begin{equation*}
  \mu_{3}=\frac{c_{14}}{e_{12}} - \frac{c_{13}e_{24}}{e_{12}e_{23}} - \frac{c_{13}c_{21}e_{34}}{e_{12}e_{23}e_{31}}.
\end{equation*}
At \(\mu_{3}=0\), \(w_{\max}\) of \(M_{1}^{(3)}\) has a zero entry, and so does not satisfy Podvigina's third stability condition.

The admissibility of the fixed points \(\vartheta_{34}^{*}\) and \(\vartheta_{34}^{-}\) is the complicating factor for the analysis of the tournament network's projected map. The projected map of the tournament network, like the \(\Delta\)-clique network, is continuous (see \cref{sec:anal:ssec:tournament:sssec:cont}). Therefore, border-collision bifurcations of fixed points of \(f_{34}\) also occur at \(\nu_{4}=0\) and \(\mu_{3}=0\). However, which fixed point, \(\vartheta_{34}^{*}\) or \(\vartheta_{34}^{-}\), changes admissibility, and which direction it crosses which switching manifold, depends on whether the fold bifurcation of \(\vartheta_{34}^{*}\) and \(\vartheta_{34}^{-}\) occurs in \(\Theta_{3}\), \(\Theta_{34}\), \(\Theta_{4}\), or outside \(S\), and also depends on the concavity of \(f_{34}\) in \(\Theta_{34}\); that is, on whether \(\zeta_{34}<0\) or \(\zeta_{34}>0\).

We focus on the case where \(\zeta_{34}<0\), and so the function \(f_{34}\) is concave down in \(\Theta_{34}\). We assume, when \(\omega_{34}>0\), that \(\tau_{34}>\min\left\{2,1+\delta_{34}\right\}\). We will also assume \(\delta_{3}>1\) and \(\delta_{4}>1\). We show various examples of the projected map of the tournament network in \Cref{fig:tournament_proj_map}. For visual clarity, we plot each function only in its domain of definition, and plot multiple examples before and after certain border-collision bifurcations in the same plot, distinguished by lines of different shades.

To aid in the understanding of the dynamics we describe below, we present a bifurcation set in \Cref{fig:tournament_bif_set}. We also provide four qualitatively different bifurcation diagrams in \Cref{fig:tournament_bif_diag}, which capture the dynamics we have described here.

\Cref{fig:tournament_proj_map:subfig:pre_fold} shows the projected map with \(\omega_{34}<0\). We see that, when \(\omega_{34}<0\), only one of the \(\cyc{[3]}\) or \(\cyc{[4]}\) cycles is stable, and that all trajectories near the network are asymptotic to this cycle. Trajectories therefore switch between cycles, and can switch from a length-three cycle to the \(\cyc{[34]}\) cycle, before switching again to the stable length-three cycle. If parameter values are such that \(f_{3}(\vartheta_{s}^{+})<\vartheta_{s}^{-}\), or \(f_{4}(\vartheta_{s}^{-})>\vartheta_{s}^{+}\), trajectories may switch from one length-three cycle to the other, without cycling around the \(\cyc{[34]}\) cycle.

\Cref{fig:tournament_proj_map:subfig:admiss_bcb} shows the possible dynamics if the fold bifurcation of \(f_{34}\) is admissible, and so occurs in \(\Theta_{34}\). First, if \(\nu_{4}<0\) and \(\mu_{3}>0\), all trajectories are asymptotic to the \(\cyc{[3]}\) cycle, and it is the only stable cycle. If \(\nu_{4}>0\) and \(\mu_{3}>0\), the length-three cycles are bistable, and \(\Theta_{34}\) is partitioned by the unstable fixed point \(\vartheta_{34}^{-}\) into subsets of points that are mapped into \(\Theta_{3}\) or \(\Theta_{4}\). If \(\nu_{4}>0\) and \(\mu_{3}<0\), there is bistability of the \(\cyc{[4]}\) and \(\cyc{[34]}\) cycles. The \(\cyc{[3]}\) cycle is unstable, and all trajectories in its domain are asymptotic to the \(\cyc{[34]}\) cycle. Since \(\nu_{4}>0\), \(\vartheta_{34}^{-}\) is admissible, and \(\Theta_{34}\) is partitioned by the unstable fixed point \(\vartheta_{34}^{-}\) into subsets of points that are mapped into \(\Theta_{4}\) and those that are asymptotic to the \(\theta_{34}^{*}\). Lastly, if \(\nu_{4}<0\) and \(\mu_{3}<0\), only the \(\cyc{[34]}\) cycle is stable, and all trajectories are asymptotic to it. Therefore, under the various cases, trajectories can switch once or twice in total, but only once between two given cycles. These results can be seen in \Cref{fig:tournament_bif_diag:subfig:admiss_I,fig:tournament_bif_diag:subfig:admiss_II}, which differ by the admissibility of \(\vartheta_{34}^{*}\) between the two border-collision bifurcations. In \Cref{fig:tournament_bif_diag:subfig:admiss_I}, between the two border-collisions, \(\cyc{[3]}\) and \(\cyc{[4]}\) are bistable, and, in \Cref{fig:tournament_bif_diag:subfig:admiss_II}, only \(\cyc{[34]}\) is stable.

In \Cref{fig:tournament_proj_map:subfig:theta_3_bcb}, we show the possible dynamics if the fold bifurcation of \(f_{34}\) is virtual, and occurs in \(\Theta_{3}\). Under these conditions, the \(\cyc{[34]}\) cycle is always unstable, and there is bistability of the length-three cycles if both \(\mu_{3}>0\) and \(\nu_{4}>0\). In this case, \(\Theta_{34}\) is again partitioned into subsets of points that are mapped into \(\Theta_{3}\) or \(\Theta_{4}\). If \(\nu_{4}<0\), all trajectories are asymptotic to \(\cyc{[3]}\), and if \(\mu_{3}<0\), all trajectories are asymptotic to \(\cyc{[4]}\). We can have switching once or twice in total, and again only once between two given cycles. These results can be observed in \Cref{fig:tournament_bif_diag:subfig:Theta3_bcb}.

In \Cref{fig:tournament_proj_map:subfig:theta_4_bcb}, we show the possible dynamics if the fold bifurcation of \(f_{34}\) is virtual, and occurs in \(\Theta_{4}\), or outside of \(S\). Under these conditions, any cycle can be stable, but at most one cycle is stable for a given set of parameter values. The stable cycle is \(\cyc{[3]}\) if \(\mu_{3}>0\), \(\cyc{[4]}\) if \(\nu_{4}>0\), and it is \(\cyc{[34]}\) if \(\nu_{4}<0\) and \(\mu_{3}<0\). All trajectories are asymptotic to the stable cycle, and we again have switching once or twice in total, and only once between two given cycles. These results can be observed in \Cref{fig:tournament_bif_diag:subfig:Theta4_bcb}.

A similar analysis can be done in the case of \(\zeta_{34}>0\).\footnote{For a visual indication of the results, look at \Cref{fig:tournament_proj_map} upside down, swapping the yellow and blue colours to swap the \(\Theta_{3}\) and \(\Theta_{4}\) domains, and also swapping the captions of \subref{fig:tournament_proj_map:subfig:theta_3_bcb} and \subref{fig:tournament_proj_map:subfig:theta_4_bcb}.}

\subsubsection{Continuity of the projected map of the tournament network at \texorpdfstring{\(\vartheta_{s}^{+}\)}{} and \texorpdfstring{\(\vartheta_{s}^{-}\)}{}}\label{sec:anal:ssec:tournament:sssec:cont}

In the case of the tournament network, we find at \(\vartheta_{s}^{+}\) that
\begin{equation*}
  \lim_{\vartheta\searrow\vartheta_{s}^{+}}f_{3}(\vartheta)=\frac{-c_{13}}{c_{13} + c_{14}}
  \qquad\textrm{ and }\qquad
  \lim_{\vartheta\nearrow\vartheta_{s}^{+}}f_{34}(\vartheta)=\frac{-c_{13}}{c_{13} + c_{14}},
\end{equation*}
and we find at \(\vartheta_{s}^{-}\) that
\begin{equation*}
  \lim_{\vartheta\searrow\vartheta_{s}^{-}}f_{34}(\vartheta)=\frac{-(c_{42}c_{13} + e_{12}c_{43})}{c_{42}\left(c_{13} + c_{14}\right) + e_{12}c_{43}},
  \qquad\textrm{ and }\qquad
  \lim_{\vartheta\nearrow\vartheta_{s}^{-}}f_{4}(\vartheta)=\frac{-(c_{42}c_{13} + e_{12}c_{43})}{c_{42}\left(c_{13} + c_{14}\right) + e_{12}c_{43}}.
\end{equation*}
Therefore, the projected map of the tournament network is continuous at both \(\vartheta_{s}^{+}\) and \(\vartheta_{s}^{-}\), as seen in all examples in \Cref{fig:tournament_proj_map}.

\subsubsection{Numerical examples}

\input{figure_14.tex}
For the tournament network, we consider the system of ODEs

\begin{equation}\label{eqn:tournament_dyn_sys}
  \begin{aligned}
    \dot{x}_{1}&=x_{1}\left(1-\left\lVert x\right\rVert_{2}^{2}-c_{21}x_{2}^{2}+e_{31}x_{3}^{2}+e_{41}x_{4}^{2}\right),\\
    \dot{x}_{2}&=x_{2}\left(1-\left\lVert x\right\rVert_{2}^{2}+e_{12}x_{1}^{2}-c_{32}x_{3}^{2}-c_{42}x_{4}^{2}\right),\\
    \dot{x}_{3}&=x_{3}\left(1-\left\lVert x\right\rVert_{2}^{2}-c_{13}x_{1}^{2}+e_{23}x_{2}^{2}-c_{43}x_{4}^{2}\right),\\
    \dot{x}_{4}&=x_{4}\left(1-\left\lVert x\right\rVert_{2}^{2}-c_{14}x_{1}^{2}+e_{24}x_{2}^{2}+e_{34}x_{3}^{2}\right),
  \end{aligned}
\end{equation}
where \(c_{jk}\) and \(e_{jk}\) are positive real numbers. System \cref{eqn:tournament_dyn_sys} is \(\Z_{2}^{4}\)-equivariant, and contains a heteroclinic network equivalent to the tournament network. All four equilibria of the network lie at unit distance from the origin on the four coordinate axes.

Again we use the change of coordinates \(X_{i}\coloneq\log x_{j}\), and integrate the resulting ODEs \(\dot{X}_{j}\). We show one example in \Cref{fig:tournament_example}. This trajectory has itinerary \((231)^{2}(2341)^{2}(241)^{\infty}\), and therefore switches twice between cycles.
The onset of switching from the \(\cyc{[3]}\) cycle to the \(\cyc{[34]}\) observed in \cref{fig:tournament_example} is the result of a border-collision bifurcation, as \(\mu_{3}<0\), whereas the onset of switching from the \(\cyc{[34]}\) cycle to the \(\cyc{[4]}\) cycle is a result of a fold bifurcation, as \(\omega_{34}<0\).

%% file: figure_11.tex
\begin{figure}[p!]
  \centering
  \begin{subfigure}[t]{0.49\linewidth}
    \centering
    \includegraphics[width=\linewidth]{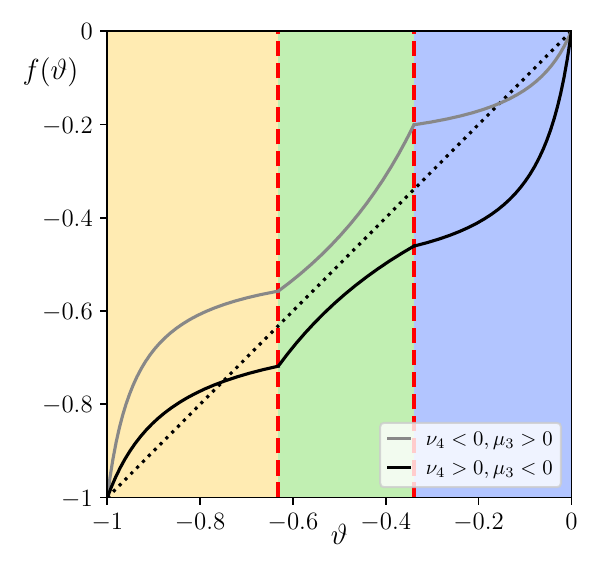}
    \vspace{-8mm}
    \caption{\(\omega_{34}<0\)}
    \label{fig:tournament_proj_map:subfig:pre_fold}
  \end{subfigure}
  \begin{subfigure}[t]{0.49\linewidth}
    \centering
    \includegraphics[width=\linewidth]{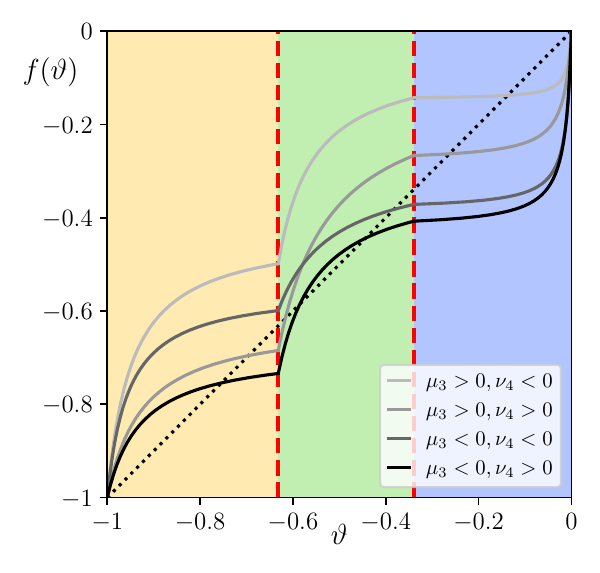}
    \vspace{-8mm}
    \caption{Admissible fold in \(\Theta_{34}\)}
    \label{fig:tournament_proj_map:subfig:admiss_bcb}
  \end{subfigure}
  \begin{subfigure}[t]{0.49\linewidth}
    \centering
    \includegraphics[width=\linewidth]{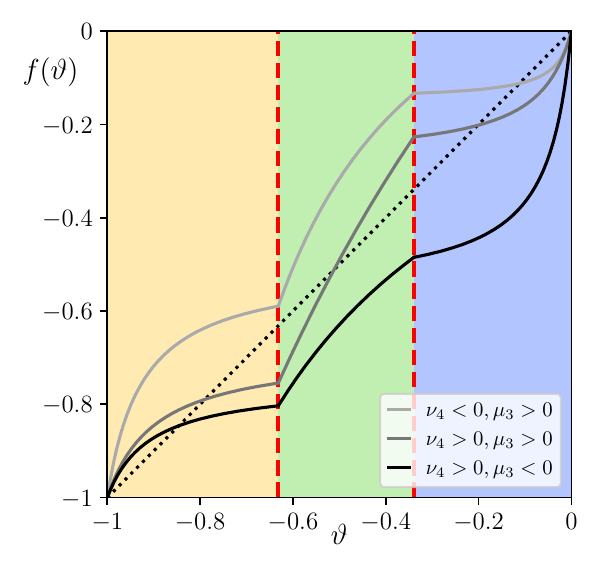}
    \vspace{-8mm}
    \caption{Virtual fold in \(\Theta_{3}\)}
    \label{fig:tournament_proj_map:subfig:theta_3_bcb}
  \end{subfigure}
  \begin{subfigure}[t]{0.49\linewidth}
    \centering
    \includegraphics[width=\linewidth]{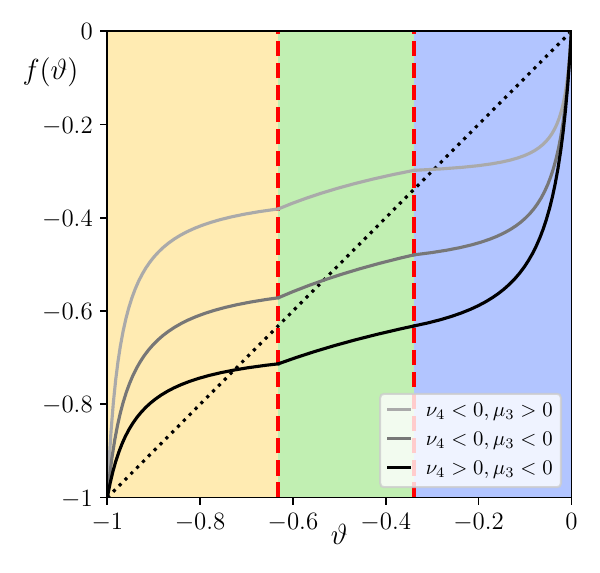}
    \vspace{-8mm}
    \caption{Virtual fold in \(\Theta_{4}\), or outside \(S\).}
    \label{fig:tournament_proj_map:subfig:theta_4_bcb}
  \end{subfigure}
  \caption{Schematic representations of the projected map of the tournament network, for different relations between parameters. In \subref{fig:tournament_proj_map:subfig:pre_fold}, the fixed points \(\vartheta_{34}^{*}\) and \(\vartheta_{34}^{-}\) do not exist. In \subref{fig:tournament_proj_map:subfig:admiss_bcb}, \subref{fig:tournament_proj_map:subfig:theta_3_bcb}, and \subref{fig:tournament_proj_map:subfig:theta_4_bcb}, \(\omega_{34}>0\) and so the fixed points exist, and were created in a fold bifurcation, the admissibility of which is as indicated. In all plots, \(\delta_{3}>1\) and \(\delta_{4}>1\). The values of \(\mu_{3}\) and \(\nu_{4}\) are as indicated by the shading of the plot. The regions shaded yellow, green, and blue correspond to \(\Theta_{4}\), \(\Theta_{34}\), and \(\Theta_{3}\), respectively. The left dashed red line is \(\vartheta_{s}^{-}\), and the right dashed red line is \(\vartheta_{s}^{+}\).}
  \label{fig:tournament_proj_map}
\end{figure}

%% file: figure_12.tex
\begin{figure}
  \centering
  \includegraphics[width=0.8\linewidth]{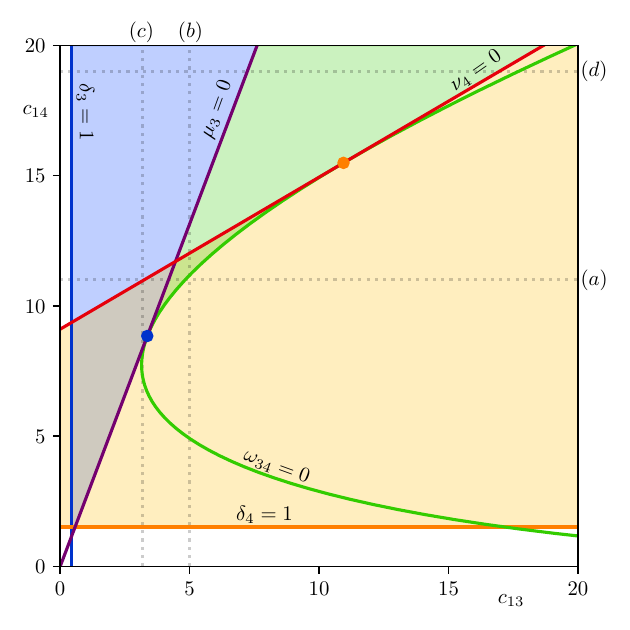}
  \caption{A bifurcation set of the tournament network, for \(\tau_{34}>\min\{2,1+\delta_{34}\}\) and \(\zeta_{34}<0\). Where the set is shaded blue, green, or orange, the \(\cyc{[3]}\), \(\cyc{[34]}\), or \(\cyc{[4]}\) cycle, respectively, is f.a.s. Coloured lines correspond to relations between parameter values as indicated. Faint dotted lines correspond as labelled to parameter values of the bifurcation diagrams in \cref{fig:tournament_bif_diag}. The blue and orange dots are the location of the codimension-\(2\) point where the fold bifurcation of the fixed points of \(f_{34}\) occurs on the switching manifolds \(\vartheta_{s}^{+}\) and \(\vartheta_{s}^{-}\), respectively.}
  \label{fig:tournament_bif_set}
\end{figure}

%% file: figure_13.tex
\begin{figure}
  \centering
  \begin{subfigure}[t]{0.49\linewidth}
    \centering
    \includegraphics[width=0.98\linewidth]{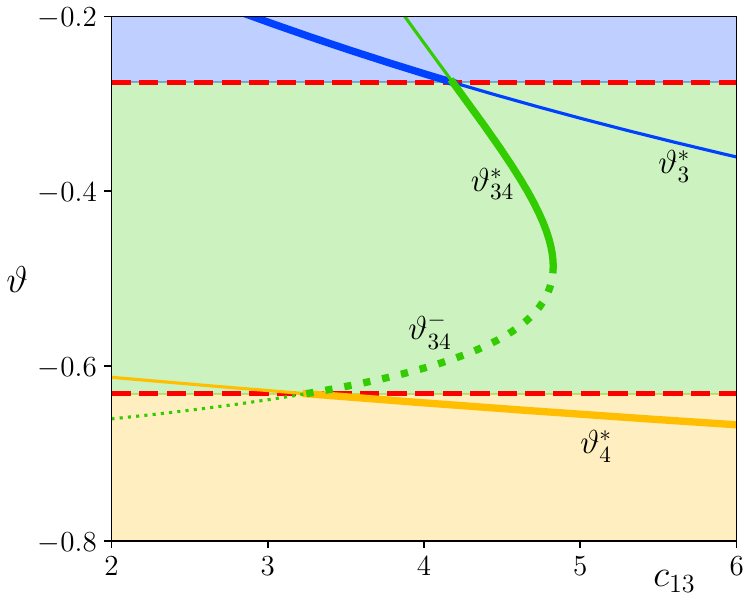}
    \vspace{-2mm}
    \caption{\(c_{14}=11\): admissible fold.}
    \label{fig:tournament_bif_diag:subfig:admiss_I}
  \end{subfigure}%
  \hfill%
  \begin{subfigure}[t]{0.49\linewidth}
    \centering
    \includegraphics[width=0.98\linewidth]{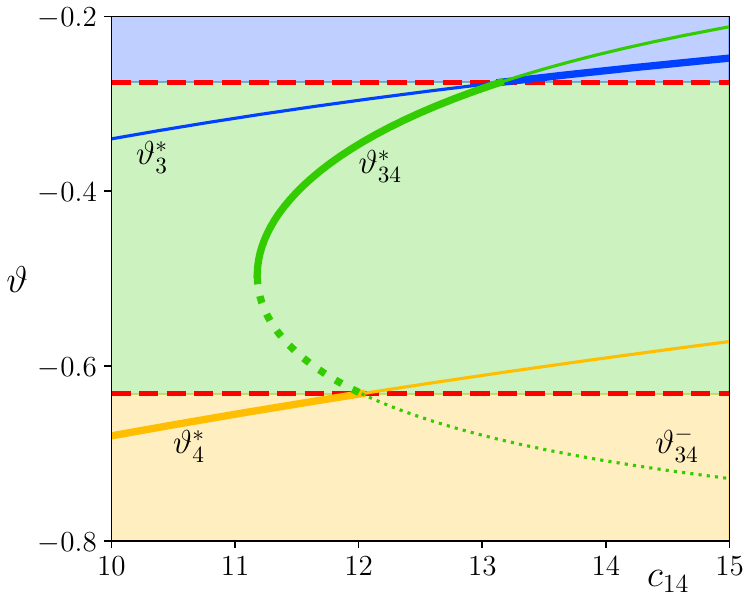}
    \vspace{-2mm}
    \caption{\(c_{13}=5\): admissible fold.}
    \label{fig:tournament_bif_diag:subfig:admiss_II}
  \end{subfigure}

  \begin{subfigure}[t]{0.49\linewidth}
    \centering
    \includegraphics[width=0.98\linewidth]{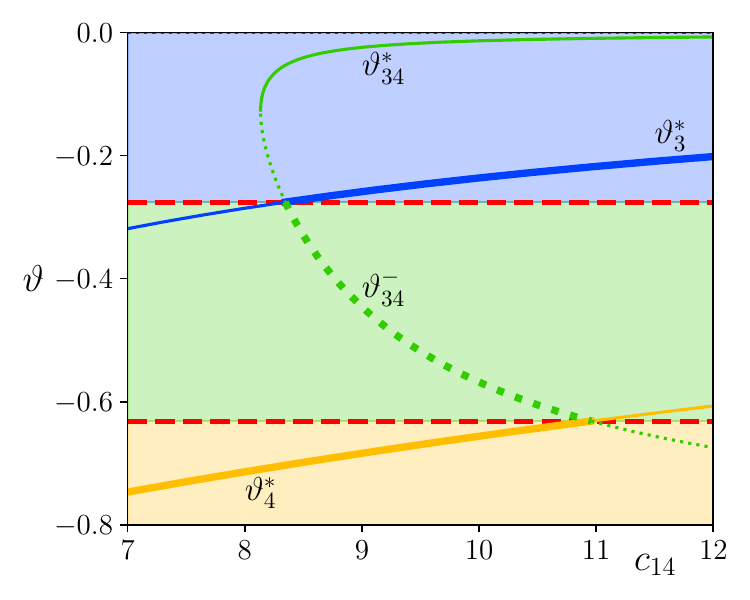}
    \vspace{-2mm}
    \caption{\(c_{13}=3.175\): virtual fold in \(\Theta_{3}\).}
    \label{fig:tournament_bif_diag:subfig:Theta3_bcb}
  \end{subfigure}%
  \hfill%
  \begin{subfigure}[t]{0.49\linewidth}
    \centering
    \includegraphics[width=0.98\linewidth]{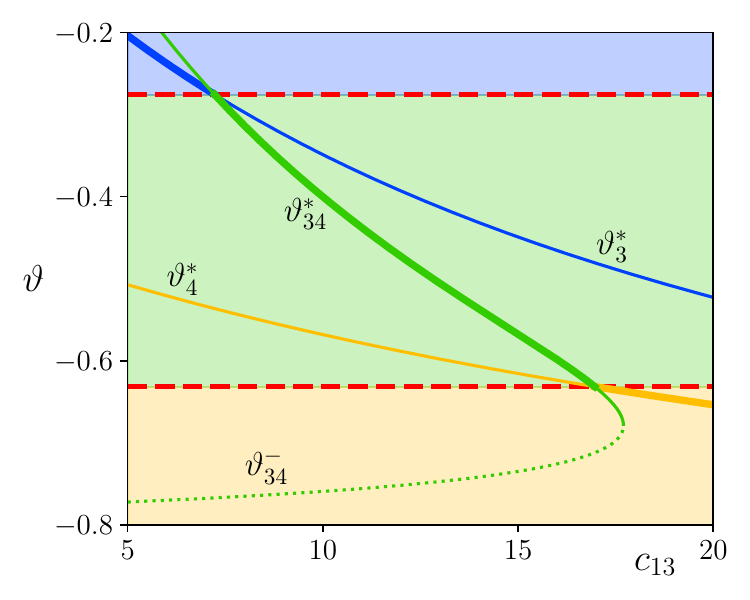}
    \vspace{-2mm}
    \caption{\(c_{14}=19\): virtual fold in \(\Theta_{4}\).}
    \label{fig:tournament_bif_diag:subfig:Theta4_bcb}
  \end{subfigure}
  \caption{Bifurcation diagrams of the fixed points of \(f_{3}\), \(f_{34}\), and \(f_{4}\), as labelled. The domain \(\Theta_{3}\), \(\Theta_{34}\), and \(\Theta_{4}\) are shaded blue, green, and orange, respectively. Stable fixed points are solid, and unstable dotted. Admissible fixed points are thick lines and virtual fixed points are thin lines. The top and bottom dashed red lines are the switching manifolds \(\vartheta_{s}^{+}\) and \(\vartheta_{s}^{-}\), respectively. All other parameters are the same as those in \Cref{fig:tournament_bif_set}.}
  \label{fig:tournament_bif_diag}
\end{figure}

%% file: figure_14.tex
\begin{figure}
  \includegraphics[width=0.46\linewidth]{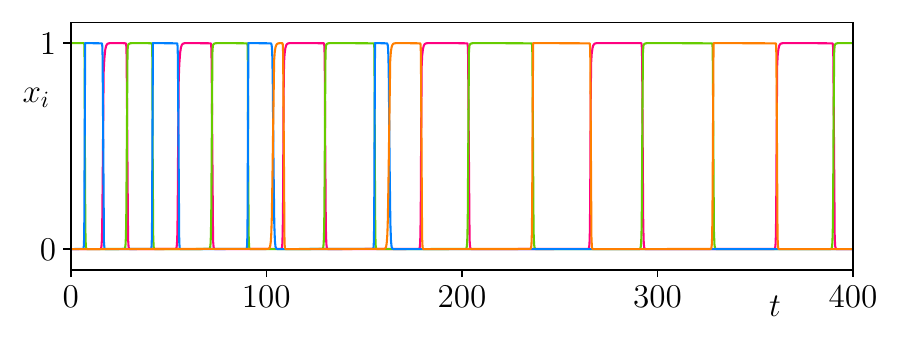}%
  \includegraphics[width=0.53\linewidth]{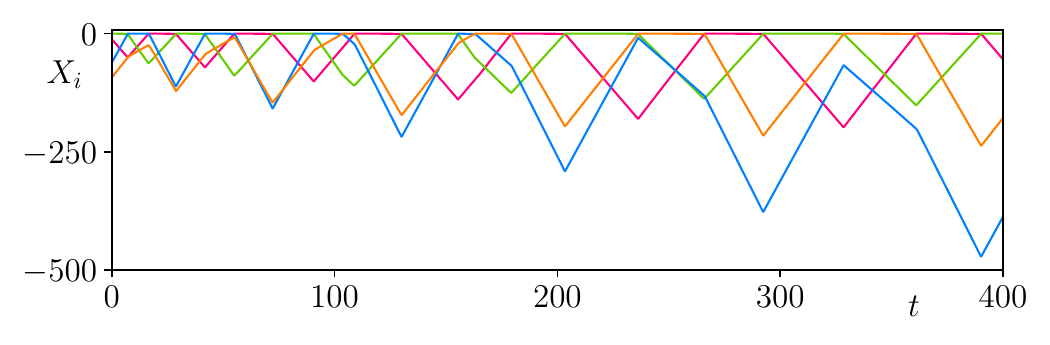}
  \caption{An example of a trajectory near the tournament network that switches from the \(\cyc{[3]}\) cycle to the \(\cyc{[34]}\) cycle, and then switches to the \(\cyc{[4]}\) cycle, to which it is asymptotic. For this example, \(\omega_{34}<0\), \(\delta_{3}>1\), \(\delta_{4}>1\), \(\mu_{3}<0\), and \(\nu_{4}>0\).}
  \label{fig:tournament_example}
\end{figure}

%% file: 5_generalised.tex
\section{Generalised projected maps and border-collision bifurcations}\label{sec:proj_map_BCB}

In this section, we generalise the analysis presented in \cref{sec:anal:ssec:kirk_silber:sssec:BCB_PodIII_equiv} to investigate border-collision bifurcations in the projected map of an arbitrary heteroclinic network

Let \(\mathcal{N}\) be a heteroclinic network. We assume this network is composed of cycles that are, in the terminology of Garrido da Silva and Castro \cite{garrido_da_silva_castro_2019}, of quasi-simple type. This condition implies that transition matrices can be defined, and that \cref{thm:pod_bif} holds. Without loss of generality, we further assume that at least two component cycles of the network have the connection \(\xi_{1}\to\xi_{2}\) in common, and that \(\xi_{2}\) is a splitting equilibrium with outgoing heteroclinic connections \(\xi_{2}\to\xi_{3}\) and \(\xi_{2}\to\xi_{4}\). We assume that the other equilibria in the network are \(\xi_{5},\dots,\xi_{n}\), and assume that all globally transverse eigenvalues are negative. Thus the eigenvalues of \(Df(\xi_{2})\) are \(-c_{21},e_{23},e_{24},-t_{25},\dots,-t_{2n}\).

The relevant basic maps are \(\varphi_{123}\colon\Gamma_{3}\subseteq\poinin[1]{2}\to\poinin[2]{3}\) and \(\varphi_{124}\colon\Gamma_{4}\subseteq\poinin[1]{2}\to\poinin[2]{4}\). As shown by Garrido da Silva and Castro \cite{garrido_da_silva_castro_2019} and Podvigina \cite{podvigina_2012}, and following the process outlined in \cref{sec:ret_proj_map:ssec:ret_maps}, the general form of the transition matrices \(m_{123}=M(\varphi_{123})\) and \(m_{124}=M(\varphi_{124})\) in the basis \(\left(X_{3},X_{4},X_{5},\dots,X_{n}\right)\) are
\begin{equation*}
  m_{123}=\begin{pmatrix}
    \frac{c_{21}}{e_{23}} & 0 & 0 & \cdots & 0\\[0.35em]
    -\frac{e_{24}}{e_{23}} & 1 & 0 & \cdots & 0 \\[0.35em]
    \frac{t_{25}}{e_{23}} & 0 & 1 & \cdots & 0 \\[0.35em]
    \vdots & \vdots & \vdots & \ddots & \vdots\\[0.35em]
    \frac{t_{2n}}{e_{23}} & 0 & 0 & \cdots & 1 \\[0.35em]
  \end{pmatrix}
  \qquad\textrm{ and }\qquad
  m_{124}=\begin{pmatrix}
     1 & -\frac{e_{23}}{e_{24}} & 0 & \cdots & 0 \\[0.35em]
     0 & \frac{c_{21}}{e_{24}} & 0 & \cdots & 0\\[0.35em]
     0 & \frac{t_{25}}{e_{24}} & 1 & \cdots & 0 \\[0.35em]
     \vdots & \vdots & \vdots & \ddots & \vdots\\[0.35em]
     0 & \frac{t_{2n}}{e_{24}} & 0 & \cdots & 1 \\[0.35em]
  \end{pmatrix}.
\end{equation*}
Note that the topology of the rest of the network does not affect these matrices.

The domains of these matrices as linear maps acting on \(\poinin[1]{2}\) in logarithmic coordinates are \(\mathcal{D}_{3}\) and \(\mathcal{D}_{4}\),  separated by the excluded subspace \(\mathcal{D}_{c}\), defined around the codimension-\(1\) switching subspace
\begin{equation*}
  W_{s}=\left\{(X_{3},X_{4},X_{5},\dots,X_{n})\in\negativeR{p}\mid e_{24}X_{3}=e_{23}X_{4}\right\}.
\end{equation*}

As we did in \cref{sec:anal:ssec:kirk_silber:sssec:BCB_PodIII_equiv,sec:anal:ssec_delta_clique:sssec:fp_exist_admiss}, we consider the image of a vector in the switching subspace under \(m_{123}\) and \(m_{124}\), which gives

\begin{equation*}
  m_{123}w_{s}=\begin{pmatrix}
    \frac{c_{21}}{e_{23}} & 0 & 0 & \cdots & 0\\[0.35em]
    -\frac{e_{24}}{e_{23}} & 1 & 0 & \cdots & 0 \\[0.35em]
    \frac{t_{25}}{e_{23}} & 0 & 1 & \cdots & 0 \\[0.35em]
    \vdots & \vdots & & \ddots & \\[0.35em]
    \frac{t_{2n}}{e_{23}} & 0 & 0 & \cdots & 1 \\[0.35em]
  \end{pmatrix}
  \begin{pmatrix}
    -1 \\[0.35em]
    -\frac{e_{24}}{e_{23}} \\[0.35em]
    X_{5} \\[0.35em]
    \vdots \\[0.35em]
    X_{n}
  \end{pmatrix}=
  \begin{pmatrix}
    -\frac{c_{21}}{e_{23}} \\[0.35em]
    0 \\[0.35em]
    -\frac{t_{25}}{e_{23}}+X_{5} \\[0.35em]
    \vdots \\[0.35em]
    -\frac{t_{2n}}{e_{23}}+X_{n}
  \end{pmatrix}.
\end{equation*}
Note that the image has a \(0\) entry, and it is straightforward to check that the image under \(m_{124}\) also has a \(0\) entry. Note also that this product holds for \textit{any} vector in the switching subspace \(W_{s}\). Therefore, if \(w_{\max}\) of \(M_{2}\) lies in \(W_{s}\), \(w_{\max}\) of \(M_{3}\) or \(M_{4}\) has a zero entry, and so does not satisfy Podvigina's third condition for f.a.s. in \cref{thm:pod_bif}. The corresponding cycle is therefore not f.a.s., and we have the following result.
\begin{theorem}\label{thm:general_BCB}
  If the vector \(w_{\max}\) of a full transition matrix \(M_{j}\) defined at a splitting equilibrium \(\xi_{j}\) lies in the switching subspace \(W_{s}\), then the corresponding heteroclinic cycle is not fragmentarily asymptotically stable.
\end{theorem}

In the projected map \(f\), the projection of \(W_{s}\) is the \(\left(n-4\right)\)-dimensional switching manifold \(\Theta_{s}\). All eigenvectors of \(M_{j}\) are fixed points of \(f\). A border-collision bifurcation occurs if the fixed point \(\vartheta^{*}\) of \(f\) that corresponds to \(w_{\max}\) of \(M_{j}\) lies in \(\Theta_{s}\). Hence, assuming that Podvigina's first two conditions for stability are satisfied, a border-collision bifurcation in the projected map corresponds to stability loss of the heteroclinic cycle, since Podvigina's third condition for stability is no longer satisfied.

A similar result holds if the splitting equilibrium has multiple incoming heteroclinic orbits, such as the bowtie network studied in \cite{castro_lohse_2016b}, or the Rock-Paper-Scissors-Lizard-Spock network studied in \cite{postlethwaite_rucklidge_2022}. However, there are more steps to constructing the projected map, which we leave for a later study of the nature of regular and irregular cycling near these and other networks. Note that, as in the case of the tournament network in \cref{sec:anal:ssec:tournament}, a cycle may have additional splitting equilibria. If so, the cycle may lose stability through a border-collision bifurcation in the switching subspace of a different splitting equilibrium. Assuming all trajectories near the network pass through the same splitting equilibrium, the pre-image of the domains of each cycle at other splitting equilibria can be taken, allowing for border-collision bifurcations to be analysed on a single cross-section. Moreover, these bifurcations, and a result analogous to \cref{thm:general_BCB}, will occur in the case of positive globally transverse eigenvalues, but these bifurcations will not result in the switching of trajectories near the network, as positive globally transverse eigenvalues only shrink the basin of attraction of the network.

%% file: 6_discussion.tex
\section{Discussion}\label{sec:disc}

In this paper, we constructed a \textit{projected map} to analyse the dynamics of trajectories near heteroclinic networks of four equilibria in \(\R^{4}\). As outlined in \cref{sec:ret_proj_map}, the projected map of a heteroclinic network is derived from the transition matrices of each cycle in the network: for the full transition matrix of each cycle defined on a splitting equilibrium, we project the linear action of the matrix onto a line \(S\) across the negative orthant. The resulting action on \(S\) is a discrete, piecewise-smooth one-dimensional map that can be more easily analysed. For the three networks we consider, all trajectories near the network pass by a common splitting equilibrium. By dividing the cross-section at the splitting equilibrium into appropriate domains, we can determine the itinerary of trajectories by studying full return maps at a single splitting equilibrium.

In the case of the Kirk--Silber network in \cref{sec:anal:ssec:ks}, our analysis reproduced the results of Kirk and Silber \cite{kirk_silber_1994}, showing that, under certain parameter values, trajectories are able to switch from one cycle to another, but only once and only in one particular direction. We also showed that trajectories may switch if one cycle is completely unstable, which occurs in a region of parameter space not considered by Kirk and Silber \cite{kirk_silber_1994}. The onset of trajectories that switch near the Kirk--Silber network occurs through a border-collision bifurcation. For the \(\Delta\)-clique network in \cref{sec:anal:ssec:delta_clique}, we considered the additional bifurcations of the length-four cycle. We showed this cycle loses stability in a bifurcation that corresponds to a fold bifurcation of the fixed points of its projected map. We demonstrated that there can be bistability of the two cycles, and that if there is bistability, there must also exist trajectories that switch from the \(\cyc{[34]}\) cycle to the \(\cyc{[4]}\) cycle. The onset of trajectories that switch is caused by the creation of the fixed points in a fold bifurcation, which produces an unstable fixed point that partitions the domain of the \(\cyc{[34]}\) cycle into trajectories that do and do not switch.

The analysis of the tournament network is more complicated, as there are many cases to consider. A full analysis of the projected map of the tournament network shows that any two component cycles can be simultaneously stable, but all three cycles cannot be simultaneously stable. If there is bistability of the \(\cyc{[3]}\) and \(\cyc{[4]}\) cycles, then the domain of the \(\cyc{[34]}\) cycle is partitioned into sets of trajectories that switch to each cycle. Again in the case of bistability of the \(\cyc{[34]}\) cycle and one of the length-three cycles, there must also exist trajectories that switch from the \(\cyc{[34]}\) cycle to the stable length-three cycle, and almost all trajectories near the unstable cycle switch to the \(\cyc{[34]}\) cycle. If there is stability of only one cycle, almost all trajectories near the network are asymptotic to that cycle, potentially cycling around one or both of the unstable cycles before switching to the stable cycle. Trajectories may in some cases switch from one length-three cycle to the other, without cycling around \(\cyc{[34]}\). Regardless of the complexity of completely classifying the dynamics of the projected map, for a given parameter values, the projected map allows us to easily determine all possible itineraries of trajectories near the network.

The analysis of trajectories using the projected map could potentially be problematic, since instability of a cycle may not be evident at the splitting equilibrium where the projected map is defined. For example, in the case of the Kirk--Silber network, the \(\cyc{[3]}\) or \(\cyc{[4]}\) cycle is unstable when \(\nu_{3}<0\) or \(\nu_{4}<0\), respectively. This bifurcation, however, occurs when \(M_{3}^{(3)}\) or \(M_{4}^{(4)}\) no longer satisfies Podvigina's third stability condition. We were able to show in \cref{prop:ks:bcb_podIII_equiv,prop:d_clique:bcb_podIII_equiv} that this stability loss is equivalent to \(w_{\max}\) of \(M_{2}^{(3)}\) or \(M_{2}^{(4)}\) lying in the switching subspace \(W_{s}\), and therefore corresponds to a border-collision bifurcation in the projected map. Thus, dynamics of trajectories and stability of cycles can be studied at \(\xi_{2}\) only, without needing to consider other equilibria. We generalised this idea in \cref{sec:proj_map_BCB}, where we showed in \cref{thm:general_BCB} that a border-collision bifurcation at a splitting equilibrium will result in a loss of stability of any quasi-simple heteroclinic cycle, including type \(Z\) cycles. This instability manifests as a breaking of Podvigina's third condition for stability at the next equilibrium in the cycle.

We have focussed on three heteroclinic networks in \(\R^{4}\) that are composed of four equilibria. The analysis shows that trajectories near a heteroclinic network of four equilibria in \(\R^{4}\) with negative globally transverse eigenvalues are always asymptotic to at most one component cycle, and can switch from one cycle to another only once. In particular, it is straightforward to show that the projected map of each network is injective and strictly increasing, and so cannot contain periodic orbits. This result confirms Brannath's conjecture \textit{``...that there should be a quite general (but non-trivial) argument...''} \cite{brannath_1994} to prove complicated switching cannot occur near these three networks.

For more complicated dynamics, we therefore need to consider networks with more equilibria, such as those studied in \cite{postlethwaite_dawes_2005,postlethwaite_rucklidge_2022,podvigina_2023}. The projected map can be defined for these and other networks, though with some additional subtleties if there are multiple splitting equilibria or multiple incoming heteroclinic orbits to the splitting equilibria. We leave the analysis of these projected maps to a future study.

%% file: 7_appendix.tex
\appendix
\section{Well-definedness of the projected map}\label{app:well_defness}

In this appendix, we show that the projected map is well-defined, and can be used to study the dynamics of trajectories near heteroclinic networks. We focus on the construction used for the Kirk--Silber network, and the same reasoning can be applied to any of the networks we consider. Recall that \(M\) is the full transition matrix for \(\xi_{2}\) defined piecewise from the matrices \(M_{3}\) and \(M_{4}\).

We first discuss why the projected map might not be well-defined. The full return map \(\Phi\colon\Gamma_{3}\cup\Gamma_{4}\to\poinin[1]{2}\) is the map from which the projected map is ultimately derived. As discussed in \cref{sec:ret_proj_map:ssec:ret_maps}, this map is not defined on the entirety of the cross-section \(\poinin[1]{2}\), due to the excluded cusp \(\Gamma_{c}\), defined around the switching curve \(\Sigma_{s}\). When moving to logarithmic coordinates, the set \(\Gamma_{c}\) becomes the subset \(\mathcal{D}_{c}\subseteq\negativeR{2}\)---defined around the switching subspace \(W_{s}\)---whose boundaries are the affine subspaces \(W_{s}^{-}\) and \(W_{s}^{+}\). These subspaces project onto the set \(S\), and so we might think that all \(\vartheta\in S\) such that \(\Pi(W_{s}^{-})<\vartheta<\Pi(W_{s}^{+})\) must be excluded from the domain of the projected map. We show here that, in fact, only \(\vartheta_{s}\), which corresponds to the switching curve \(\Sigma_{s}\), must be excluded.

Let \(\vartheta\) be a point of \(S\) that is not \(\vartheta_{s}\). This point is the projection of any vector in the linear subspace \(W(\vartheta)\) defined by \(X_{4}=-\frac{1+\vartheta}{\vartheta}X_{3}\), and this subspace is not \(W_{s}\). As \(W_{s}^{+}\) and \(W_{s}^{-}\) are parallel to \(W_{s}\), \(W(\vartheta)\) will intersect \(W_{s}^{+}\) or \(W_{s}^{-}\) at some point \(p(\vartheta)\in\negativeR{2}\), given by
\begin{equation}\label{eqn:intersection_point}
  p(\vartheta)=\begin{cases}
    p_{3}(\vartheta)=W(\vartheta)\cap W_{s}^{+}=\left(\frac{-\vartheta\log\left(1-\epsilon\right)}{\frac{e_{24}}{e_{23}}\vartheta+1+\vartheta},\frac{\left(1+\vartheta\right)\log\left(1-\epsilon\right)}{\frac{e_{24}}{e_{23}}\vartheta+1+\vartheta}\right) & \textrm{ if } \vartheta_{s}<\vartheta,\\[1em]
    p_{4}(\vartheta)=W(\vartheta)\cap W_{s}^{-}=\left(\frac{\vartheta\log\left(1-\epsilon\right)}{\vartheta+\frac{e_{23}}{e_{24}}\left(1+\vartheta\right)},\frac{-\left(1+\vartheta\right)\log\left(1-\epsilon\right)}{\vartheta+\frac{e_{23}}{e_{24}}\left(1+\vartheta\right)}\right) & \textrm{ if } \vartheta<\vartheta_{s}.
  \end{cases}
\end{equation}

\input{figure_15.tex}
If we set
\begin{equation*}
  \alpha^{*}(\vartheta)=\begin{cases}
    \alpha_{3}^{*}(\vartheta)=\frac{-\log\left(1-\epsilon\right)}{\frac{e_{24}}{e_{23}}\vartheta+1+\vartheta} & \textrm{ if } \vartheta_{s}<\vartheta,\\[0.75em]
    \alpha_{4}^{*}(\vartheta)=\frac{\log\left(1-\epsilon\right)}{\vartheta+\frac{e_{23}}{e_{24}}\left(1+\vartheta\right)} & \textrm{ if } \vartheta<\vartheta_{s},
  \end{cases}
\end{equation*}
then all vectors in the half-line \(L(\vartheta)=\left\{\alpha\left(\vartheta,-\left(1+\vartheta\right)\right)\mid\alpha>\alpha^{*}(\vartheta)\right\}\subseteq{\negativeR{2}}\) extending from, but not including, the point \(p(\vartheta)\) are in the domain of \(M\) and project onto \(\vartheta\). \Cref{fig:projection_admiss} shows examples of the process of finding the point \(p(\vartheta)\) and the half-line \(L(\vartheta)\).

Next, since \(M\) is defined by the linear action of the matrices \(M_{3}\) and \(M_{4}\), we know that, for all \(v,w\in L(\vartheta)\), \(M(v)=tM(w)\) for some \(t\). Thus, the projection of \(M(v)\) onto \(S\) is the same as that of \(M(w)\); that is, \(\Pi(M(v))=\Pi(M(w))\).

Hence, for any point \(\vartheta\in S\) such that \(\vartheta\neq\vartheta_{s}\), we define the projected map \(f\) as the projection onto \(S\) of the image under \(M\) of \(\alpha(\vartheta,-1-\vartheta)\) for any \(\alpha>\alpha^{*}(\vartheta)\):
\begin{equation}\label{eqn:proj_map_init_def}
  f\colon\vartheta\mapsto \frac{-1}{(1,1)\cdot M\left(\alpha\left(\vartheta,-\left(1+\vartheta\right)\right)\right)}e_{1}\cdot M\left(\alpha\left(\vartheta,-\left(1+\vartheta\right)\right)\right),
\end{equation}
and this mapping is well-defined. In (\ref{eqn:proj_map_init_def}), we have not cancelled the constant \(\alpha\) to emphasise that this scaling is needed so that the vector \(\alpha(\vartheta,-1-\vartheta)\) is in the domain of \(M\).

In summary, we do not have to exclude from the domain of definition of the projected map those \(\vartheta\) such that \(\Pi(W_{s}^{-})<\vartheta<\Pi(W_{s}^{+})\). We instead take the projected map \(f\) as being defined only for the projection of those vectors in the linear subspace \(W(\vartheta)\) that lie on the half-line \(L(\vartheta)\).

The following proposition allows us to use the dynamics of \(f\) to understand the dynamics of trajectories near heteroclinic networks. For this proposition, we define \(\Sigma(\vartheta)\) as the curve in \(\poinin[1]{2}\) that, in logarithmic coordinates, projects onto \(\vartheta\in S\):
\begin{equation*}
  \Sigma(\vartheta)=\left\{\left(x_{3},x_{4}\right)\in\poinin[1]{2}\mid x_{4}=x_{3}^{-\frac{1+\vartheta}{\vartheta}}\right\}.
\end{equation*}

\begin{proposition}\label{prop:valid}
  Let \(T\subseteq S\) be a subset of \(S\) for which \(\vartheta_{s}\) is not a (topological) limit point. Then for all sufficiently small neighbourhoods of the origin in \(\poinin[1]{2}\), the curve \(\Sigma(\vartheta)\) is disjoint from the excluded cusp \(\Gamma_{c}\) for all \(\vartheta\in T\).
\end{proposition}

\begin{proof}
  We begin by defining \(d^{*}\coloneq\inf\left\{d\left(t,\vartheta_{s}\right)\mid t\in T\right\}\). Since \(\vartheta_{s}\) is not a limit point of \(T\), \(d^{*}>0\). We set \(\alpha^{*}=\max\left\{\alpha^{*}(\vartheta_{s}-d),\alpha^{*}(\vartheta_{s}+d)\right\}\). Then, for all \(\vartheta\in T\), all vectors in the half-lines
  \begin{equation*}
    L(\vartheta)=\left\{\left(\alpha\vartheta,-\alpha\left(1+\vartheta\right)\right)\in\negativeR{2}\mid\alpha>\alpha^{*}\right\}\subseteq W(\vartheta)
  \end{equation*}
  are disjoint from \(\mathcal{D}_{c}\) and so are in the domain of \(M\). Therefore, all points on the curve
  \begin{equation*} 
    \Lambda(\vartheta)=\left\{\left(e^{\alpha\vartheta},e^{-\alpha\left(1+\vartheta\right)}\right)\in\poinin[1]{2}\mid\alpha>\alpha^{*}\right\}\subseteq\Sigma(\vartheta)
  \end{equation*}
  are disjoint from the excluded cusp \(\Gamma_{c}\) and so are in the domain of \(\Phi\). If
  \begin{equation*}
    r^{*}=\inf_{\vartheta\in T}\sup_{x\in\Lambda(\vartheta)}d(0,x),
  \end{equation*}
  then, for all \(r<r^{*}\), the restriction of \(\Sigma(\vartheta)\) to the open neighbourhood of radius \(r\) of the origin lies outside the excluded cusp \(\Gamma_{c}\) for all \(\vartheta\in T\).
\end{proof}

Therefore, if \(\vartheta\in\Theta_{3}\cup\Theta_{4}\) and \(\left(f^{n}(\vartheta)\right)_{n\in\N}\) is the forward orbit from \(\vartheta\), then, if this orbit is not asymptotic to \(\vartheta_{s}\), the orbit describes the dynamics of trajectories in all sufficiently small neighbourhoods of the given network.

%% file: figure_15.tex
\begin{figure}
  \centering
  \includegraphics[width=0.4\linewidth]{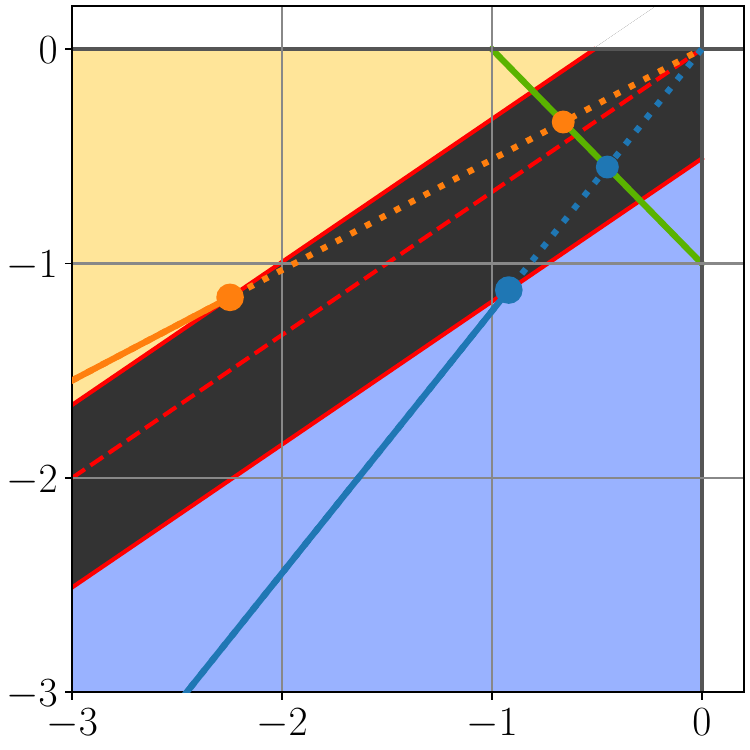}
  \caption{Two examples of points in \(S\) that are not \(\vartheta_{s}\), given by the coloured dots on \(S\), the solid green line. The linear subspace through them is shown as solid for vectors in the domain of \(M\), and dotted for vectors in \(\mathcal{D}_{c}\). The solid lines are examples of the half-lines \(L(\vartheta)\). The point \(p(\vartheta)\) (see \cref{eqn:intersection_point}) is shown as a correspondingly coloured point on the solid red lines.}
  \label{fig:projection_admiss}
\end{figure}